\documentclass{amsproc}

\usepackage[arrow, curve, frame, matrix, graph, tips]{xy}
\usepackage{amssymb}

\newtheorem{theorem}{Theorem}[section]
\newtheorem{lemma}[theorem]{Lemma}
\newtheorem{proposition}[theorem]{Proposition}
\newtheorem{corollary}[theorem]{Corollary}

\theoremstyle{definition}
\newtheorem{definition}[theorem]{Definition}
\newtheorem{example}[theorem]{Example}

\newtheorem{non-example}[theorem]{Non-Example}

\newtheorem{remark}[theorem]{Remark}

\newcommand{\op}[1]{{#1}^{\textnormal{op}}}
\newcommand{\co}[1]{{#1}^{\textnormal{co}}}
\newcommand{\rev}[1]{{#1}^{\textnormal{rev}}}
\newcommand{\coop}[1]{{#1}^{\textnormal{coop}}}
\newcommand{\iop}[1]{{#1}^{\comp}}
\newcommand{\Cat}{\textnormal{Cat}}
\newcommand{\CAT}{\textnormal{CAT}}

\newcommand{\Set}{\textnormal{Set}}
\newcommand{\PtSet}{\Set_{\bullet}}
\newcommand{\SET}{\textnormal{SET}}
\newcommand{\Span}{\textnormal{Span}}
\newcommand{\G}{\mathbb{G}}

\newcommand{\ca}[1]{\mathcal{#1}}

\newcommand{\N}{\mathbb{N}}
\newcommand{\id}{\textnormal{id}}

\newcommand{\iso}{\cong}
\newcommand{\catequiv}{\simeq}
\renewcommand{\implies}{\Rightarrow}
\renewcommand{\impliedby}{\Leftarrow}
\renewcommand{\iff}{\Leftrightarrow}
\newcommand{\comp}{\circ}

\newcommand{\ladj}{\dashv}

\newcommand{\res}{\textnormal{res}}
\newcommand{\PSh}[1]{\widehat{#1}}
\newcommand{\PSH}[1]{\CAT(\PSh {#1})}

\newcommand{\Sp}[1]{\textnormal{Sp}_{#1}}
\newcommand{\col}{\textnormal{col}}
\newcommand{\lan}{\textnormal{lan}}
\newcommand{\ran}{\textnormal{ran}}
\newcommand{\CoCts}{\textnormal{CoCts}}

\newcommand{\C}{\mathbb{C}}

\newcommand{\DFib}{\textnormal{DFib}}

\newcommand{\Cosieve}{\textnormal{Cosieve}}
\newcommand{\Fib}{\textnormal{Fib}}

\newcommand{\Sub}{\textnormal{Sub}}
\newcommand{\Sh}{\textnormal{Sh}}
\newcommand{\el}{\textnormal{el}}

\newcommand{\Hom}{\textnormal{Hom}}

\newcommand{\colim}{\textnormal{colim}}

\newcommand{\TriTwoCell}[7]{\xymatrix{
{#1} \ar[rr]^-{#5} \ar[dr]_{#4} \save \POS?="dom" \restore
&& {#2} \ar[dl]^{#6} \save \POS?="cod" \restore \\ & {#3}
\POS "dom"; "cod" **@{} ?(.35) \ar@{=>}^{#7} ?(.65)
}}

\newcommand{\TwoDiagRel}[3]{\begin{xy}
(0,0);<2em,0em>:<0em,2em>::
(-1,0)*-!R{\xybox{#1}};
(1,0)*-!L{\xybox{#3}}
**@{} ?*{#2}
\end{xy}}

\newcommand{\LaxSq}[9]{\xymatrix{{#1} \ar[r]^-{#6}
\ar[d]_{#5} \save \POS?="dom" \restore
& {#2} \ar[d]^{#7} \save \POS?="cod" \restore \\
{#3} \ar[r]_-{#8} & {#4}
\POS "dom"; "cod" **@{} ?(.35) \ar@{=>}^{#9} ?(.65)}}

\newcommand{\MNDmorphism}[6]{\xymatrix{{#1} \ar[r]^-{#5} \ar[d]_{#2}
& {#3} \save \POS="dom" \restore \ar[d]^{#4} \\
{#1} \save \POS="cod" \restore \ar[r]_-{#5} & {#3}
\POS "dom"; "cod" **@{} ?(.35) \ar@{=>}^{#6} ?(.65)}}

\newcommand{\PbSq}[8]{\xymatrix{{#1} \ar[d]_{#5}
\save \POS?(.3)="lpb" \restore
\ar[r]^-{#6} \save \POS?(.3)="tpb" \restore
& {#2} \ar[d]^{#7} \save \POS?(.3)="rpb" \restore \\
{#4} \ar[r]_-{#8} \save \POS?(.3)="bpb" \restore & {#3}
\POS "rpb"; "lpb" **@{}; ?!{"bpb";"tpb"}="cpb" **@{}; ? **@{-};
"tpb"; "cpb" **@{}; ? **@{-}}}

\newcommand{\rel}{\SelectTips{cm}{}\object@{/}}

\begin{document}

\title{Strict 2-toposes}

\author{Mark Weber}
\email{weber@math.uqam.ca}
\thanks{}


%
\bibliographystyle{amsalpha}
\maketitle
\begin{abstract}
A 2-categorical generalisation of the notion of elementary
topos is provided, and some of the properties
of the yoneda structure \cite{SW78} it generates are explored.
Examples relevant to the globular approach to higher dimensional
category theory are discussed.
This paper also contains some expository material on the
theory of fibrations internal to a finitely complete 2-category
\cite{Str74} and provides a self-contained development
of the necessary background material on yoneda structures.
\end{abstract}

\section{Introduction}

The idea of an internal category is due to Ehresmann \cite{Ehr58}.
This notion has experienced somewhat of a resurgence
in recent times because of developments in higher category theory.
For example in the work of Michael Batanin \cite{Bat98} \cite{Bat98b}
\cite{Bat02} \cite{Bat03} internal category theory is the backdrop
for the theory of higher operads.
In the work of John Baez and his collaborators
\cite{BL04} \cite{BC04} \cite{BS05}
we see internal categories as fundamental
in the process of categorifying differential geometry
and gauge theory with a view to applications in physics.

This paper is about doing category theory internally
with a particular focus on the theory of colimits.
It was motivated by the need to manipulate internal
colimits more easily in order to push forward the theory
of higher operads. In \cite{OpII} the
results and notions of the present paper are used
to bring all the operad theory of \cite{Bat98}
to the level of generality of \cite{OpI} so that
the theory of higher symmetric operads is facilitated.
This will then be used in future work on reconciling
the notions of higher dimensional categories in
\cite{BD98a} and \cite{Bat98}, and in an operadic
exploration of the stabilisation hypothesis \cite{BD95}.

In the 1970's internal colimits were understood from
a 2-categorical perspective in the work of Ross Street and Bob Walters
\cite{SW78}. In this paper the concept of a
\emph{yoneda structure} on a 2-category was discovered;
inspired largely by the work of Bill Lawvere
on the foundational importance of
the category of categories \cite{LawDct}.
Logical motivations notwithstanding,
the perspective of this paper is that
the point of having a yoneda structure
on a 2-category $\ca K$, is that one can then say
what it means for an object
of $\ca K$ to be \emph{cocomplete} in such a 
way that the theory of cocompleteness develops in $\ca K$
as in ordinary category theory.
The resulting theory of Street and Walters
clarifies colimits in both enriched and internal
category theory, and is surprisingly simple.

In this paper we focus attention on the yoneda structures
that arise in internal category theory. These
yoneda structures satisfy some additional properties
described in definition(\ref{yon-str}).
We recall the resulting theory of colimits that arises in this setting
in section(\ref{sec:yon-str}).

The examples of interest for us involve a 2-category $\ca K$
equipped with an object $\Omega \in \ca K$ which plays
the role of an internal category of sets for $\ca K$.
In the paradigmatic example we consider 
the category $\Set$ of small sets,
another category $\SET$ of sets which contains the arrows
of $\Set$ as an object, and $\CAT$ as the 2-category
of categories internal to $\SET$.
For this example $\ca K$ is $\CAT$ and $\Omega$ is $\Set$.
In higher category theory one obtains another important example
by taking $\ca K$ to be the 2-category of globular categories
and $\Omega$ to be the globular category of higher spans
of sets \cite{Bat98} \cite{Str00}.

Yoneda structures arise from the setting alluded to
in the previous paragraph because $\Omega$ satisfies
a property analogous to that enjoyed by the
subobject classifier of an elementary topos.
That this is so in the paradigmatic example described above
is an important observation of Bill Lawvere.
A 2-categorical expression of this property
is provided in the work of Ross Street by the
notion of a \emph{fibrational cosmos} \cite{Str74b} \cite{Str80}.
The \emph{2-toposes} of this work are simply
cosmoses whose underlying 2-category is cartesian closed
and comes equipped with a duality involution{\footnotemark{\footnotetext{
Less general notions of 2-toposes were considered
in \cite{Pen74} and \cite{Bou74},
in which $\Omega$ is used to classify cosieves
(see example(\ref{cdc-subobj}))
for a cartesian closed finitely complete 2-category $\ca K$.}}}.
This notion is isolated here because it is easier
to exhibit examples of 2-toposes and to explore their properties.
We exploit this to understand better some of the yoneda structures
that arise.
In particular, part of any yoneda structure
is a presheaf construction: an assignment of an object $\PSh A$
for any \emph{admissible} object $A \in \ca K$
to be regarded as the object of presheaves on $A$.
The main results of this paper describe when the yoneda structures
that arise from 2-toposes have presheaves which are cocomplete.

Our basic references for background on 2-categories is \cite{KS74}
and \cite{Str80b}.
Another important background article is \cite{Str74}
although efforts are made in this paper to keep
the exposition relatively self-contained.
The 2-categorical background pertinent to this work is
collected in section(\ref{sec:2cat-pre}), and the definition
of a duality involution is provided in subsection(\ref{subsec:2sdf-dual}).

In section(\ref{sec:2topos}) the notion of 2-topos is defined and the
basic examples are presented, and in section(\ref{2-topos->yoneda})
the yoneda structure arising from a 2-topos is described.
Not all 2-toposes give yoneda structures with cocomplete presheaves
as we see in example(\ref{psh-not-complete}).
Part of the axiomatics of a yoneda structure is a right ideal
of arrows called \emph{admissible maps}.
Section(\ref{sec:char-adm}) also provides a basic result to
help characterise the admissible maps in some of our examples.
Section(\ref{sec:psh-coco}) develops the results on presheaf
cocompleteness.
In section(\ref{sec:omcc}) we exhibit $\Omega$ as a cartesian
closed object of $\ca K$ under some hypotheses on the 2-topos
$\ca K$ (which include the cocompleteness of $\Omega$).
Applied to the globular category of higher spans,
the results of this paper say that this globular category is
the small globular colimit completion of $1$,
and that it is cartesian closed as a globular category.

In all the examples considered in this work $\ca K$
is a 2-category of categories internal to some nice
category, and so one might wonder why bother with
a 2-categorical abstraction?
One reason for this is that the theory just comes out easier
when things are expressed this way. Lax and pseudo pullbacks
are very useful things. However the main reason
is that it is 2-toposes \emph{together} with nice 2-monads
on them which provide a conceptual basis for the theory of operads
\cite{OpII}, and many of these nice 2-monads do not arise
from nice monads on the categories in which
we internalise{\footnotemark{\footnotetext{
The most basic example of this is the 2-monad on $\CAT$
whose (strict) algebras are symmetric (strict) monoidal categories.}}}.
This is especially so when one wishes to study
weakly symmetric higher dimensional monoidal categories.

Consistent with the notation for a yoneda structure,
for a category $\C$ we denote by $\PSh {\C}$ the category
of presheaves on $\C$, that is the functor category
$[\op \C, \Set]$. When working with presheaves
we adopt the standard practises of writing $C$ for
the representable $\C(-,C) \in \PSh {\C}$
and of not differentiating between
an element $x \in X(C)$ and the corresponding
map $x:C{\rightarrow}X$ in $\PSh {\C}$.
We denote by $\PSH {\C}$ the functor 2-category $[\op \C, \CAT]$
which consists of functors $\op {\C} \rightarrow \CAT$,
natural transformations between them and modifications between those.
We adopt the standard notations for the various duals
of a 2-category $\ca K$: $\op {\ca K}$ is obtained from $\ca K$
by reversing just the 1-cells, $\co {\ca K}$ is obtained by reversing
just the 2-cells, and $\coop {\ca K}$ is obtained by reversing
both the 1-cells and the 2-cells.

The title of this paper
is \emph{strict} 2-toposes because the structure
we consider on a 2-category is not invariant under biequivalence.
The reason for this is that strict 2-categorical limits
(especially pullbacks) play a central role in this theory.
So the 2-toposes of this work are different in spirit
to the higher toposes considered in homotopical algebraic geometry
\cite{TV05} which are a type of homotopy-invariant structure.

\section{2-categorical preliminaries}\label{sec:2cat-pre}

\subsection{Finitely complete $2$-categories}
Recall that a $2$-category $\ca K$ is \emph{finitely complete}
when it admits all limits weighted by $2$-functors $I:J{\rightarrow}\CAT$,
in the sense of $\CAT$-enriched category theory \cite{Kel82},
such that the set of $2$-cells of $J$ and the sets of arrows of the $I(j)$
for $j \in J$, are all finite.
This means that for such $I$ and $T:J{\rightarrow}{\ca K}$
there is an object $\lim(I,T)$ of $\ca K$ and isomorphisms of categories
\[ \ca K(X,\lim(I,T)) \iso [J,\CAT](I,\ca K(X,T)) \]
$2$-natural in $X$.
From \cite{Str76} one has the result that
$\ca K$ is finitely complete iff it admits
finite conical limits and cotensors with the ordinal $2$.
Slightly less efficiently $\ca K$ is finitely complete
iff it has a terminal object, pullbacks and lax pullbacks
(also known as comma objects).
Thus the discussion of finitely complete $2$-categories differs from
the discussion of finitely complete categories only because
for $2$-categories there are different types of pullback:
one can consider the lax pullback, pseudo pullback or pullback
of a pair of arrows $f$ and $g$ depending on whether one is considering
squares \[ \LaxSq P A B C {} {} g f {\phi} \]
in $\ca K$ such that $\phi$ is a general $2$-cell, invertible or an identity
respectively.
In addition to these variations one can also weaken the universal property
and consider bicategorical limits \cite{Str80b} but in this paper
we shall only consider the stronger notion.

Recall that the \emph{lax pullback} of
\[ \xymatrix{A \ar[r]^-{f} & B & C \ar[l]_-{g}} \]
in a $2$-category $\ca K$
consists of \[ \LaxSq P B A C p q g f {\lambda} \]
in $\ca K$ universal in the following sense:
\begin{enumerate}
\item  given \[ \LaxSq X C A B h k g f {\phi} \]
there is a unique $\delta:X{\rightarrow}{f/g}$ such that $p\delta=h$,
$q\delta=k$ and $\lambda\delta=\phi$.
\item  for $\delta_1$ and $\delta_2:X{\rightarrow}f/g$,
given $\alpha:p\delta_1{\rightarrow}p\delta_2$
and $\gamma:q\delta_1{\rightarrow}q\delta_2$ such that
\[ \xymatrix{{fp\delta_1} \ar[r]^-{\lambda\delta_1} \ar[d]_{f\alpha}
\save \POS?="domeq" \restore
& gq\delta_1 \ar[d]^{g\beta} \save \POS?="codeq" \restore \\
fp\delta_2 \ar[r]_-{\lambda\delta_2} & gq\delta_2
\POS "domeq"; "codeq" **@{}; ?*{=}} \]
there is a unique $\pi:\delta_1{\rightarrow}\delta_2$ such that
$p\pi=\alpha$ and $q\pi=\gamma$.
\end{enumerate}
It is standard notation when fixing a choice of lax pullback
to write $f/g$ in the place of $P$.
The \emph{pseudo pullback} $f/_{\iso}g$ in $\ca K$ is defined
in the same way except that the $2$-cells $\phi$ and $\lambda$
in the above definition are invertible,
and the \emph{pullback} as $f/_{=}g$ may be defined in the same way
as above except that $\phi$ and $\lambda$ are identities.
\begin{example}
Given functors \[ \xymatrix{A \ar[r]^-{f} & B & C \ar[l]_-{g}} \]
one can define the category $f/g$ as follows:
\begin{itemize}
\item  objects are triples $(a,h:fa{\rightarrow}gc,c)$ where $a \in A$,
$c \in C$, and $h \in B$.
\item  an arrow $(a,h,c){\rightarrow}(a',{h}',c')$ is a pair
$(\alpha:a{\rightarrow}a',\gamma:c{\rightarrow}c')$ such that
${h}'f(\alpha)=g(\gamma)h$.
\item  composition and identities induced
from the category structures of $A$, $B$ and $C$.
\end{itemize}
and one has
\[ \LaxSq {f/g} C A B p q g f {\lambda} \]
where $p(a,h,c)=a$, $q(a,h,c)=c$ and $\lambda_{(a,h,c)}=h$,
satisfying the universal property of a lax pullback.
One obtains $f/_{\iso}g$ as the full subcategory of $f/g$
consisting of the $(a,h,c)$ such that $h$ is invertible,
and $f/_{=}g$ as the full subcategory of $f/g$
consisting of the $(a,h,c)$ such that $h$ is an identity arrow.
Conversely given lax pullbacks in $\CAT$ one can define
general lax pullbacks as follows.
A square \[ \LaxSq P C A B p q g f {\lambda} \] in a $2$-category $\ca K$
exhibits $P$ as $f/g$ iff for all $X \in \ca K$
the functor \[ \ca K(X,P){\rightarrow}\ca K(X,f)/\ca K(X,g) \]
induced by $\ca K(X,\lambda)$ is an isomorphism of categories.
This is called the \emph{representable} definition
of lax pullbacks in $\ca K$.
Similarly, pseudo pullbacks, pullbacks and indeed all
$2$-categorical limits can be defined representably.
\end{example}
\begin{example}\label{2-slice}
Let $\ca K$ be a 2-category and $A \in \ca K$.
It is standard to denote by $\ca K/A$ the 2-category
formed as the lax pullback of
\[ \xymatrix{{\ca K} \ar[r]^-{1_{\ca K}} & {\ca K} & 1 \ar[l]_-{A}} \]
in the 2-category of 2-categories, 2-functors
and 2-natural transformations.
Its explicit description is similar to its categorical
analogue in that the objects of $\ca K/A$ are arrows $f:X{\rightarrow}A$,
and a morphism $f_1{\rightarrow}f_2$ is an arrow $g$ of $\ca K$
such that $f_2g=f_1$.
However a $2$-cell $\gamma:g_1{\implies}g_2$ of $\ca K/A$
is a $2$-cell $\gamma:g_1{\implies}g_2$ of $\ca K$
such that $f_2\gamma$ is an identity.
Recall \cite{CJ} that the 2-functor $\ca K/A{\rightarrow}\ca K$
whose object map takes the domain of an arrow into $A$,
creates any connected limits that exist in $\ca K$.
When $\ca K$ is finitely complete and $f:A{\rightarrow}B$
in $\ca K$,
the processes of taking the pullback, pseudo pullback
and lax pullback along $f$ provide $2$-functors
$\ca K/B{\rightarrow}\ca K/A$.
\end{example}
\begin{example}\label{laxpb-overC}
Let $\C$ be a category. Then a lax pullback $f/g$ of
\[ \xymatrix{A \ar[dr]_{\alpha} \save \POS?="domeq1" \restore
\ar[r]^-{f} & B \ar[d]|{\beta} \save \POS?="codeq1" \restore
\save \POS?="domeq2" \restore & C \ar[l]_-{g}
\ar[dl]^{\gamma} \save \POS?="codeq2" \restore \\ & {\C}
\POS "domeq1"; "codeq1" **@{}; ?*{=}
\POS "domeq2"; "codeq2" **@{}; ?*{=}} \]
in $\CAT/{\C}$ can be specified as follows:
the domain category is the full subcategory of $f/g$ as defined in $\CAT$
as in the previous example, consisting of the $(a,h,c)$
such that $\beta(h)$ is an identity; and the functor into $\C$
sends such $(a,h,c)$ to $\alpha(a)=\gamma(c)$ in $\C$.
Notice that while the domain $2$-functor
$\CAT/{\C}{\rightarrow}\CAT$ preserves pullbacks,
from the description of lax pullbacks in $\CAT/{\C}$ given here,
it does \emph{not} preserve lax pullbacks in general.
\end{example}
\begin{example}\label{internal-laxpb}
Let $\ca E$ be a category with pullbacks.
Recall the 2-category $\Cat(\ca E)$ of categories
internal to $\ca E$.
For $A \in \Cat(\ca E)$ it is standard to denote by
\[ \xymatrix{{A_0} \ar[r]|{i}
& {A_1} \ar@<-1ex>[l]_{s} \ar@<1ex>[l]^{t} \ar@<1ex>[r] \ar@<-1ex>[r]
& {A_2} \ar@<-2ex>[l] \ar[l]|{c} \ar@<2ex>[l]} \]
the 2-truncated simplicial diagram determined by $A$: $A_0$
is the object of objects of $A$, $A_1$ the object of arrows,
$s$ the source or domain map, $t$ the target or codomain map,
$i$ the map that provides identity arrows, $A_2$ the object
of composable pairs in $A$ obtained by pulling back $s$ along $t$,
and $c$ the composition map.
Using pullbacks alone
one can construct pullbacks, pseudo pullbacks and lax pullbacks
in $\Cat(\ca E)$. One can either do this directly, or by
interpretting the explicit description of these constructions
in $\CAT$ in the internal language of $\ca E$ (see \cite{JohElpI}).
A nice consequence of this is that for any pullback preserving
functor ${\ca E}{\rightarrow}{\ca E}'$, the 2-functor
$\Cat(\ca E){\rightarrow}\Cat({\ca E}')$ it induces preserves
pullbacks, pseudo pullbacks and lax pullbacks.
\end{example}
Lax pullbacks satisfy the same composition and cancellation
properties as pullbacks do in ordinary category theory.
That is, given
\[ \xymatrix{X \ar[r] \ar[d] \save \POS?="domeq" \restore & A \ar[r]^-{v}
\ar[d]|{u} \save \POS?="codeq" \restore \save \POS?="dom" \restore
& B \ar[d]^{g} \save \POS?="cod" \restore \\
Y \ar[r]_-{h} & C \ar[r]_-{f} & D
\POS "domeq"; "codeq" **@{} ?*{=}
\POS "dom"; "cod" **@{} ?(.35) \ar@{=>}^{\lambda} ?(.65)} \]
such that $\lambda$ exhibits $A$ as $f/g$,
then the composite square exhibits $X$ as $fh/g$ iff the front square is a
pullback.
Similarly for pseudo pullbacks and pullbacks in any $2$-category.
So for example, one can use this observation to
obtain all lax pullbacks
from pullbacks and lax pullbacks of identity arrows.
We adopt some standard notational abuses with regards to identity arrows:
$f/1_A$ is written $f/A$, $1_A/g$ is written $A/g$ and $1_A/1_A$ is written
$A/A$.

Some other notation: we shall denote a terminal object of
$\ca K$ by $1$ and for $A \in \ca K$ the unique
map $A{\rightarrow}1$ by $t_A$. 

\subsection{Fibrations}\label{fib}
The concept of a fibration between categories
is due to Grothendieck \cite{Gro70}.
Fibrations can be defined internal to
any finitely complete $2$-category $\ca K$.
In fact there are two approaches:
one can work $2$-categorically and follow \cite{Str74},
or one can regard $\ca K$ as a bicategory following \cite{Str80b}.
For the case $\ca K=\CAT$ the $2$-categorical definition
of fibration coincides with that of Grothendieck \cite{Gray66}
whereas the bicategorical definition of fibration
is more general.
In this paper we shall consider only the stronger notion.

Let $f:A{\rightarrow}B$ be a functor.
A morphism $\alpha:a_1{\rightarrow}a_2$ in $A$ is \emph{$f$-cartesian}
when for all $\alpha_1$ and $\beta$ as shown:
\[ \xymatrix{{fa_1} \ar[r]^-{f\alpha} & {fa_2} \\
{fa_3} \ar[u]^{\beta} \save \POS?(.65)="domeq" \restore
\ar[ur]_{f\alpha_1} \save \POS?(.65)="codeq" \restore
\POS "domeq"; "codeq" **@{}; ?*{=}} \]
there is a unique $\gamma:a_3{\rightarrow}a_1$ such that
$f\gamma=\beta$ and $\alpha\gamma=\alpha_1$.
The basic facts about $f$-cartesian morphisms are:
\begin{enumerate}
\item  if
$\xymatrix@1{{a_1} \ar[r]^-{\alpha_1} & {a_2} \ar[r]^-{\alpha_2} & {a_3}}$
are in $A$ and $\alpha_2$ is $f$-cartesian, then
$\alpha_1$ is $f$-cartesian iff $\alpha_2\alpha_1$ is $f$-cartesian.
\item  isomorphisms are $f$-cartesian for any $f$.
\item  if $\alpha$ is $f$-cartesian and $f\alpha$ is an isomorphism
then $\alpha$ is an isomorphism.
\end{enumerate}
A \emph{cartesian lift} of the pair $(\beta:b{\rightarrow}fa, a)$
is an $f$-cartesian morphism $\alpha:a_2{\rightarrow}a$
such that $f\alpha=\beta$.
Grothendieck defined $f$ to be a \emph{fibration}
when every $(\beta:b{\rightarrow}fa, a)$ has a cartesian lift.

Let $\ca K$ be a finitely complete $2$-category and
$f:A{\rightarrow}B$ be in $\ca K$.
A $2$-cell
\[ \xymatrix{X \ar@/^{1pc}/[rr]^-{a_1} \save \POS?="domal" \restore
\ar@/_{1pc}/[rr]_-{a_2} \save \POS?="codal" \restore && A
\POS "domal"; "codal" **@{}; ?(.25) \ar@{=>}^{\alpha} ?(.75)} \]
is \emph{$f$-cartesian} when for all $g:Y{\rightarrow}X$ in $\ca K$,
${\alpha}g \in \ca K(Y,A)$ is $\ca K(Y,f)$-cartesian.
One then defines $f$ to be a \emph{fibration} when for all
\[ \TriTwoCell X A B b a f {\beta} \]
there exists $f$-cartesian $\overline{\beta}:c{\implies}a$ so that
$fc=b$ and $f\overline{\beta}=\beta$.
The following well known result is fundamental and easy to prove.

\begin{theorem}\label{fib-basic-prop}
In any finitely complete 2-category $\ca K$:
\begin{enumerate}
\item  The composite of fibrations is a fibration.
\item  The pullback of a fibration along any map is a fibration.
\end{enumerate}
\end{theorem}

It is natural to consider applying the cartesian lifting criterion
of a fibration $f$ to the $2$-cell from a lax pullback involving $f$.
When one does this, one is lead to theorem(\ref{Chev-crit}).
Given $f:A{\rightarrow}B$ and $g:X{\rightarrow}B$
in a finitely complete $2$-category $\ca K$,
we shall denote by $i:g/_=f{\rightarrow}g/f$ the map induced
by the universal property of $g/f$ and the identity cell
\[ \LaxSq {g/_=f} A X B {} {} f g {\id} \]
of a defining pullback square for $g/_=f$.
\begin{lemma}\label{adj-eq}
Let $\ca K$ be a $2$-category and
$i:X{\rightarrow}Y$ in $\ca K$.
Then to give $i \ladj r$ with invertible unit,
it suffices to give $\varepsilon:ir{\implies}1_Y$
such that $r\varepsilon$ and ${\varepsilon}i$ are invertible
and $1_X{\iso}ri$.
\end{lemma}
\begin{theorem}\label{Chev-crit}
\cite{Str74}
Let $\ca K$ be a finitely complete $2$-category
and $f:A{\rightarrow}B$ in $\ca K$.
Then the following statements are equivalent:
\begin{enumerate}
\item  $f$ is a fibration.
\label{Chev-fib}
\item  For all $g:X{\rightarrow}B$, the map $i:g/_=f{\rightarrow}g/f$
has a right adjoint in $\ca K/X$ with invertible unit.
\label{Chev-genlaxpb}
\item  The map $i:A{\rightarrow}B/f$
has a right adjoint in $\ca K/B$ with invertible unit.
\label{Chev-speclaxpb}
\end{enumerate}
\end{theorem}
\begin{proof}
(\ref{Chev-fib})$\implies$(\ref{Chev-genlaxpb}):
Apply the above definition
to $\phi=\lambda$ the defining lax pullback for $g/f$, to obtain
\[ \TwoDiagRel
{\LaxSq {g/f} A X B {p_1} {q_1} f g {\lambda}}
{=}
{\xymatrix{{g/f} \ar[dd]_{p_1} \save \POS?="domeql" \restore
\ar[dr]_{r} \ar[rr]^-{q_1} \save \POS?="codlamb" \restore
&& A \ar[dd]^{f} \save \POS?(.75)="codeqr" \restore
\\ & {g/_=f} \save \POS="domlamb" \restore \save \POS="codeql" \restore \ar[dl]^{p_2}
\save \POS?="domeqr" \restore \ar[ur]_{q_2} \\ X \ar[rr]_-{g} && B
\POS "domlamb"; "codlamb" **@{}; ?(.5) \ar@{=>}^{\overline{\lambda}} ?(.75)
\POS "domeql"; "codeql" **@{}; ?*{=} \POS "domeqr"; "codeqr" **@{}; ?(.6)*{=}}} \]
where $\overline\lambda$ is $f$-cartesian, and so by the universal property
of $g/f$ we obtain $\varepsilon:ir{\implies}1$ such that $p_1\varepsilon=\id$
and $q_1\varepsilon=\overline{\lambda}$.
By lemma(\ref{adj-eq}) it suffices to show that
${\varepsilon}i$ and $r\varepsilon$ are invertible,
and that $ri{\iso}1$ in $\ca K/X$.
Observe that $q{\varepsilon}i=\overline{\lambda}i$ is an $f$-cartesian lift
of an identity cell, and thus is invertible.
Since $p{\varepsilon}i=\id$, the $2$-cell ${\varepsilon}i$ is invertible
by the universal property of $g/f$.
So we have $\overline{\lambda}i:q_2ri{\iso}q_2$ and $p_2ri=p_2$,
and by the above defining diagram/equation of $\overline{\lambda}$
precomposed with $i$, together with the universal property of $g/_=f$,
one obtains $\eta:ri{\iso}1$ such that $p\eta=\id$.
Finally to see that $r\varepsilon$ is invertible,
by the universal property of $g/_=f$, it suffices to show that
$q_2r\varepsilon$ is invertible, since $p_2r\varepsilon=p_1\varepsilon=\id$.
Note that $fq_2r\varepsilon=gp_2r\varepsilon=gp_1\varepsilon=\id$,
and so it suffices to show that $q_2r\varepsilon$ is $f$-cartesian.
For this we note that
\[ \xymatrix{{q_2rir} \ar[d]_{\overline{\lambda}i_r} \ar[r]^-{q_2r\varepsilon}
& {q_2r} \ar[d]^{\overline{\lambda}} \\ {q_1ir} \ar[r]_-{q_1\varepsilon} & {q_1}} \]
commutes, and that since $\overline{\lambda}=q_1\varepsilon$ is $f$-cartesian,
$q_2r\varepsilon$ is $f$-cartesian also.\\
(\ref{Chev-genlaxpb})$\implies$(\ref{Chev-speclaxpb}):
just take $g=1_B$.\\
(\ref{Chev-speclaxpb})$\implies$(\ref{Chev-fib}):
By the universal property of $B/f$ it suffices to verify the defining property of a fibration
as defined above for the case $\beta=\lambda$, the defining lax pullback cell for $B/f$.
We have
\[ \xymatrix{& {B/f} \ar[dl]_{r} \ar[d]^{1} \save \POS?="codeps" \restore \\
A \save \POS="domeps" \restore \ar[r]_-{i} & {B/f} \ar[d]_{p} \save \POS?="domlam" \restore
\ar[r]^-{q} & A \ar[d]^{f} \save \POS?="codlam" \restore \\ & B \ar[r]_-{1} & B
\POS "domeps"; "codeps" **@{}; ?(.55) \ar@{=>}^{\varepsilon} ?(.75)
\POS "domlam"; "codlam" **@{}; ?(.35) \ar@{=>}^{\lambda} ?(.65)} \]
where $\varepsilon$ is the counit for the adjunction $i \ladj r$ in $\ca K/B$.
Since ${\lambda}i=\id$ this composite evaluates to $fq\varepsilon$,
and since $p\varepsilon=\id$ this composite also evaluates to $\lambda$
whence $\lambda=fq\varepsilon$.
Thus it suffices to show that $q\varepsilon$ is $f$-cartesian.
To this end let $\delta:s{\rightarrow}q$ and $\gamma:fs{\rightarrow}fr$
such that $f\delta=f(q\varepsilon)\gamma$.
Another way to express this last equation, since $\lambda=fq\varepsilon$
and $\lambda{is}=\id$, is that
\[ \xymatrix{{pis} \ar[d]_{\lambda{is}} \ar[r]^-{\gamma} & p \ar[d]^{\lambda} \\
{fqis} \ar[r]_-{f\delta} & {fq}} \]
commutes, and so by the universal property of $B/f$, there is a unique $\beta:is{\rightarrow}1$
such that $p\beta=\gamma$ and $q\beta=\delta$.
Since $\varepsilon$ exhibits $r$ as a right lifting of $1$ along $i$
(see example(\ref{adj-univ}) below), there is a unique $\alpha$ such that
\[ \TwoDiagRel
{\xymatrix{{B/f} \ar[dr]_{1} \save \POS?="codbet" \restore \ar[rr]^-{s}
&& A \ar[dl]^{i} \save \POS?="dombet" \restore \\ & {B/f}
\POS "dombet"; "codbet" **@{}; ?(.35) \ar@{=>}_{\beta} ?(.65)}}
{=}
{\xymatrix{{B/f} \ar[dr]_{1} \save \POS?="codeps" \restore
\ar[rr]|-{r} \save \POS?="codal" \restore
\ar@/^{2pc}/[rr]^-{s} \save \POS?="domal" \restore
&& A \ar[dl]^{i} \save \POS?="domeps" \restore \\ & {B/f}
\POS "domal"; "codal" **@{}; ?(.3) \ar@{=>}_{\alpha} ?(.7)
\POS "domeps"; "codeps" **@{}; ?(.35) \ar@{=>}_{\varepsilon} ?(.65)}} \]
Post composing this last equation with $p$ gives $\gamma=f\alpha$,
and post composing it with $q$ gives $\delta=(q\varepsilon)\alpha$.
Conversely $\alpha$ is clearly the unique $2$-cell satisfying these equations.
\end{proof}
A functor $p:E{\rightarrow}B$ is an \emph{opfibration}
when the functor $\op p:{\op E}{\rightarrow}{\op B}$ is a fibration.
These functors were originally called cofibrations by Grothendieck,
and indeed they are fibrations in $\co {\CAT}$,
however beginning with Gray \cite{Gray66} the term opfibration was used instead
so as not to give topologists the wrong idea: it is the fibrations
in $\op {\CAT}$ which are more like what a topologist would call
a cofibration, and both fibrations and opfibrations with their
lifting properties which correspond intuitively to topologists' fibrations.
Similarly one refers to \emph{opcartesian liftings},
defines opfibrations in an arbitrary 2-category $\ca K$ as
fibrations in $\co {\ca K}$, and interprets the above results
in $\co {\ca K}$ when working with this dual notion.
By the characterisation of fibrations given in
theorem(\ref{Chev-crit})(\ref{Chev-speclaxpb}) and its dual
one has the following immediate corollary.
\begin{corollary}\label{cor:pres-fib}
If a 2-functor between finitely complete 2-categories preserves
lax pullbacks,
then it preserves fibrations and opfibrations.
\end{corollary}
\begin{example}\label{fib-cat}
In this example we revisit fibrations in $\CAT$ from the present general
point of view and recall the Grothendieck construction for later use.
To see that fibrations in $\CAT$ are indeed Grothendieck fibrations
let $f:A{\rightarrow}B$ be a functor.
To give a right adjoint to $i:A{\rightarrow}B/f$ over $B$,
one must give for each pair $(\beta:b{\rightarrow}fa,a)$
an object $r(\beta,a)$ in $A$ such that $fr(\beta,a)=b$,
and an arrow $\varepsilon_{(\beta,a)}:r(\beta,a){\rightarrow}a$
such that $f\varepsilon_{\beta,a}=\beta$ and so that the map
\[ (\id,\varepsilon_{(\beta,a)}) : (\id,r(\beta,a)) \rightarrow
(\beta, a) \]
in $B/f$ has the universal property of the counit of an adjunction.
This universal property amounts to $\varepsilon_{(\beta,a)}$
being $f$-cartesian.
Thus a right adjoint to $i:A{\rightarrow}B/f$ over $B$
amounts to a choice of cartesian liftings for $f$.
A \emph{cleavage} is the standard terminology for
such a choice of cartesian liftings,
and a fibration together with a cleavage
is known as a \emph{cloven} fibration.
Given a pseudo functor $X:\op \C{\rightarrow}\CAT$
the Grothendieck construction
produces an associated cloven fibration over $\C$.
Define the category $\el(X)$ as follows:
\begin{itemize}
\item  objects: are pairs $(x,C)$ where $C \in \C$ and $x \in X(C)$.
\item  an arrow $(x_1,C_1){\rightarrow}(x_2,C_2)$
is a pair $(\alpha,\beta)$ where $\beta:C_1{\rightarrow}C_2$ in $\C$ and
$\alpha:x_1{\rightarrow}X\beta(x_2)$.
\item  compositions: formed in the evident way using the pseudo functor
coherence cells and composition in $\C$.
\end{itemize}
define the functor into $\C$ to be the obvious forgetful functor,
and note that this functor has an obvious cleavage.
This construction provides a $2$-equivalence
$\Fib(\C)\catequiv{\Hom(\op \C,\CAT)}${\footnotemark{\footnotetext{
$\Hom(\op \C,\CAT)$ is the $2$-category
of pseudo functors $\op \C{\rightarrow}\CAT$,
pseudo natural tranformations between them,
and modifications between those
and $\Fib(\C)$ is
the sub-2-category of $\CAT/{\C}$ consisting of fibrations
and functors over $\C$ that preserve cartesian arrows.}}}.
A useful perspective on this last fact is that the $2$-functor
\[ \el : \PSH {\C} \rightarrow \CAT \]
factors as
\[ \xymatrix{{\PSH {\C}} \ar[r] & {{\CAT}/{\C}} \ar[r] & {\CAT}} \]
where the first 2-functor is monadic
and the second takes the domain of a functor into $\C$.
The induced 2-monad on on ${\CAT}/{\C}$ has functor
part given by taking lax pullbacks along $1_{\C}$
and the 2-category of pseudo algebras for this 2-monad
is 2-equivalent to $\Fib(\C)$.
The general idea of seeing fibrations as algebras of a 2-monad
in this way was developed in \cite{Str74}.
Note also that this factorisation of $\el$
implies that $\el$ preserves all connected conical limits \cite{CJ}.
\end{example}
\begin{example}\label{fib-overC}
For a category $\C$ we shall now consider fibrations in $\CAT/{\C}$.
In example(\ref{laxpb-overC}) we described
explicitly the construction of lax pullbacks in $\CAT/{\C}$.
Thus one can unpack the definition of fibration in $\CAT/{\C}$
in much the same way as we did for $\CAT$ in the previous example.
When one does this one finds that
\[ \xymatrix{A \ar[dr]_{\alpha} \save \POS?="domeq" \restore \ar[rr]^-{f}
&& B \ar[dl]^{\beta} \save \POS?="codeq" \restore \\ & {\C}
\POS "domeq"; "codeq" **@{}; ?*{=}} \]
is a fibration in $\CAT/{\C}$
iff every $(\phi:b{\rightarrow}fa,a)$ such that $\beta(\phi)=\id$
has a cartesian lift.
In this way the idea that $f$ be ``locally'' a fibration,
or in other words a
``fibration on the fibres of $\beta$'',
is formalised by saying that $f$ is a fibration in $\CAT/{\C}$.
In particular notice that the condition that $f$ be a fibration over $\C$
is in general weaker than the condition that $f$ be a fibration,
although they are equivalent when $f$ lives in $\Fib(\C)$
\cite{Ben85} \cite{Her99} as we shall now recall.
For such an $f$ let $\phi:b{\rightarrow}fa$.
To obtain a cartesian lift for $\phi$
one first takes a cartesian lift $\phi_1:a_1{\rightarrow}a$ of
$\beta(\phi)$. Since $f$ preserves cartesian arrows $\phi$
factors uniquely as
\[ \xymatrix{b \ar[r]^-{\phi_2} & {fa_1} \ar[r]^-{f\phi_1} & {fa}} \]
where $\beta(\phi_2)=\id$. Since $f$ is a fibration in $\CAT/{\C}$
one can take a cartesian lift $\phi_3$ of $\phi_2$,
and then the composite
\[ \xymatrix{{a_2} \ar[r]^-{\phi_3} & {a_1} \ar[r]^-{\phi_1} & a} \]
is a cartesian lift for $\phi$.
\end{example}
%

\subsection{2-sided discrete fibrations and duality involutions}
\label{subsec:2sdf-dual}
In this subsection $\ca K$ is a finitely complete 2-category.
We shall now recall the notion of 2-sided discrete fibration
of \cite{Str74} and an associated yoneda lemma,
and then define a notion of duality involution for $\ca K$.
Recall first the bicategory $\Span(\ca K)$:
\begin{itemize}
\item  objects: those of $\ca K$.
\item  a morphism $A{\rightarrow}B$ in $\Span(\ca K)$ is a triple $(d,E,c)$
\[ \xymatrix{A & E \ar[l]_-{d} \ar[r]^-{c} & B} \]
and is called a \emph{span} from $A$ to $B$.
When the maps $d$ and $c$ are understood
we may abuse notation a little and refer to the above span as $E$.
\item  a $2$-cell between spans is a commutative diagram
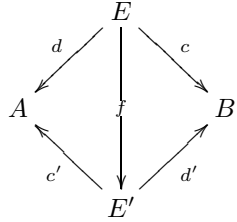
\[ \xymatrix{& E \ar[dl]_{d} \ar[dr]^{c} \ar[dd]|{f} \\ A && B \\
& {E'} \ar[ul]^{c'} \ar[ur]_{d'}} \]
Composition of $2$-cells is just composition in $\ca K$.
\item  composition of $1$-cells is obtained by pulling back:
the composite of $(d_1,E_1,c_1)$ and $(d_2,E_2,c_2)$
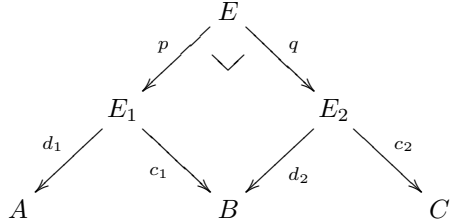
\[ \xymatrix{&& E \ar[dl]_{p} \save \POS?(.3)="lpb" \restore
\ar[dr]^{q} \save \POS?(.3)="tpb" \restore \\
& E_1 \ar[dl]_{d_1} \ar[dr]_{c_1} \save \POS?(.3)="bpb" \restore
&& E_2 \ar[dl]^{d_2} \save \POS?(.3)="rpb" \restore \ar[dr]^{c_2}
\\ A && B && C
\POS "rpb"; "lpb" **@{}; ?!{"bpb";"tpb"}="cpb" **@{}; ? **@{-};
"tpb"; "cpb" **@{}; ? **@{-}} \]
is $(d_1p,E,c_2q)$, and this composite may be denoted as
$E_2 \comp E_1$.
\end{itemize}
Notice that the homs $\Span(\ca K)(A,B)$
being isomorphic to the underlying category of
${\ca K}/{(A \times B)}$ are in fact 2-categories,
and that pullback-composition is 2-functorial.
Given a map $f:A{\rightarrow}B$ in $\ca K$,
we have spans $(1_A,A,f)$ and $(f,A,1_A)$
which we denote by $f$ and $\rev f$ respectively.
Note also that $f \ladj \rev f$ in $\Span(\ca K)$
and that for any span $(d,E,c)$
one has $E \iso \rev c \comp d$.

A span $(d,E,c)$ from $A$ to $B$
in $\CAT$ is a \emph{discrete fibration from $A$ to $B$}
when it satisfies:
\begin{enumerate}
\item  $\forall f:a{\rightarrow}d(e)$, ${\exists}! \overline{f}$
in $E$ with codomain $e$ such that $d(\overline{f})=f$ and
$c(\overline{f})=1_{c(e)}$.
\item  $\forall g:c(e){\rightarrow}b$, ${\exists}! \overline{g}$
in $E$ with domain $e$ such that $d(\overline{f})=1_d(e)$ and
$c(\overline{f})=f$.
\item  $\forall h:e_1{\rightarrow}e_2$ in $E$, the composite
$\overline{d(h)} \comp \overline{c(h)}$ is defined and equal to $h$.
\end{enumerate}
More generally, a span $(d,E,c)$ from $A$ to $B$
in $\ca K$ is a \emph{discrete fibration from $A$ to $B$}
when for all $X \in \ca K$
\[ \xymatrix{{\ca K(X,A)} & {\ca K(X,E)} \ar[l]_-{\ca K(X,d)}
\ar[r]^-{\ca K(X,c)} & {\ca K(X,B)}} \]
is a discrete fibration from $\ca K(X,A)$ to $\ca K(X,B)$.
We define a category $\DFib(\ca K)(A,B)$ whose objects are
discrete fibrations from $A$ to $B$ in $\ca K$,
and morphisms are morphisms
of the underlying spans{\footnotemark{\footnotetext{
A 2-cell between maps of discrete fibrations is necessarily
an identity}}}.
In particular a map $p:E{\rightarrow}B$ is a \emph{discrete fibration}
when the span
\[ \xymatrix{B & E \ar[l]_-{p} \ar[r] & 1} \]
is a discrete fibration,
and a \emph{discrete opfibration} when the span
\[ \xymatrix{1 & E \ar[l] \ar[r]^-{p} & B} \]
is a discrete fibration.
\begin{theorem}\label{thm-2-sided-fib}
In a finitely complete 2-category $\ca K$:
\begin{enumerate}
\item  Given a discrete fibration $E$ from $A$ to $B$ and a map
$f:C{\rightarrow}A$, the span $E \comp f$ is a discrete fibration
from $C$ to $A$.
\label{2sdf-I}
\item  Given a discrete fibration $E$ from $A$ to $B$ and a map
$g:D{\rightarrow}A$, the span $\rev g \comp E$ is a discrete fibration
from $A$ to $D$.
\label{2sdf-II}
\item  If $(d,E,c)$ is a discrete fibration from $A$ to $B$
then $d$ is a fibration and $c$ is an opfibration.
\label{2sdf-III}
\item  If $f:A{\rightarrow}C$ and $g:B{\rightarrow}C$,
then the span from $A$ to $B$ obtained from the lax pullback $f/g$
is a discrete fibration from $A$ to $B$.
\label{2sdf-IV}
\end{enumerate}
\end{theorem}
\begin{proof}
Because of the representability of the notions
involved, to prove (\ref{2sdf-I}), (\ref{2sdf-II}) and (\ref{2sdf-IV})
it suffices to verify each statement in the case $\ca K=\CAT$,
and in this case the direct verifications are straight forward.
As for (\ref{2sdf-III})
it suffices to show that $d$ is a fibration, because
this result interpretted in $\co {\ca K}$ says that $c$ is an opfibration.
Consider
\[ \xymatrix{&& E \ar[dl]^{d} \save \POS?="codphi" \restore \ar[dr]^{c} \\
X \ar[r]_-{a}
\ar@/^{2pc}/[urr]^-{e} \save \POS?(.3)="domphi" \restore  & A && B
\POS "domphi"; "codphi" **@{}; ?(.35) \ar@{=>}^{\phi} ?(.65)} \]
and take $\varepsilon:e_1{\rightarrow}e$ to be the unique 2-cell
such that $d\varepsilon=\phi$ and $c\varepsilon=\id$.
We will now show that $\varepsilon$ is $d$-cartesian.
Let $z:Z{\rightarrow}X$, $g:e_2{\rightarrow}ez$ and
$f:pe_2{\rightarrow}az$ be such that $pg=(\phi{z})f$.
We must produce $\delta:e_2{\rightarrow}e_1z$
unique such that $d\delta=f$ and $(\varepsilon{z})\delta=g$.
Since $(d,E,c)$ is a discrete fibration we can factor $g$
uniquely as
\[ \xymatrix{{e_2} \ar[r]^-{g_1} & {e_3} \ar[r]^-{g_2} & {ez}} \]
where
\[ \begin{array}{cccc}
{dg_1=\id} & {cg_1=cg} & {dg_2=dg} & {cg_2=\id} \end{array} \]
and we can take $h:e_4{\rightarrow}e_1z$ unique such that
$dh=f$ and $ch=\id$.
Since $d((\varepsilon{z})h)=(\phi{z})f=dg$
and $c((\varepsilon{z})h)=\id$,
we have $g_2=(\varepsilon{z})h$ and so $e_3=e_4$.
Thus we take $\delta$ to be the composite
\[ \xymatrix{{e_2} \ar[r]^-{g_1} & {e_3} \ar[r]^-{h} & {e_1z}} \]
because $d\delta=(dh)(dg_1)=f$ and
$(\varepsilon{z})\delta=(\varepsilon{z})hg_1=g_2g_1=g$.
To see that $\delta$ is unique, let $\delta':e_2{\rightarrow}e_1z$
be such that $d\delta'=f$ and $(\varepsilon{z})\delta'=g$.
Factor $\delta'$ uniquely as
\[ \xymatrix{{e_2} \ar[r]^-{\delta_1} & {e_5} \ar[r]^-{\delta_2}
& {e_1z}} \]
where
\[ \begin{array}{cccc}
{d\delta_1=\id} & {c\delta_1=c\delta'} & {d\delta_2=f} & {c\delta_2=\id}
\end{array} \]
Then $\delta_1$ is forced to be $g_1$ since
$d\delta_1=\id=dg_1$ and $c\delta_1=c\delta'=c((\varepsilon{z})\delta')
=cg=cg_1$,
and $\delta_2$ is forced to be $h$ since
$d\delta_2=f=dh$ and $c\delta_2=\id=ch$.
\end{proof}
In particular notice that by (\ref{2sdf-III})
discrete fibrations are fibrations
and discrete opfibrations are opfibrations as one would hope.

We now recall an analogue of the yoneda lemma
for 2-sided discrete fibrations.
Let $f:A{\rightarrow}B$ be in $\ca K$ a finitely complete
2-category. Denote by
\[ \TwoDiagRel
{\LaxSq {f/B} B A B p q {1_B} f {\lambda}}
{}
{\LaxSq {B/f} A B B {p'} {q'} f {1_B} {{\lambda}'}}
\]
the defining lax pullback squares for $f/B$ and $B/f$,
define $i_f:A{\rightarrow}f/B$ to be the unique map such that
$pi_f=1_A$, $qi_f=f$ and $\lambda{i_f}=\id$,
and $j_f:A{\rightarrow}B/f$ to be the unique map such that
$p'j_f=f$, $q'j_f=1_A$ and ${\lambda}'{j_f}=\id$.
\begin{theorem}\label{yon}
Let $\ca K$ be a finitely complete 2-category and $f:A{\rightarrow}B$
be in $\ca K$.
\begin{enumerate}
\item  (yoneda lemma):
for any span $(d_1,E_1,c_1)$ from $X$ to $A$
and discrete fibration $(d_2,E_2,c_2)$ from $X$ to $B$,
a map of spans
\[ f/B \comp E_1 \rightarrow E_2 \]
is determined uniquely by its composite with
\[ i_f \comp \id : f \comp E_1 \rightarrow f/B \comp E_1 \]
\label{yonI}
\item  (coyoneda lemma): for any span $(d_1,E_1,c_1)$ from $A$ to $X$
and discrete fibration $(d_2,E_2,c_2)$ from $B$ to $X$,
a map of spans
\[ E_1 \comp B/f \rightarrow E_2 \]
is determined uniquely by its composite with
\[ \id \comp j_f : E_1 \comp \rev f \rightarrow E_1 \comp B/f \]
\label{yonII}
\end{enumerate}
\end{theorem}
\begin{proof}
The coyoneda lemma is the yoneda lemma in $\co {\ca K}$,
so it suffices to prove the yoneda lemma.
By the representability of the notions involved
it suffices to prove this result for the case $\ca K=\CAT$.
In this case
the head of the span $f/B \comp E_1$ can be described as follows:
\begin{itemize}
\item  objects are 4-tuples
$(e,a,\beta:fa{\rightarrow}b,b)$
where $e \in E_1$ and $c_1e=a$.
\item  an arrow
\[ (\varepsilon,\alpha,\beta) :
(e_1,a_1,\beta_1,b_1) \rightarrow (e_2,a_2,\beta_2,b_2) \]
consists of maps $\varepsilon:e_1{\rightarrow}e_2$,
$\alpha:a_1{\rightarrow}a_2$
and $\beta:b_1{\rightarrow}b_2$,
such that $c_1\varepsilon=\alpha$ and $\beta\beta_1=\beta_2f(\alpha)$.
\end{itemize}
and the left and right legs of the span send $(\varepsilon,\alpha,\beta)$
described above to $d_1\varepsilon$ and $\beta$ respectively.
The image of $i_f \comp \id$ is the full subcategory
given by the $(e,a,\beta,b)$
such that $\beta=\id$.
Let $\phi:f/B \comp E_1 \rightarrow E_2$.
For any object $(e,a,\beta,b)$ we have a map
\[ (1,1,\beta) : (e,a,\id,fa) \rightarrow (e,a,\beta,b) \]
and since $d_2\phi(1,1,\beta)=\id$ and $c_2\phi(1,1,\beta)=\beta$,
$\phi(e,a,\beta,b)$ and $\phi(1,1,\beta)$
are defined as the unique ``right'' lift of $(\phi(e,a,1_{fa},fa),\beta)$
since $E_2$ is a discrete fibration.
For any arrow $(\varepsilon,\alpha,\beta)$ as above, we have
a commutative square
\[ \xymatrix{{\phi(e_1,a_1,\id,fa_1)} \ar[r]^-{\phi(1,1,\beta_1)}
\ar[d]_{\phi(\varepsilon,\alpha,f\alpha)}
& {\phi(e_1,a_1,\beta_1,b_1)} \ar[d]^{\phi(\varepsilon,\alpha,\beta)} \\
{\phi(e_2,a_2,\id,fa_2)} \ar[r]_-{\phi(1,1,\beta_2)}
& {\phi(e_2,a_2,\beta_2,b_2)}} \]
in $E_2$, but from the proof of theorem(\ref{thm-2-sided-fib})
we know that $\phi(1,1,\beta_1)$ is $d_2$-opcartesian
and so the rest of the above square is determined uniquely by
$\phi(\varepsilon,\alpha,f\alpha)$.
\end{proof}
The yoneda lemma of \cite{Str74} is theorem(\ref{yon})(\ref{yonI})
in the case where $E_1$ is the identity span.
We use an analogue of this more general result
in the proof of theorem(\ref{thm:om-cc}).
\begin{corollary}\label{cor:yon}
Let
\[ \xymatrix{A \ar[rr]^-{f} \ar[dr]_{g} && B \ar[dl]^{h} \\ & C} \]
be in $\ca K$. Then composition with $i_f$
\[ \xymatrix{{\ca K(f/B,C)(gp,hq)} \ar[r]^-{(-){\comp}i_f} &
{\ca K(A,C)(g,hf)}} \]
is a bijection with inverse given by pasting with $\lambda$.
\end{corollary}
\begin{proof}
2-cells $gp{\rightarrow}hq$ are in bijection with span maps
$f/B{\rightarrow}g/h$ by the definition of $g/h$.
By theorem(\ref{yon}) these are in bijection with
span maps $f{\rightarrow}g/h$ which by the definition of $g/h$
are in bijection with $2$-cells $g{\rightarrow}hf$.
This composite bijection is clearly $(-){\comp}i_f$.
By the definition of $i_f$ the inverse of this bijection
is given by pasting with $\lambda$.
\end{proof}
In order to define duality involutions
we shall now return to our discussion of $\Span(\ca K)$.
The operation
\[ \TwoDiagRel
{\xymatrix{B & E \ar[l]_-{c} \ar[r]^-{d} & A}}
{\mapsto}
{\xymatrix{A & E \ar[l]_-{d} \ar[r]^-{c} & B}} \]
of reversing spans is part of an isomorphism of bicategories
\[ \rev {(-)}:\op {\Span(\ca K)} \rightarrow \Span(\ca K) \]
which is the identity on objects
and provides isomorphisms of 2-categories on the homs.
Cartesian product on $\ca K$ extends to a tensor product
on $\Span(\ca K)$, and this is compatible with $\rev {(-)}$
because the operation
\[ \TwoDiagRel
{\xymatrix{{A{\times}B} & E \ar[l]_-{(p,q)} \ar[r]^-{r} & C}}
{\mapsto}
{\xymatrix{A & E \ar[l]_-{p} \ar[r]^-{(q,r)} & {B{\times}C}}} \]
provides the object part of isomorphisms of 2-categories
\begin{equation}\label{span-duality}
\Span(\ca K)(A{\times}B,C) \iso
\Span(\ca K)(A,{\rev B}{\times}C)
\end{equation}
which are pseudo natural in $A$,
$B$ and $C${\footnotemark{\footnotetext{
It's the sense in which these isomorphisms are natural
in $B$ which neccessitates the  $\rev {(-)}$ here.
The objects $B$ and $\rev B$ are of course equal.}}}.
Notice that $\rev {(-)}$ does not restrict
to 2-sided discrete fibrations: for instance
the map $d$ in a discrete fibration $(d,E,c)$
is a fibration although not in general an opfibration.
However from theorem(\ref{thm-2-sided-fib})(\ref{2sdf-I})-(\ref{2sdf-II})
and using the lifting properties from the definition
of 2-sided discrete fibration,
one can verify that the formation of
categories of 2-sided discrete fibrations provides a pseudo functor
\[ \DFib(\ca K)(-,-) : {\coop {\ca K}} \times {\op {\ca K}}
\rightarrow \CAT \]
the pseudoness arising because of span composition.
\begin{definition}\label{duality}
Let $\ca K$ be a finitely complete $2$-category.
A \emph{duality involution} for $\ca K$ consists of a $2$-functor
\[ {\iop {(-)}}: {\co {\ca K}} \rightarrow {\ca K} \]
such that $\iop {({\iop {(-)}})}=\id$ and
for all $A$, $B$ and $C$ in $\ca K$ equivalences of categories
\[ \DFib(\ca K)(A \times B, C) \catequiv
\DFib(\ca K)(A,{\iop B} \times C) \]
pseudo natural in $A$, $B$ and $C$.
\end{definition}
\begin{example}\label{op-CAT}
Let $(d,E,c)$ be a discrete fibration from $A$ to $B$ in $\CAT$.
Then for each $a \in A$ and $b \in B$ one can define
\[ E(a,b) = \{e \in E : \textnormal{\,\,\,
$de=a$ and $ce=b$}\} \]
which is contravariant in $a$ and covariant in $b$ and so defines
a functor $\op A \times B \rightarrow \SET$.
On the other hand given $P:\op A \times B \rightarrow \SET$
one can define a discrete fibration from $A$ to $B$ whose domain
is the category
\begin{itemize}
\item  objects: triples $(a,e,b)$ where $a \in A$, $b \in B$ and
$e \in P(a,b)$.
\item  a morphism $(a_1,e_1,b_1){\rightarrow}(a_2,e_2,b_2)$
consists of $\alpha:a_2{\rightarrow}a_1$ and $\beta:b_1{\rightarrow}b_2$
such that $P(\alpha,\beta)(e_1)=e_2$.
\end{itemize}
These constructions define the object maps of equivalences of categories
\[ \DFib(\CAT)(A,B) \catequiv [\op A \times B, \SET] \]
which are pseudo natural in $A$ and $B$.
Using these equivalences and the isomorphisms
\[ [\op {(A \times B)} \times C, \SET]
\iso [\op A \times (\op B \times C), \SET] \]
one can exhibit $\op {(-)}$ as a duality involution for $\CAT$.
\end{example}
\begin{example}\label{op-internal-CAT}
Let $\ca E$ be a category with finite limits.
Define $\iop A$ for $A \in \Cat(\ca E)$
by interchanging the source and target maps.
For $A, B \in \Cat(\ca E)$ there is a forgetful functor
\[ \DFib(\Cat(\ca E))(A,B) \rightarrow \Span(\ca E)(A_0,B_0) \]
which is monadic. The relevant monad on $\Span(\ca E)(A_0,B_0)$
is given by span composition $B \comp - \comp A$ regarding
$A$ and $B$ as monads in $\Span(\ca E)$.
One can verify these facts directly in the case $\ca E{=}\SET$.
Using the internal language of $\ca E$ (again see \cite{JohElpI})
this verification may in fact be interpretted in $\ca E$
to give a proof of the general result.
See also \cite{Str80} for a derivation of these facts
from the viewpoint of more general 2-sided fibrations.
From this monadicity and the description of the relevant monad,
the isomorphisms (\ref{span-duality}) above
(in the case $\ca K{=}\ca E$)
lift through the forgetful functors
to provide equivalences of categories
\[ \DFib(\Cat(\ca E))(A \times B, C) \catequiv
\DFib(\Cat(\ca E))(A,{\iop B} \times C) \]
pseudo natural in $A$, $B$ and $C$.
\end{example}

\subsection{Left extensions and left liftings}
As we shall see later in this paper,
to express the cocompleteness of an object $A$
of a finitely complete $2$-category $\ca K$,
one requires pointwise left extensions and left liftings.
We recall these notions in this subsection as well as some results
about left extending along fully faithful maps and opfibrations. 

A $2$-cell 
\[ \TriTwoCell A B C f g h {\phi} \]
in a $2$-category $\ca K$
exhibits $h$ as a \emph{left extension} of $f$ along $g$ when $\forall k$
pasting with $\phi$:
\[ \TwoDiagRel
{*{{\kappa}:h{\implies}k}}
{\mapsto}
{\xymatrix{A \ar[rr]^-{g} \ar[dr]_{f} \save \POS?="domphi" \restore
&& B \ar[dl]|{h} \save \POS?="codphi" \restore \save \POS?="domkap" \restore
\ar@/^3pc/[dl]^{k} \save \POS?="codkap" \restore \\ & C
\POS "domphi"; "codphi" **@{} ?(.35) \ar@{=>}^{\phi} ?(.65)
\POS "domkap"; "codkap" **@{} ?(.35) \ar@{=>}^{\kappa} ?(.65)}} \]
provides a bijection between $2$-cells $h{\implies}k$ and $2$-cells $f{\implies}kg$.
This left extension is \emph{preserved} by $j:C{\rightarrow}D$
when $j{\phi}$ exhibits $jh$ as a left extension of $jf$ along $g$,
and is \emph{absolute} when it is preserved by all arrows out of $C$.
The $2$-cell $\phi$ exhibits $g$ as a \emph{left lifting} of $f$ along $h$
when $\forall k$ pasting with $\phi$:
\[ \TwoDiagRel
{*{\kappa:g \implies k}}
{\mapsto}
{\xymatrix{A \ar[rr]|{g} \save \POS?="domkap" \restore
\ar@/^2pc/[rr]^{k} \save \POS?="codkap" \restore
\ar[dr]_{f} \save \POS?="domphi" \restore
&& B \ar[dl]^{h} \save \POS?="codphi" \restore \\ & C
\POS "domphi"; "codphi" **@{} ?(.35) \ar@{=>}^{\phi} ?(.65)
\POS "domkap"; "codkap" **@{} ?(.25) \ar@{=>}^{\kappa} ?(.75)}} \]
provides a bijection between $2$-cells $g{\implies}k$ and $2$-cells $f{\implies}hk$.
This left lefting is \emph{respected} by $j:D{\rightarrow}A$
when ${\phi}j$ exhibits $gj$ as a left lifting of $fj$ along $h$
and is absolute when it is respected by all arrows into $C$.
Clearly a left extension in $\ca K$ is
a left lifting in $\op {\ca K}$.
Applying the above definitions to the $2$-category $\co {\ca K}$
gives the notions of right extension and right lifting.
Some basic elementary facts regarding
left extensions in $\ca K$ are:
\begin{enumerate}
\item  If $f$ is an isomorphism and $\phi:y{\rightarrow}xf$
then $\phi$ is a left extension along $f$ iff $\phi$
is an isomorphism.
\item  If $\phi$ and $\phi'$ are left extensions along $f$,
and $\beta$ is the unique $2$-cell making
\[ \xymatrix{y \ar[r]^-{\phi} \ar[d]_{\alpha} & {xf} \ar[d]^{{\beta}f} \\
y' \ar[r]_-{{\phi}'} & {x'f}} \]
commute; then if $\alpha$ is an isomorphism so is $\beta$.
\item  If $\phi_1$ is a left extension along $g$
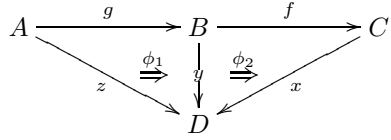
\[ \xymatrix{A \ar[rr]^-{g} \ar[drr]_{z} \save \POS?="domphi1" \restore
&& B \ar[rr]^-{f} \ar[d]|{y} \save \POS?="codphi1" \restore
\save \POS?="domphi2" \restore
&& C \ar[dll]^{x} \save \POS?="codphi2" \restore \\ && D
\POS "domphi1"; "codphi1" **@{}; ?(.35) \ar@{=>}^{\phi_1} ?(.65)
\POS "domphi2"; "codphi2" **@{}; ?(.35) \ar@{=>}^{\phi_2} ?(.65)} \]
then $\phi_2:y{\rightarrow}xf$ is a left extension along $f$
iff the composite is a left extension along $fg$.
\end{enumerate}
and applying these observations to the various duals of $\ca K$
gives analogous statements for left liftings, right extensions and right liftings.
\begin{example}\label{adj-univ}
The unit of an adjunction provides the most basic example of a left lifting
and left extension, and this expresses its
universal nature $2$-categorically: let
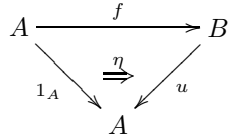
\[ \TriTwoCell A B A {1_A} f u {\eta} \]
be a $2$-cell in a $2$-category $\ca K$,
then it is an elementary exercise to show that
the following are equivalent:
\begin{enumerate}
\item  $\eta$ is the unit of an adjunction $f \ladj u$.
\item  $\eta$ exhibits $u$ as a left extension of $1_A$ along $f$
and this left extension is preserved by $f$.\label{unit-lext}
\item  $\eta$ exhibits $u$ as a left extension of $1_A$ along $f$
and this left extension is absolute.\label{unit-abslext}
\end{enumerate}
Applying this observation in $\op {\ca K}$ one obtains
two further equivalent conditions for $\eta$ in terms of left liftings,
and considering $\co {\ca K}$ one obtains
the corresponding conditions for the counit of an adjunction
in terms of right extensions and right liftings.
In these terms it is easy to see that left adjoints
preserve left extensions.
For suppose that \[ \TriTwoCell X Y A g h k {\phi} \]
exhibits $k$ as a left extension of $g$ along $h$.
Then composition with $\eta$ gives a bijection between 2-cells
$fk{\rightarrow}r$ and 2-cells $k{\rightarrow}ur$
since $\eta$ is an absolute left lifting.
Since $\phi$ is a left extension pasting with it gives a bijection
between such 2-cells and 2-cells $g{\rightarrow}urh$.
Pasting $\eta$ gives a bijection between these and 2-cells
$fg{\rightarrow}rh$.
The composite bijection is composition with $f\phi$,
and so $f\phi$ is indeed a left extension.
\end{example}
\begin{example}\label{ff}
Another basic example of absolute left liftings is provided by
fully faithful maps.
A map $f:A{\rightarrow}B$ in $\ca K$ is \emph{fully faithful}
when for all $X \in \ca K$, $\ca K(X,f)$ is a fully faithful
functor.
It is almost a tautology that
$f$ is fully faithful iff the identity cell
\[ \TriTwoCell A A B f {1_A} f {\id} \]
exhibits $1_A$ as an absolute left lifting of $f$ along itself.
\end{example}
\begin{example}\label{lexfib}
Let $f:A{\rightarrow}B$ be an opfibration and $g:C{\rightarrow}B$
in a finitely complete $2$-category $\ca K$.
In the proof of theorem(\ref{Chev-crit})
carried out in $\co {\ca K}$ so that it applies to opfibrations,
we factored the lax pullback through the pullback
\[ \TwoDiagRel
{\LaxSq {f/g} X A B {p_1} {q_1} g f {\lambda}}
{=}
{\xymatrix{{f/g} \ar[dd]_{p_1} \save \POS?="domlamb" \restore
\ar[dr]_{r} \ar[rr]^-{q_1} \save \POS?="codeqt" \restore
&& X \ar[dd]^{g} \save \POS?(.75)="codeqr" \restore
\\ & {f/_=g} \save \POS="domeqt" \restore
\save \POS="codlamb" \restore \ar[dl]^{p_2}
\save \POS?="domeqr" \restore \ar[ur]_{q_2} \\ A \ar[rr]_-{f} && B
\POS "domlamb"; "codlamb" **@{};
?(.25) \ar@{=>}^{\overline{\lambda}} ?(.5)
\POS "domeqt"; "codeqt" **@{}; ?*{=}
\POS "domeqr"; "codeqr" **@{}; ?(.6)*{=}}} \]
to obtain the $f$-opcartesian $2$-cell $\overline{\lambda}$.
Since the unique $\eta:1{\rightarrow}ir$
such that $p_1\eta=\overline{\lambda}$
and $q_1\eta=\id$ is the unit of an adjunction $r \ladj i$ by
the proof of theorem(\ref{Chev-crit}),
$\eta$ exhibits $i$ as an absolute left extension along
$r$ by the example(\ref{adj-univ}).
Thus $\overline{\lambda}$ exhibits $p_2$
as an absolute left extension of $p_1$ along $r$.
This is the key observation for theorem(\ref{cofib-lex}) below.
\end{example}
As for the description of limits and colimits, the basic example is to take
a functor $f:A{\rightarrow}B$ and then
a left extension of $f$ along $t_A:A{\rightarrow}1$ is a colimit for $f$ in $B$,
and a right extension of $f$ along $t_A$ is a limit for $f$ in $B$.
In practise to compute a left extension of $L$ of $f$ along $g:A{\rightarrow}C$
using colimits in $B$, one has the famous formula
\[ L(c) = \colim(\xymatrix{{f/c} \ar[r] & A \ar[r]^{f} & B}) \]
due to Bill Lawvere. Diagramatically we have
\[ \xymatrix{{f/c} \ar[rr] \ar[d]_{p} \save \POS?="domlam" \restore
&& 1 \ar[d]^{c} \save \POS?="codlam" \restore \\
A \ar[rr]|{g} \ar[dr]_{f} \save \POS?="domphi" \restore
&& C \ar[dl]^{L} \save \POS?="codphi" \restore \\ & B
\POS "domlam"; "codlam" **@{} ?(.42) \ar@{=>}^{\lambda} ?(.58)
\POS "domphi"; "codphi" **@{} ?(.35) \ar@{=>}^{\phi} ?(.65)} \]
in $\CAT$
where $\lambda$ is a lax pullback, $\phi$ exhibits $L$ as a left extension,
and Lawvere's formula means that the composite $2$-cell exhibits
$Lc$ as a left extension along $fp$.
This basic example suggests that the left extensions that arise
in mathematical practise remain left extensions when pasted with lax pullbacks
in this way.

To this end a $2$-cell
\[ \TriTwoCell A C B f g h {\phi} \]
in a finitely complete $2$-category $\ca K$ is defined to be a
\emph{pointwise left extension} when for all $c:X{\rightarrow}C$,
the composite
\[ \xymatrix{{g/c} \ar[rr]^-{q} \ar[d]_{p} \save \POS?="domlam" \restore
&& X \ar[d]^{c} \save \POS?="codlam" \restore \\
A \ar[rr]|{g} \ar[dr]_{f} \save \POS?="domphi" \restore
&& C \ar[dl]^{h} \save \POS?="codphi" \restore \\ & B
\POS "domlam"; "codlam" **@{} ?(.42) \ar@{=>}^{\lambda} ?(.58)
\POS "domphi"; "codphi" **@{} ?(.35) \ar@{=>}^{\phi} ?(.65)} \]
exhibits $hc$ as a left extension of $fp$ along $q$.
Such a $\phi$ is automatically a left extension \cite{Str74}:
to see this apply the definition
with $c=1_C$ and use lemma(\ref{lem:ple-ff}) below.

The remainder of this section is devoted to explaining how pointwise
left extending along fully faithful maps and opfibrations is well-behaved.

We explained what it means for $g:A{\rightarrow}B$ to be fully faithful
in example(\ref{ff}). When $\ca K$ has finite limits one can form
$H_g:A/A{\rightarrow}g/g$ as the unique map such that
\[  \begin{array}{lcccr}
p_2H_g=p_1 && q_2H_g=q_1 && \lambda_2H_g=g\lambda_1
\end{array} \]
where
\[ \TwoDiagRel
{\LaxSq {A/A} A A A {p_1} {q_1} {1_A} {1_A} {\lambda_1}}
{}
{\LaxSq {g/g} A A B {p_2} {q_2} g g {\lambda_2}} \]
are lax pullbacks. Then $g$ is fully faithful iff
$H_g$ is an isomorphism: this statemement is easily seen as true
in the case $\ca K=\CAT$ by direct inspection, and the general
case follows by a representable argument.
\begin{lemma}\label{lem:ple-ff}
The bijection of corollary(\ref{cor:yon})
sends 2-cells $gp{\rightarrow}hq$ which exhibit $h$ as a left extension
along $q$, to 2-cells which exhibit $h$ as a left extension along $f$.
\end{lemma}
\begin{proof}
This bijection clearly respects composition with $2$-cells
$h{\rightarrow}h'$.
\end{proof}
\begin{proposition}\label{ple-ff}
\cite{Str74}
If  \[ \TriTwoCell A B C f g h {\phi} \]
exhibits $h$ as a pointwise left extension of $f$ along $g$
and $g$ is fully faithful,
then $\phi$ is invertible.
\end{proposition}
\begin{proof}
Since $H_g$ is invertible the composite
\[ \xymatrix{{A/A} \ar@/_{1.5pc}/[ddr]_{p_1}
\save \POS?="domeq1" \restore
\ar[dr]|{H_g} \ar@/^{1pc}/[drrr]^{q_1} \save \POS?="codeq2" \restore \\
& {g/g} \save \POS="codeq1" \restore \save \POS="domeq2" \restore
\ar[d]|{p_2} \save \POS?="domlam2" \restore \ar[rr]|{q_2}
&& A \ar[d]^{g} \save \POS?="codlam2" \restore \\
& A \ar[rr]|{g} \ar[dr]_{f} \save \POS?="domphi" \restore
&& B \ar[dl]^{h} \save \POS?="codphi" \restore \\ && C
\POS "domlam2"; "codlam2" **@{} ?(.42) \ar@{=>}^{\lambda_2} ?(.58)
\POS "domphi"; "codphi" **@{} ?(.35) \ar@{=>}^{\phi} ?(.65)
\POS "domeq1"; "codeq1" **@{}; ?*{=}
\POS "domeq2"; "codeq2" **@{}; ?*{=}} \]
exhibits $hg$ as a left extension of $fp_1$ along $q_1$,
and so by lemma(\ref{lem:ple-ff}) $\phi$ exhibits $f$ as a
left extension along $1_A$, whence $\phi$ is an isomorphism.
\end{proof}
When taking a pointwise left extension along an opfibration,
lax pullbacks may be replaced by pullbacks for the sake of computation.
\begin{theorem}\label{cofib-lex}
\cite{Str74}
Let $\ca K$ be a finitely complete $2$-category and
\[ \TriTwoCell A C B f g h {\phi} \]
be in $\ca K$ such that $g$ is an opfibration.
Then $\phi$ exhibits $h$ as a pointwise left extension of $f$ along $g$
iff for all $c:X{\rightarrow}C$ the $2$-cell $\phi{p}$ where
\[ \PbSq P X C A p q c g \]
exhibits $hc$ as a left extension of $fp$ along $q$.
\end{theorem}
\begin{proof}
For $c:X{\rightarrow}C$ we have
\[ \TwoDiagRel
{\xymatrix{{g/c} \ar[rr]^-{q_1} \ar[d]_{p_1} \save \POS?="domlam" \restore
&& X \ar[d]^{c} \save \POS?="codlam" \restore \\
A \ar[rr]|{g} \ar[dr]_{f} \save \POS?="domphi" \restore
&& C \ar[dl]^{t} \save \POS?="codphi" \restore \\ & B
\POS "domlam"; "codlam" **@{} ?(.42) \ar@{=>}^{\lambda} ?(.58)
\POS "domphi"; "codphi" **@{} ?(.35) \ar@{=>}^{\tau} ?(.65)}}
{=}
{\xymatrix{{g/c} \ar[dr]|{r} \ar[rr]^-{q_1} \save \POS?="codeqt" \restore
\ar[dd]_{p_1} \save \POS?="domlam" \restore && X \ar[dd]^{c}
\save \POS?(.75)="codeqb" \restore \\ & P
\save \POS="domeqt" \restore \save \POS="codlam" \restore
\ar[dl]^{p} \save \POS?="domeqb" \restore
\ar[ur]_{q} \\ A \ar[rr]|{g} \ar[dr]_{f} \save \POS?="domphi" \restore
&& C \ar[dl]^{h}  \save \POS?="codphi" \restore \\ & B
\POS "domeqt"; "codeqt" **@{}; ?(.65)*{=}
\POS "domeqb"; "codeqb" **@{}; ?(.65)*{=} 
\POS "domlam"; "codlam" **@{} ?(.25) \ar@{=>}^{\overline{\lambda}} ?(.5)
\POS "domphi"; "codphi" **@{} ?(.35) \ar@{=>}^{\phi} ?(.65)}} \]
by the factorisation of lax pullbacks described in
example(\ref{lexfib})
where $\overline{\lambda}$ exhibits $p$ as an absolute left extension
of $p_1$ along $r$.
Thus $f\overline{\lambda}$ exhibits $fp$ as a left extension of $fp_1$ along $r$.
Thus by the elementary properties of left extensions,
the composite of $\phi$ and $\lambda$ exhibits $hc$ as a left extension along $q_1$
iff ${\phi}p$ exhibits $hc$ as a left extension along $q$.
\end{proof}
The next result expresses the
consequence that \emph{pointwise} left extensions are ``closed''
under the operation of pasting with lax pullback squares,
and pointwise left extensions along opfibrations are in addition
closed under pasting with pullback squares.
\begin{corollary}\label{all-ple}
\cite{Str74}
Let $\ca K$ be a finitely complete $2$-category and
\[ \TwoDiagRel {\TriTwoCell A C B f g h {\phi}} {}
{\LaxSq P X A C p q c g {\lambda}} \]
be in $\ca K$. Suppose that $\phi$ exhibits $h$ as a pointwise
left extension of $f$ along $g$.
\begin{enumerate}
\item  If $\lambda$ exhibits $P$ as $g/c$ then the composite
of $\lambda$ and $\phi$ exhibits $hc$
as a pointwise left extension along $q$.
\label{ple-ple}
\item  If $\lambda$ exhibits $P$ as the pullback
of $g$ and $c$ and $g$ is an opfibration,
then the composite of $\lambda$ and $\phi$ exhibits $hc$
as a pointwise left extension along $q$.
\label{ple-pb}
\end{enumerate}
\end{corollary}
\begin{proof}
In both (\ref{ple-ple}) and (\ref{ple-pb}) note that $q$ is an opfibration.
For (\ref{ple-ple}): obtain $q$ by first lax pulling back $g$ along $1_C$,
which is the free opfibration on $g$,
and then pulling back the result along $c$,
which is an opfibration since opfibrations are pullback stable.
For (\ref{ple-pb}): $q$ is an opfibration
since opfibrations are pullback stable.
Thus both results follow immediately from
the previous theorem and the basic properties
of pullbacks and lax pullbacks.
\end{proof}

\section{Yoneda Structures}\label{sec:yon-str}

There are two main sources of examples of yoneda structures:
(1) internal category theory; and (2) enriched
category theory.
Those that arise from internal category theory
are more nicely behaved. They satisfy a further axiom,
axiom(3*) of \cite{SW78}, and the left extensions
that form part of the axiomatics are all pointwise in the sense
discussed above.
We restrict our attention to this special case
in the following definition,
and later isolate a further condition
-- cocompleteness of presheaves.
\begin{definition}\label{yon-str}
A \emph{good yoneda structure} on a finitely complete $2$-category
$\ca K$ consists of
\begin{enumerate}
\item  A right ideal of $1$-cells
called \emph{admissible} maps.
\item  An object $A$ of $\ca K$ is admissible
when $1_A$ is so.{\footnotemark{\footnotetext{
Notice that if $B$ is admissible and $f:A{\rightarrow}B$,
then $f$ is admissible since the locally small maps
form a right ideal and of course $f=1_Bf$.}}}
For each such object, an object $\PSh A$ and an admissible map
\[ y_A : A \rightarrow \PSh A \]
is provided.
\item  For each $f:A{\rightarrow}B$ with both $A$ and $f$ admissible,
an arrow $B(f,1)$ and a $2$-cell
\[ \TriTwoCell A B {\PSh A} {y_A} f {B(f,1)} {\chi^f} \]
is provided.
\end{enumerate}
This data must satisfy the following axioms:
\begin{enumerate}
\item  $\chi^f$ exhibits $f$ as an absolute left lifting of $y_A$
through $B(f,1)$.
\item  If $A$ and $f:A{\rightarrow}B$ are admissible and
\[ \TriTwoCell A B {\PSh A} {y_A} f g {\phi}, \]
exhibits $f$ as an absolute left lifting of $y_A$ along $g$,
then $\phi$ exhibits $g$ as a pointwise left extension of $y$ along $f$.
\label{ax3*}
\end{enumerate}
\end{definition}
A yoneda structure in the sense of \cite{SW78} is defined as above
except: one replaces axiom(\ref{ax3*})
by an axiom which asserts only that $\chi^f$
is a left extension along $f$,
and proposition(\ref{compn-chi}) below is taken as an axiom.
Moreover in \cite{SW78} the hypothesis
that $\ca K$ have finite limits is not needed.
\begin{lemma}\label{ax2-iff}
Let $A,f,g,\phi$ be given as in axiom(\ref{ax3*}).
If $\phi$ exhibits
$g$ as a left extension of $y_A$ along $f$;
then this left extension is pointwise
and $\phi$ exhibits
$f$ as an absolute left lifting of $y_A$ along $g$.
\end{lemma}
\begin{proof}
Suppose that $\phi$ is a left extension.
Then since $f$ is admissible there is an isomorphism $B(f,1){\iso}g$
whose composite with $\chi^f$ is $\phi$, whence
$\phi$ is an absolute left lifting since $\chi^f$ is.
Moreover $\phi$ is a pointwise left extension since $\chi^f$ is.
\end{proof}
\begin{example}\label{CAT-yon}
We shall see later in example(\ref{CAT-yonII}),
that $\CAT$ has the following good yoneda structure,
which is \emph{the} basic example.
A functor $f:A{\rightarrow}B$ is admissible when
$\forall a,b$, the homset $B(fa,b)$ is small.
Admissible objects are locally small categories.
For such a category $A$, $\PSh A = [\op A, \Set]$ and $y_A$
is the yoneda embedding. For a functor $f:A{\rightarrow}B$
such that $A$ and $f$ are admissible, $B(f,1)(b)(a) = B(fa,b)$ and
$(\chi^f)_a:A(-,a){\rightarrow}B(f-,fa)$ is given by the arrow maps of $f$.
\end{example}
Inspired by this example one should regard the 2-cell $\chi^f$
as ``the arrow maps of f'' and indeed the notation
is selected to encourage this idea.
Given $f:A{\rightarrow}B$ with $A$ and $f$ admissible
and $x:X{\rightarrow}B$ we denote by $B(f,x)$
the composite $B(f,1)x$.
When in addition $B$ is admissible,
denote by $\res_f$ the map ${\PSh A}(y_Bf,1)$.
Given $z:Z{\rightarrow}A$ with $Z$ admissible
the following proposition
justifies denoting the composite $\res_zB(f,1)$ as $B(fz,1)$.
In fact the next two results enforce these hom-set interpretations.
\begin{proposition}\label{compn-chi}
\cite{SW78}
\begin{enumerate}
\item  If $A$ is admissible then
\[ \TriTwoCell A {\PSh A} {\PSh A} {y_A} {y_A} 1 {\id} \]
exhibits $1$ as a left extension of $y_A$ along $y_A$.
\label{yon-dense}
\item  Let \[ \xymatrix{A \ar[r]^-{f} & B \ar[r]^-{g} & C} \]
where $A,B,f$ and $g$ are admissible, then the composite
\[ \xymatrix{A \ar[r]^-{f} \ar[d]_{y} \save \POS?="domchif" \restore
& B \ar[d]|{y} \save \POS?="codchif" \restore
\save \POS?="domchig" \restore
\ar[r]^-{g} & C \ar[dl]^{B(g,1)} \save \POS?="codchig" \restore \\
{\PSh A} & {\PSh B} \ar[l]^-{\res_f}
\POS "domchif"; "codchif" **@{}; ?(.35) \ar@{=>}^{\chi^{yf}} ?(.65)
\POS "domchig"; "codchig" **@{}; ?(.35) \ar@{=>}^{\chi^g} ?(.65)} \]
exhibits $res_fB(g,1)$ as a left extension of $y_A$ along $gf$.
\label{yon-comp}
\end{enumerate}
\end{proposition}
\begin{proof}
Trivially, the identity $2$-cell in (\ref{yon-dense})
exhibits $y_A$ as an absolute left lifting of $y_A$ along
$1_{\PSh A}$.
By the composability of left liftings,
the composite $2$-cell in (\ref{yon-comp}) exhibits
$gf$ as an absolute left lifting of $y_A$.
The result follows from axiom(\ref{ax3*}).
\end{proof}
\begin{corollary}\label{ff-chi}
\cite{SW78}
\begin{enumerate}
\item  If $A$ is admissible then $y_A$ is fully faithful.
\label{ff-chi1}
\item  For $f:A{\rightarrow}B$ with $A$ and $f$ admissible,
$f$ is fully faithful iff $\chi^f$ is an isomorphism.
\label{ff-chi2}
\end{enumerate}
\end{corollary}
\begin{proof}
(\ref{ff-chi1}):
The identity cell
\[ \TriTwoCell A A {\PSh A} {y_A} {1_A} {y_A} {\id} \]
is trivially a left extension and thus an absolute
left lifting by lemma(\ref{ax2-iff}).
The result follows by example(\ref{ff}).\\
(\ref{ff-chi2})($\implies$): by proposition(\ref{ple-ff}).\\
(\ref{ff-chi2})($\impliedby$):
Since $\chi^f$ is an isomorphism and $y_A$ is fully faithful,
$\chi^f$ exhibits $1_A$ as an absolute left lifting of $y_A$ along
$B(f,f)$. Thus in the following diagram
\[ \xymatrix{A \ar[dd]_{y_A} \save \POS?(.75)="domchif" \restore
\ar[rr]^-{1_A} 
\ar[ddrr]|{f} \save \POS?(.75)="codchif" \restore
\save \POS?(.25)="domphi" \restore
&& A \ar[dd]^{f} \save \POS?(.25)="codphi" \restore
\\ \\ {\PSh A} && B \ar[ll]^-{B(f,1)} 
\POS "domchif"; "codchif" **@{} ?(.35) \ar@{=>}^{\chi^f} ?(.65)
\POS "domphi"; "codphi" **@{} ?(.35) \ar@{=>}^{\id} ?(.65)} \]
$\id$ exhibits $1_A$ as an absolute left lifting along itself,
and so the result follows by example(\ref{ff}).
\end{proof}
\noindent The notation $\res_f$ is also meant to be evocative:
in the case of example(\ref{CAT-yon}) $\res_f$ corresponds
to restriction along (that is, precomposition with) $\op f$.
So one is immediately lead to ask when $\res_f$
has left and right adjoints. The case of the right adjoints
is the easiest and we shall describe this now.
The treatment of the left adjoints is
provided in theorem(\ref{free-cocomp})
at the end of this section.
\begin{definition}\label{yon-small}
An object $C$ of $\ca K$ is \emph{small}
when both $C$ and $\PSh C$ are admissible.
\end{definition}
Given $f:A{\rightarrow}B$ with $A$ small and $B$ admissible
one can form $\PSh A(B(f,1),1)$ which we shall denote by $\ran_f$.
Let $\eta$ be the unique $2$-cell such that
\[ \TwoDiagRel
{\xymatrix{B \ar[dr]_{y} \save \POS?="domid" \restore
\ar[r]^-{y} & {\PSh B} \ar[d]|{1_{\PSh B}}
\save \POS?="codid" \restore \save \POS?="dometa" \restore
\ar[r]^-{\res_f} & {\PSh A} \ar[dl]^{\ran_f}
\save \POS?="codeta" \restore \\ & {\PSh B}
\POS "domid"; "codid" **@{} ?(.35) \ar@{=>}^{\id} ?(.65)
\POS "dometa"; "codeta" **@{} ?(.35) \ar@{=>}^{\eta} ?(.65)}}
{=} {*{\chi^{\res_fy}}} \]
Then $\eta$ exhibits $\ran_f$ as a left extension along $\res_f$.
Now observe that the composite
\[ \xymatrix{A \ar[r]^-{f} \ar[dr]_{y} \save \POS?="domchif" \restore
& B \ar[dr]_{y} \save \POS?="codchif" \restore \save \POS?="domid" \restore
\ar[r]^-{y} & {\PSh B}
\ar[d]|{1_{\PSh B}} \save \POS?="codid" \restore
\save \POS?="dometa" \restore
\ar[r]^-{\res_f} & {\PSh A} \ar[dl]^{\ran_f}
\save \POS?="codeta" \restore \\
& {\PSh A} & {\PSh B} \ar[l]^-{\res_f}
\POS "domchif"; "codchif" **@{} ?(.35) \ar@{=>}^{\chi^{yf}} ?(.65)
\POS "domid"; "codid" **@{} ?(.35) \ar@{=>}^{\id} ?(.65)
\POS "dometa"; "codeta" **@{} ?(.35) \ar@{=>}^{\eta} ?(.65)} \]
exhibits $\res_f\ran_f$ as a left extension
by proposition(\ref{compn-chi}),
and so the left extension $\eta$ is preserved by $\res_f$.
We have proved:
\begin{proposition}\label{ran}
\cite{SW78}
Given $f:A{\rightarrow}B$ with $A$ small and $B$ admissible,
$\res_f$ has right adjoint $\ran_f$.
\end{proposition}
The remainder of this subsection is devoted to
the discussion of internal colimits.
First we shall give the abstract definitions of colimits
and cocompleteness. In preparation for this consider
\[ \xymatrix@1{C \ar[r]^-{f} & A \ar[r]^-{g} & B} \]
with $C,f,g$ admissible, denote by $\chi^g_f$ the unique
2-cell such that
\[ \TwoDiagRel
{\xymatrix{C \ar[r]^-{f} \ar[dr]_{y_C}
\save \POS?="domchif" \restore
& A \ar[r]^-{g}
\ar[d] \save \POS?="codchif" \restore
\save \POS?="domchigf" \restore
& B \ar[dl]^{B(gf,1)} \save \POS?="codchigf" \restore
\\ & {\PSh C}
\POS "domchif"; "codchif" **@{}
?(.35) \ar@{=>}^{\chi^f} ?(.65)
\POS "domchigf"; "codchigf" **@{}
?(.35) \ar@{=>}^{\chi^g_f} ?(.65)}}
{=} {*{\chi^{gf}}} \]
and note that $\chi^g_f$ exhibits $B(gf,1)$ as a left extension
of $A(f,1)$ along $g$.
\begin{definition}\label{int-colimits}
\cite{SW78}
Let $\ca K$ have a good yoneda structure.
\begin{enumerate}
\item  Consider $i:M{\rightarrow}{\PSh C}$ and $f:C{\rightarrow}A$
with $M$ and $f$ admissible, and $C$ small.
A \emph{colimit} of $f$ \emph{weighted by} $i$ consists of
an admissible map $\col(i,f):M{\rightarrow}A$ together with
\[ \TriTwoCell M A {\PSh C} i {\col(i,f)} {A(f,1)} {\eta} \]
which exhibits $\col(i,f)$ as an absolute left lifting of
$i$ through $A(f,1)$.
\item  $\col(i,f)$ is \emph{preserved by} an admissible map
$g:A{\rightarrow}B$ when the composite 2-cell
\[ \xymatrix{M \ar[r]^-{\col(i,f)} \ar[dr]_{i}
\save \POS?="domchif" \restore & A \ar[r]^-{g}
\ar[d] \save \POS?="codchif" \restore
\save \POS?="domchigf" \restore
& B \ar[dl]^{B(gf,1)} \save \POS?="codchigf" \restore \\ & {\PSh C}
\POS "domchif"; "codchif" **@{} ?(.35) \ar@{=>}^{\eta} ?(.65)
\POS "domchigf"; "codchigf" **@{} ?(.35) \ar@{=>}^{\chi^g_f} ?(.65)} \]
exhibits $g(\col(i,f))$ as the absolute left lifting of $i$
along $B(gf,1)$.
\label{pres-col}
\item  $A$ is \emph{small cocomplete} when $\col(i,f)$ exists
for all $i:M{\rightarrow}{\PSh C}$ and $f:C{\rightarrow}A$
with $M$ and $f$ admissible and $C$ small.
\item  $g:A{\rightarrow}B$ is \emph{cocontinuous} when
for all $i:M{\rightarrow}{\PSh C}$ and $f:C{\rightarrow}A$
such that $\col(i,f)$ exists, $\col(i,f)$ is preserved by $g$.
\end{enumerate}
\end{definition}
\noindent
At first glance the above definition may seem
very abstract so we shall now reconcile it with the
corresponding familiar notions.
The key procedure that enables one to do this
is a bijection between 2-cells $\phi$ and ${\phi}'$ so that
\begin{equation}\label{prime-def}
\TwoDiagRel
{\xymatrix{A \ar[dd]_{y} \save \POS?(.75)="domchif" \restore \ar[rr]^-{g} 
\ar[ddrr]|{f} \save \POS?(.75)="codchif" \restore \save \POS?(.25)="domphi" \restore
&& B \ar[dd]^{h} \save \POS?(.25)="codphi" \restore
\\ \\ {\PSh A} && C \ar[ll]^-{C(f,1)} 
\POS "domchif"; "codchif" **@{} ?(.35) \ar@{=>}^{\chi^f} ?(.65)
\POS "domphi"; "codphi" **@{} ?(.35) \ar@{=>}^{\phi} ?(.65)}}
=
{\xymatrix{A \ar[dd]_{y} \save \POS?(.25)="domchig" \restore \ar[rr]^-{g} && B
\ar[ddll]|{B(g,1)} \save \POS?(.25)="codchig" \restore
\save \POS?(.75)="domphip" \restore
\ar[dd]^{h} \save \POS?(.75)="codphip" \restore
\\ \\ {\PSh A} && C \ar[ll]^-{C(f,1)}
\POS "domchig"; "codchig" **@{} ?(.35) \ar@{=>}^{\chi^g} ?(.65)
\POS "domphip"; "codphip" **@{} ?(.35) \ar@{=>}^{{\phi}'} ?(.65)}}
\end{equation}
where $A$, $f$ and $g$ are assumed admissible.
Equation(\ref{prime-def}) does indeed establish a bijection
$\phi \mapsto {\phi}'$ because $\chi^g$ is a left extension
and $\chi^f$ is a left lifting.
\begin{lemma}\label{priming}
\begin{enumerate}
\item  $\phi$ exhibits $h$ as a left extension of $f$ along $g$
iff $\phi'$ exhibits $h$ as a left lifting of $B(g,1)$ along $C(f,1)$.
\label{char-lex}
\item  $\phi$ exhibits $h$ as a pointwise left extension of $f$ along $g$
iff $\phi'$ exhibits $h$ as an absolute left lifting
of $B(g,1)$ along $C(f,1)$.
\label{char-plex}
\item  $\phi$ exhibits $g$ as an absolute left lifting of $f$ along $h$
iff $\phi'$ is an isomorphism.
\label{char-all}
\end{enumerate}
\end{lemma}
\begin{proof}
(\ref{char-lex}): the result follows since
the bijection $\phi{\mapsto}\phi'$
respects composition with 2-cells $h{\rightarrow}k$.
\\
(\ref{char-plex}): For any $b:X{\rightarrow}B$ we have
\[ \TwoDiagRel
{\xymatrix{{g/b} \ar[rr]^-{p_X} \ar[dd]_{p_A}
\save \POS?="domlam" \restore
&& X \ar[dd]^{b} \save \POS?="codlam" \restore \\ \\
A \ar[dd]_{y} \save \POS?(.75)="domchif" \restore \ar[rr]^-{g} 
\ar[ddrr]|{f} \save \POS?(.75)="codchif" \restore
\save \POS?(.25)="domphi" \restore
&& B \ar[dd]^{h} \save \POS?(.25)="codphi" \restore
\\ \\ {\PSh A} && C \ar[ll]^-{C(f,1)} 
\POS "domlam"; "codlam" **@{} ?(.4) \ar@{=>}^{\lambda} ?(.6)
\POS "domchif"; "codchif" **@{} ?(.35) \ar@{=>}^{\chi^f} ?(.65)
\POS "domphi"; "codphi" **@{} ?(.35) \ar@{=>}^{\phi} ?(.65)}}
=
{\xymatrix{{g/b} \ar[rr]^-{p_X} \ar[dd]_{p_A}
\save \POS?="domlam" \restore
&& X \ar[dd]^{b} \save \POS?="codlam" \restore \\ \\
A \ar[dd]_{y} \save \POS?(.25)="domchig" \restore \ar[rr]^-{g} && B
\ar[ddll]|{B(g,1)} \save \POS?(.25)="codchig" \restore
\save \POS?(.75)="domphip" \restore
\ar[dd]^{h} \save \POS?(.75)="codphip" \restore
\\ \\ {\PSh A} && C \ar[ll]^-{C(f,1)} 
\POS "domlam"; "codlam" **@{} ?(.4) \ar@{=>}^{\lambda} ?(.6)
\POS "domchig"; "codchig" **@{} ?(.35) \ar@{=>}^{\chi^g} ?(.65)
\POS "domphip"; "codphip" **@{} ?(.35) \ar@{=>}^{{\phi}'} ?(.65)}} \]
Now ${\phi}'$ is an absolute left lifting along $C(f,1)$
iff for all $b$ and $c:X{\rightarrow}C$,
pasting with ${\phi}'$ gives a bijection
$\ca K(X,C)(bh,c) \iso \ca K(X,\PSh A)(B(g,b),C(f,c))$.
However the composite of $\chi^g$ and $\lambda$
is a left extension along $p_X$, so this is the same as
saying that pasting with the composite $2$-cell (on either
side of the above equation) gives a bijection
$\ca K(X,C)(bh,c) \iso \ca K(g/b,\PSh A)(y_Ap_A,C(f,cp_X))$
for all $b$ and $c$.
However $\chi^f$ is an absolute left lifting along $B(f,1)$,
so this is the same as saying that pasting with the composite of $\phi$
and $\lambda$ gives a bijection
$\ca K(X,C)(bh,c) \iso \ca K(X,\PSh A)(fp_A,cp_X)$ for all $b$ and $c$,
and this is by definition the statement that $\phi$
is a pointwise left extension.
\\
(\ref{char-all}): $\phi$ is an absolute left lifting iff
composite 2-cell in equation(\ref{prime-def}) exhibits
$g$ as an absolute left lifting of $f$ along $C(f,h)$,
which by axiom(\ref{ax3*}) and lemma(\ref{ax2-iff}) is true iff
this composite exhibits $C(f,h)$ as a left extension of $y_A$ along $g$,
which is true iff $\phi'$ is an isomorphism.
\end{proof}
\noindent
Armed with this procedure we can now
reconcile definition(\ref{int-colimits})
with the usual notion of weighted colimit in terms of the
hom-set notation.
\begin{corollary}\label{col-rec}
Consider $i:M{\rightarrow}{\PSh C}$ and $f:C{\rightarrow}A$
with $M$ and $f$ admissible, and $C$ small.
Then the colimit of $f$ weighted by $i$ exists iff
there is an admissible map $\col(i,f):M{\rightarrow}A$
together with an isomorphism
\[ A(\col(i,f),1) \iso {\PSh C}(i,C(f,1)). \]
\end{corollary}
\begin{proof}
By lemma(\ref{priming})(\ref{char-all})
the isomorphism is given by $\eta'$ where $\eta$ is the defining
2-cell of the weighted colimit.
\end{proof}
\noindent
In the usual theory of weighted colimits one can express
pointwise left extensions as weighted colimits in a canonical way.
This fact is immediate in the present setting.
\begin{corollary}\label{ple-rec}
Let \[ \TriTwoCell C B A f g h {\phi} \]
such that $C$ is small and $B,f,h$ are admissible.
Then $\phi$ exhibits $h$ as a pointwise left extension of $f$ along $g$
iff the colimit of $f$ weighted by $B(g,1)$ exists
and there is an isomorphism \[ h \iso \col(B(g,1),f). \]
\end{corollary}
\begin{proof}
By lemma(\ref{priming})(\ref{char-plex})-(\ref{char-all}),
the isomorphism is given by $\phi''$.
\end{proof}
\begin{example}\label{psh-colreps}
Another basic fact about colimits in ordinary category theory
is that every presheaf is a colimit of representables.
To appreciate the general analogue of this
consider $i:M{\rightarrow}{\PSh C}$ where $C$ is admissible.
Then by proposition(\ref{compn-chi})(\ref{yon-dense})
the lax pullback square
\[ \LaxSq {y_C/i} M C {\PSh C} p q i {y_C} {\lambda} \]
exhibits $i$ as a pointwise left extension.
The well known case of this is for the yoneda structure
of example(\ref{CAT-yon}) with $M=1$ and $C$ small. 
\end{example}
\noindent
Another feature of the theory of colimits for $\CAT$
is that weighted colimits may be replaced by the more
commonly used conical ones. This too has an analogue
in a good yoneda structure in corollary(\ref{weighted-as-conical}) below.
Consider $i:M{\rightarrow}{\PSh C}$ and $f:C{\rightarrow}A$
with $M$ and $f$ admissible, and $C$ small.
Since $\lambda$ of the previous example is a left extension
and $\chi^f$ a left lifting, the equations
\[ \TwoDiagRel
{\xymatrix{{y_C/i} \ar[r]^-{q} \ar[d]_{p} \save \POS?="domlam" \restore
& M \ar[r]^-{k} \ar[d]|{i} \save \POS?="codlam" \restore
\save \POS?="dometa" \restore
& A \ar@/^{1pc}/[dl]^{A(f,1)} \save \POS?(.35)="codeta" \restore \\
C \ar[r]_-{y_C} & {\PSh C}
\POS "domlam"; "codlam" **@{}; ?(.35) \ar@{=>}^{\lambda} ?(.65)
\POS "dometa"; "codeta" **@{}; ?(.35) \ar@{=>}^{\eta} ?(.65)}}
{=}
{\xymatrix{{y_C/i} \ar[rr]^-{q} \ar[d]_{p} \save \POS?="domphi" \restore
&& M \ar[d]^{k} \save \POS?="codphi" \restore \\ C \ar[rr]|-{f}
\ar[dr]_{y_C} \save \POS?="domchif" \restore
&& A \ar[dl]^{A(f,1)} \save \POS?="codchif" \restore \\ & C
\POS "domphi"; "codphi" **@{}; ?(.42) \ar@{=>}^{{\eta}'} ?(.58)
\POS "domchif"; "codchif" **@{}; ?(.35) \ar@{=>}^{\chi^f} ?(.65)}} \]
determine a bijection between 2-cells $\eta$ and 2-cells ${\eta}'$.
\begin{lemma}\label{conical-priming}
\begin{enumerate}
\item  $\eta$ exhibits $k$ as a left lifting of $i$ along $A(f,1)$
iff $\eta'$ exhibits $k$ as a left extension of $fp$ along $q$.
\label{cpI}
\item  $\eta$ exhibits $k$ as an absolute left lifting of $i$
along $A(f,1)$ iff $\eta'$ exhibits $k$ as a pointwise left extension of
$fp$ along $q$.
\label{cpII}
\end{enumerate}
\end{lemma}
\begin{proof}
(\ref{cpI}): the bijection $\eta \mapsto \eta'$
clearly respects composition with 2-cells $k{\rightarrow}k'$. \\
(\ref{cpII}): since $q$ is an opfibration,
$\eta'$ is a pointwise left extension
iff for every $x:X{\rightarrow}M$ the composite
\[ \xymatrix{{y_C/ix} \ar[d]_{p_x} \save \POS?(.3)="lpb" \restore
\ar[r]^-{q_x} \save \POS?(.3)="tpb" \restore
& X \ar[d]^{x} \save \POS?(.3)="rpb" \restore \\
{y_C/i} \ar[r]|-{q} \save \POS?(.3)="bpb" \restore
\ar[d]_{p} \save \POS?="domphi" \restore
& M \ar[d]^{k} \save \POS?="codphi" \restore \\
C \ar[r]_-{f} & A
\POS "domphi"; "codphi" **@{}; ?(.35) \ar@{=>}^{{\eta}'} ?(.65)
\POS "rpb"; "lpb" **@{}; ?!{"bpb";"tpb"}="cpb" **@{}; ? **@{-};
"tpb"; "cpb" **@{}; ? **@{-}} \]
exhibits $kx$ as a left extension along $q_x$ by theorem(\ref{cofib-lex}),
which is equivalent by (\ref{cpI}) to $\eta{x}$ exhibiting $kx$
as a left lifting along $A(f,1)$.
\end{proof}
\noindent
An immediate consequence of lemma(\ref{conical-priming})(\ref{cpII}) is:
\begin{corollary}\label{weighted-as-conical}
Let $\ca K$ be a finitely complete 2-category equipped
with a good yoneda structure. Let
\[ \begin{array}{lcr}
{i:M{\rightarrow}{\PSh C}} && {f:C{\rightarrow}A}
\end{array} \]
be in $\ca K$
where $C$ is small and $M$ and $f$ are admissible.
Then the colimit of $f$ weighted by $i$ exists
iff there is $k$ admissible and
\[ \LaxSq {y_C/i} M C A p q k f {\phi} \]
(where $p$ and $q$ come from a lax pullback square defining $y_C/i$)
such that $\phi$ exhibits $k$ as a pointwise
left extension of $fp$ along $q$.
When this is the case $k \iso \col(i,f)$.
\end{corollary}
From example(\ref{adj-univ}) we know that
left adjoints preserve left extensions,
thus they preserve pointwise left extensions,
and so by the last corollary weighted colimits (as one would hope).
On the other hand let $f:C{\rightarrow}X$ be such that $C$ is small
and $f$ is admissible and cocontinuous.
Define $\eta$ as the unique $2$-cell satisfying
\[ \TwoDiagRel
{\xymatrix{C \ar[dr]_{y} \save \POS?="domid" \restore
\ar[r]^-{y} & {\PSh C} \ar[d]|{1_{\PSh C}}
\save \POS?="codid" \restore \save \POS?="dometa" \restore
\ar[r]^-{f} & {X} \ar[dl]^{X(fy,1)} \save \POS?="codeta" \restore \\
& {\PSh C}
\POS "domid"; "codid" **@{} ?(.35) \ar@{=>}^{\id} ?(.65)
\POS "dometa"; "codeta" **@{} ?(.35) \ar@{=>}^{\eta} ?(.65)}}
{=} {*{\chi^{fy}}} \]
so $\eta$ exhibits $X(fy,1)$ as a left extension along $f$.
Note $\chi^{fy}$ and $\id$ are actually pointwise left extensions,
and so are preserved by $f$ they are expressible as weighted colimits
which $f$ preserves by hypothesis.
Thus the left extension $\eta$ is preserved by $f$, and so by
example(\ref{adj-univ}) we have proved:
\begin{lemma}\label{lem-fc}\cite{SW78}
If $f:{\PSh C}{\rightarrow}X$ such that $C$ is small
and $f$ is admissible and cocontinuous
then $f$ has a right adjoint $X(fy,1)$.
\end{lemma}
\begin{definition}\label{psh-coco}
A good yoneda structure
is said to have \emph{cocomplete presheaves} when
for every small object $A$, the object $\PSh A$ is cocomplete.
\end{definition}
\begin{theorem}\label{free-cocomp}
Suppose $\ca K$ has a good yoneda structure with cocomplete presheaves.
\begin{enumerate}
\item  If $C$ is small then $\PSh C$ is the colimit completion of $C$;
in the sense that given $X$ admissible and cocomplete, the adjunction
\[ \xymatrix{{\ca K(C,X)} \ar@/^1pc/[rr] \save \POS?="top" \restore
&& {\ca K({\PSh C},X)} \ar@/^1pc/[ll] \save \POS?="bot" \restore
\POS "top"; "bot" **@{} ?*{\perp}} \]
given by restriction and left extension along
$y_C:C{\rightarrow}{\PSh C}$,
restricts to an equivalence of categories
\[ {\ca K(C,X)} \catequiv \CoCts({\PSh C},X) \]
where $\CoCts({\PSh C},X)$ denotes the full subcategory of
$\ca K({\PSh C},X)$ consisting of the cocontinuous maps.
\label{free-cocomp1}
\item  If $f:A{\rightarrow}B$ with $A$ small and $B$ admissible
then $\res_f$ has left adjoint $\lan_f$ defined as the left extension of
$fy_B$ along $y_A$.
\label{free-cocomp2}
\end{enumerate}
\end{theorem}
\begin{proof}
(\ref{free-cocomp1}):
The asserted adjunction exists since $\ca K$'s yoneda
structure has cocomplete presheaves.
If $f:{\PSh C}{\rightarrow}X$ is cocontinuous then it preserves
the left extension
\[ \TriTwoCell C {\PSh C} {\PSh C} y y 1 {\id} \]
so that it is the left extension along $y_C$ of $fy_C$.
Conversely, if $f$ arises as a left extension along $y_C$,
then it has a right adjoint by lemma(\ref{lem-fc}),
and so is cocontinuous.\\
(\ref{free-cocomp2}):
Defining $\lan_f:A{\rightarrow}B$ in this way note that
since $y_A$ is fully faithful we have
$\lan_fy_A \iso y_Bf$ and so by lemma(\ref{lem-fc})
$\lan_f \ladj \res_f$.
\end{proof}
%

\section{2-toposes: Definition and Examples}
\label{sec:2topos}

Recall that a category $\ca E$ is an
\emph{elementary topos} when
\begin{enumerate}
\item  it has finite limits,
\item  is cartesian closed, and
\item  has a subobject classifier.
\end{enumerate}
Finitely complete $2$-categories have been discussed
at length in section(\ref{sec:2cat-pre}),
and the definition of cartesian closed $2$-category
is the obvious direct analogue of its one dimensional counterpart:
for each $A \in \ca K$, $(- \times A):{\ca K}{\rightarrow}{\ca K}$
has a right $2$-adjoint, which we denote as $[A,-]$.

The $2$-categorical generalisation of subobject classifier
is based on an important idea of Bill Lawvere regarding $\CAT$: that
\emph{the category of sets is a generalised object of truth values}.
This idea was expressed 2-categorically in the work of
Ross Street \cite{Str74b} \cite{Str80b}
as part of the notion of a fibrational cosmos.
Following Street's approach, but allowing ourselves the luxury
of a cartesian closed $2$-category with a duality involution,
one obtains the notion of $2$-topos discussed here.

First we recall ordinary subobject classifiers.
Let $\ca E$ be a locally small finitely complete category,
and consider the functor
\[ \Sub_0 : \op {\ca E} \rightarrow \SET \]
which sends $E \in \ca E$ to the set of subobjects of $E$,
and is given on arrows by pulling back.
A \emph{subobject classifier} is a representing
object for $\Sub_0$.
This amounts to a monomorphism
$\tau_0:\Omega_{{\bullet}0}{\rightarrowtail}\Omega_0$ in $\ca E$ such that
for all $E \in \ca E$ the function
\[ {\ca E}(E,\Omega_0) \rightarrow \Sub_0(E) \]
given by pulling back $\tau_0$, is a bijection.
Since $\ca E$ is locally small $\Sub_0$ actually lands in $\Set$.
Moreover it is easily verified (or see \cite{MM} page 33 for example)
that $\Omega_{{\bullet}0}$ is a terminal object, so we will denote
it as $1$.
Subobjects form not just a set but a poset, and
as we shall recall below in example(\ref{cdc-subobj}), $\tau_0$
is the object part of an internal functor $\tau:1{\rightarrow}\Omega$.

The category $\Set$ as an object of $\CAT$ is analogous to $\Omega$.
To any functor $f:A{\rightarrow}\Set$
a version of the Grothendieck construction associates to $f$
a discrete opfibration $G(f):e(f){\rightarrow}A$.
Explicitly the category $e(f)$ can be described as having objects
pairs $(x,a)$ where $a \in A$ and $x \in f(a)$,
and having arrows $(x_1,a_1){\rightarrow}(x_2,a_2)$
which are maps $\alpha:a_1{\rightarrow}a_2$ in $A$
such that $f(\alpha)(x_1)=x_2$; and $G(f)(x,a)=a$.
Any discrete opfibration into $A$ with small fibres arises in this way,
and so we have a fully faithful functor
\[ G: \CAT(A,\Set) \rightarrow \DFib(1,A) \]
whose image consists of those discrete opfibrations with small fibres.
To complete the analogy with subobject classifiers notice that $G$
can be described more efficiently
\[ \PbSq {e(f)} {\PtSet} {\Set} A {G(f)} {} U f \]
as the process of pulling back the forgetful functor
$U:\PtSet{\rightarrow}\Set$ from the category $\PtSet$ of pointed sets
and functions which preserve the base point.
The forgetful functor $U$ is a discrete opfibration with small fibres,
and we have just recalled the sense in which it is the universal such.

Let $p:E{\rightarrow}B$ be a discrete opfibration in a finitely complete
$2$-category $\ca K$. Then just as in the case of $U$ above one has
functors
\[ G_{p,A} : \ca K(A,B) \rightarrow \DFib(1,A) \]
given by pulling back $p$.
The effect of $G_{p,A}$ on morphisms is depicted as follows:
given $\phi:f{\implies}g$ we have
\[ \xymatrix{& {e(f)} \save \POS="domphib" \restore
\ar[dl]_{G\phi} \ar[dddl]|(.65){Gf} \ar[dr] \\
{e(g)} \ar[dd]_{Gg} \ar@/_{1pc}/[rr]
\save \POS?="codphib" \restore && E \ar[dd]^{p} \\ \\
A \ar@/^{1pc}/[rr]^-{f} \save \POS?="domphi" \restore
\ar@/_{1pc}/[rr]_-{g} \save \POS?="codphi" \restore && B
\POS "domphi"; "codphi" **@{}; ?(.3) \ar@{=>}^{\phi} ?(.7)
\POS "domphib"; "codphib" **@{}; ?(.5) \ar@{=>}^{\overline{\phi}} ?(.7)} \]
where $\overline{\phi}$ is the unique lifting of $\phi$, that is,
$G(g)G(\phi)=G(f)$ and $\phi{G(g)}=p\overline{\phi}$.
\begin{definition}\label{class-dcofib}
Let $\ca K$ be a finitely complete $2$-category.
A discrete opfibration $\tau:\Omega_{\bullet}{\rightarrow}\Omega$
in $\ca K$ is \emph{classifying} when the functors
\[ G_{\tau,A} : \ca K(A,\Omega) \rightarrow \DFib(1,A) \]
are fully faithful for all $A \in \ca K$.
\end{definition}
\begin{example}\label{cdc-cat}
By the discussion preceding definition(\ref{class-dcofib})
$U:\PtSet{\rightarrow}\Set$
is a classifying discrete opfibration for $\CAT$.
There are many others.
For example let $\lambda$ be a regular cardinal,
and in the above discussion replace $\Set$ by $\Set_{\lambda}$,
the category of sets of cardinality less than $\lambda$.
Then the analogous forgetful functor
$U_{\lambda}:\Set_{\bullet\lambda}{\rightarrow}\Set_{\lambda}$
is also a classifying discrete opfibration. In this case
the image of
\[ G_{U_{\lambda},A} : \CAT(A,\Set_{\lambda}) \rightarrow \DFib(1,A) \]
consists of the discrete opfibrations with fibres of
cardinality less than $\lambda$.
These examples illustrate that each classifying discrete opfibration
provides a notion of smallness, and when one encodes category theory
internally with the aid of a given classifying discrete opfibration
as we do below, one has automatically accounted for size issues.
Notice also that the case $\lambda=2$ gives the
ordinary subobject classifier for $\SET$.
\end{example}
\begin{example}\label{cdc-subobj}
Let $\ca E$ be an elementary topos with subobject classifier
denoted as $\tau_0:1{\rightarrow}\Omega_0$.
We shall now see that $\tau_0$ is the object map of a classifying
discrete opfibration $\tau:1{\rightarrow}\Omega$.
For any $E \in \ca E$, $\Sub_0(E)$ has a maximum $1_E$,
and binary meets given by pullback
\[ \xymatrix{1 \ar[r]^-{1_E} & {\Sub_0(E)}
& {\Sub_0(E) \times \Sub_0(E)} \ar[l]_-{\wedge}} \]
and this structure is natural in $E$, and so by the
yoneda lemma induces
\[ \xymatrix{1 \ar[r]^-{\tau_0} & {\Omega_0} & {\Omega_0 \times \Omega_0}
\ar[l]_-{\wedge}} \]
in $\ca E$, satisfying the diagrams that express internally
that $\tau_0$ and $\wedge$ are the unit and multiplication for
a commutative monoid structure on $\Omega_0$ for which
every element is idempotent.
The poset structure is a consequence.
For subobjects $S_1$ and $S_2$ of $E$,
one has $S_1{\leq}S_2$ iff $S_1{\wedge}S_2=S_1$,
and so the corresponding poset structure on $\Omega$ induced by this
is expressed in $\ca E$ by an equaliser
\[ \xymatrix{{\Omega_1} \ar[r]^-{\leq} & {\Omega_0 \times \Omega_0}
\ar@<1ex>[r]^-{\pi_1} \ar@<-1ex>[r]_-{\wedge} & {\Omega_0}}. \]
That is $\Omega_1$ is the object of arrows of the internal
category (in fact poset) $\Omega$.
From this explicit description it is immediate that:
\begin{enumerate}
\item  an internal functor $f:X{\rightarrow}\Omega$ amounts
to a map $f_0:X_0{\rightarrow}\Omega_0$
such that $f_0s{\leq}f_0t$, where $s$ and $t$
are the source and target maps for the internal category $X$; and
\label{intfun-om}
\item  an internal natural transformation $\phi:f{\implies}g:X{\rightarrow}\Omega$
is unique it exists, and does so iff $f_0{\leq}g_0$.
\end{enumerate}
In particular $\tau:1{\rightarrow}\Omega$ is the internal functor
corresponding to $\tau_0$ by (\ref{intfun-om}).
For each $X$ the functor
\[ \xymatrix{{1=\ca E(X,1)} \ar[rr]^-{\ca E(X,\tau)}
&& {\ca E(X,\Omega)\iso{\Sub(X)}}} \]
picks out the maximum element of $\Sub(X)$,
and so is clearly a discrete opfibration.
Thus $\tau$ is a discrete opfibration.
We call a discrete opfibration $p:E{\rightarrow}X$ in $\Cat(\ca E)$
a \emph{cosieve} when $p_0$ is a monomorphism;
the full subcategory of $\DFib(1,X)$ given by the cosieves into $X$
is denoted by $\Cosieve(X)$.
For example $\tau:1{\rightarrow}\Omega$ is a cosieve.
Since cosieves are pullback stable
$G_{\tau,X}$ factors through $\Cosieve(X)$
and we shall now see that it provides an equivalence
$\Cat(\ca E)(X,\Omega) \catequiv \Cosieve(X)$.
Given a cosieve $p:E{\rightarrow}X$ define
$(Fp)_0:X_0{\rightarrow}\Omega_0$
to be the unique map corresponding to the subobject $p_0$
\[ \xymatrix{{E_1} \ar@<1ex>[rr]^-{s} \ar@<-1ex>[rr]|-{t} \ar[dd]_{p_1}
\ar@{.>}[dr] && {E_0} \ar[dd]^{p_0}
\save \POS?(.6)="rpb1" \restore \save \POS?(.15)="lpb2" \restore
\ar[rr] \save \POS?(.15)="tpb2" \restore && 1 \ar[dd]^{\tau_0}
\save \POS?(.15)="rpb2" \restore \\
& {E'_1} \ar[dl]_{p'_1} \save \POS?(.2)="lpb1" \restore
\ar[ur] \save \POS?(.2)="tpb1" \restore \\
{X_1} \ar@<1ex>[rr]^-{s} \ar@<-1ex>[rr]_-{t}
\save \POS?(.6)="bpb1" \restore && {X_0} \ar[rr]_-{(Fp)_0}
\save \POS?(.15)="bpb2" \restore && {\Omega_0}
\POS "rpb1"; "lpb1" **@{}; ?!{"bpb1";"tpb1"}="cpb1" **@{}; ? **@{-};
"tpb1"; "cpb1" **@{}; ? **@{-}
\POS "rpb2"; "lpb2" **@{}; ?!{"bpb2";"tpb2"}="cpb2" **@{}; ? **@{-};
"tpb2"; "cpb2" **@{}; ? **@{-}} \]
Since $p$ is a discrete opfibration, the square involving $p_1$,
$p_0$ and the source maps (both denoted as $s$) is a pullback,
and so the monomorphism $p_1$ is classified by $(Fp)_0s$.
Denoting by $E'_1$ the pullback of $p_0$ along $t$,
the resulting monomorphism $p'_1$ is classified by $(Fp)_0t$,
$p_1$ factors through $p'_1$
(as indicated by the dotted arrow) and so $(Fp)_0s{\leq}(Fp)_0t$.
Thus $(Fp)_0$ is the object map of a unique internal functor $Fp$.
If the cosieve $p$ factors through another cosieve $q$,
then by construction we have $Fp{\leq}Fq$.
By definition, $FG_{\tau,X}S=S$ for any subobject $S$
and $G_{\tau,X}Fp\iso{p}$ for any cosieve $p$.
\end{example}
\begin{proposition}\label{cdc-pb}
Let $\ca K$ be a finitely complete $2$-category,
$\tau:\Omega_{\bullet}{\rightarrow}\Omega$ be a classifying
discrete opfibration in $\ca K$,
and $f:A{\rightarrow}\Omega$.
Then $\tau_f$ defined by
\[ \PbSq {A_{\bullet}} {\Omega_{\bullet}} {\Omega} A
{\tau_f} {} {\tau} f \]
is a classifying discrete opfibration
iff $f$ is fully faithful.
\end{proposition}
\begin{proof}
By definition there is a natural isomorphism
\[ \xymatrix{{{\ca K}(X,A)} \ar[rr]^-{{\ca K}(X,f)}
\ar[dr]_{G_{{\tau_f},X}}
\save \POS?="domiso" \restore && {{\ca K}(X,\Omega)} \ar[dl]^{G_{\tau,X}}
\save \POS?="codiso" \restore \\ & {\DFib(1,X)}
\POS "domiso"; "codiso" **@{}; ?*{\iso}} \]
for all $X$, and $G_{\tau,X}$ is fully faithful.
Thus $G_{\tau_f,X}$ is fully faithful iff ${\ca K}(X,f)$ is.
\end{proof}
\begin{example}\label{cdc-LTT}
Let $\ca E$ be an elementary topos with subobject classifier
denoted as $\tau_0:1{\rightarrow}\Omega_0$.
A Lawvere-Tierney topology on $\ca E$
amounts to an idempotent monad $j$ on $\Omega$
in $\Cat(\ca E)$, such that $j$ preserves $\wedge$.
As a monad in a finitely complete $2$-category
we can take its Eilenberg-Moore object \cite{Str76},
part of which is the forgetful arrow
$u^j:\Omega^j{\rightarrow}\Omega$.
Since $j$ is idempotent, $u^j$ is fully faithful
and by proposition(\ref{cdc-pb})
\[ \PbSq {\Omega^j_{\bullet}} 1 {\Omega}
{\Omega^j} {\tau_j} {} {\tau} {u^j} \]
$\tau_j$ is a classifying discrete opfibration,
and it is easily verified that ${\Omega^j_{\bullet}}$
is the terminal object
(by the same argument that establishes $\Omega_{\bullet}=1$
for subobject classifiers).
The image of $G_{\tau_j,X}$
consists of those cosieves $p:E{\rightarrow}X$
such that $p_0$ is a monomorphism which is $j$-dense
in the sense of topos theory.
In fact $\tau_{j0}:1{\rightarrow}\Omega^j_0$
is the subobject classifier for the topos $\Sh_j(\ca E)$
of sheaves for the topology $j$,
and thus by example(\ref{cdc-subobj}) $\tau_j$
is also a classifying discrete opfibration
for $\Cat(\Sh_j(\ca E))$.
This example could also have been obtained from example(\ref{cdc-subobj})
via sheafification and the following theorem.
\end{example}
\begin{theorem}\label{cdc-from-2adj}
Let
\[ \xymatrix{{\ca A} \ar@<-1.5ex>[rr]_-{S} \save \POS?="bot" \restore
&& {\ca B} \ar@<-1.5ex>[ll]_-{E} \save \POS?="top" \restore
\POS "top"; "bot" **@{}; ?(.5)*{\perp}} \]
be a $2$-adjunction, $\ca A$ and $\ca B$ be finitely complete $2$-categories,
and $\tau:\Omega_{\bullet}{\rightarrow}{\Omega}$ be a classifying
discrete opfibration in $\ca A$.
We write $\eta$ and $\varepsilon$ for the unit and counit of $E \ladj S$.
If
\begin{enumerate}
\item  $E$ preserves pullbacks, and
\item  the naturality square
\[ \xymatrix{{ES\Omega_{\bullet}}
\ar[r]^-{\varepsilon_{\Omega_{\bullet}}} \ar[d]_{ES\tau}
& {\Omega_{\bullet}} \ar[d]^{\tau} \\
{ES\Omega} \ar[r]_-{\varepsilon_{\Omega}} & {\Omega}} \]
is a pullback
\end{enumerate}
then $S\tau$ is a classifying discrete opfibration.
\end{theorem}
\begin{proof}
Since discrete opfibrations can be defined representably,
they are preserved by right $2$-adjoints, and so it suffices
to show that $S\tau$ is classifying.
We abuse notation and denote $G_{S\tau,B}$ simply by $G$, and
for $f:B{\rightarrow}S\Omega$ depict $Gf$ by
\[ \PbSq {e(f)} {S\Omega_{\bullet}} {S\Omega} B {G(f)} {} {S\tau} f \]
as in the discussion preceding definition(\ref{class-dcofib}).
Suppose that $f,g:B{\rightarrow}S\Omega$ and $\phi:e(f){\rightarrow}e(g)$
such that $G(g)\phi=G(f)$.
We must show that there is a unique $\phi':f{\rightarrow}g$ such that
$G(\phi')=\phi$.
\[ \xymatrix{& {Ee(f)}
\ar[dl]_{E\phi} \ar[dddl]|(.65){EGf} \ar[dr]^{Eq} \\
{Ee(g)} \ar[dd]_{EGg} \ar@/_{1pc}/[rr]^-{Ep}
&& {ES\Omega_{\bullet}} \ar[dd]|{ES\tau}
\ar[rr]^-{\varepsilon_{\Omega_{\bullet}}}
&& {\Omega_{\bullet}} \ar[dd]^{\tau} \\ \\
{SEB} \ar@/^{1pc}/[rr]^-{Ef} \save \POS?="domphi" \restore
\ar@/_{1pc}/[rr]_-{Eg} \save \POS?="codphi" \restore && {ES\Omega}
\ar[rr]_-{\varepsilon_{\Omega}} && {\Omega}
\POS "domphi"; "codphi" **@{}; ?(.3) \ar@{=>}^{E\phi'} ?(.7)} \]
Since $E$ preserves pullbacks,
the naturality square in the above diagram is a pullback
and $\tau$ is a classifying discrete opfibration,
$E\phi$ induces a unique
$\phi_2:\varepsilon_{\Omega}Ef{\rightarrow}\varepsilon_{\Omega}Eg$
such that $G_{\tau,EB}(\phi_2)=E\phi$.
By adjointness there is a unique $\phi'$ such that
$\varepsilon_{\Omega}E\phi'=\phi_2$.
Now we must see that $G(\phi')=\phi$.
Let $\phi_3:\varepsilon_{\Omega_{\bullet}}Eq{\rightarrow}
\varepsilon_{\Omega_{\bullet}}E(p\phi)$
be the unique lifting of $\phi_2$ through $\tau$.
By adjointness there is a unique
$\phi_4:q{\rightarrow}p\phi$ such that $E\phi_4=\phi_3$,
and by the uniqueness clause of the
$2$-dimensional universal property for the
naturality pullback square, the diagram
\[ \xymatrix{& {Ee(f)} \save \POS="domphib" \restore
\ar[dl]_{E\phi} \ar[dddl]|(.65){EGf} \ar[dr]^{Eq} \\
{Ee(g)} \ar[dd]_{EGg} \ar@/_{1pc}/[rr]^-{Ep}
\save \POS?="codphib" \restore && {ES\Omega_{\bullet}} \ar[dd]^{ES\tau} \\ \\
{EB} \ar@/^{1pc}/[rr]^-{Ef} \save \POS?="domphi" \restore
\ar@/_{1pc}/[rr]_-{Eg} \save \POS?="codphi" \restore && {ES\Omega}
\POS "domphi"; "codphi" **@{}; ?(.3) \ar@{=>}^{E\phi'} ?(.7)
\POS "domphib"; "codphib" **@{}; ?(.5) \ar@{=>}^{E\phi_4} ?(.7)} \]
commutes, and so by adjointness
\[ \xymatrix{& {e(f)} \save \POS="domphib" \restore
\ar[dl]_{\phi} \ar[dddl]|(.65){Gf} \ar[dr]^{q} \\
{e(g)} \ar[dd]_{Gg} \ar@/_{1pc}/[rr]^-{p}
\save \POS?="codphib" \restore && {S\Omega_{\bullet}} \ar[dd]^{S\tau} \\ \\
{B} \ar@/^{1pc}/[rr]^-{f} \save \POS?="domphi" \restore
\ar@/_{1pc}/[rr]_-{g} \save \POS?="codphi" \restore && {S\Omega}
\POS "domphi"; "codphi" **@{}; ?(.3) \ar@{=>}^{\phi'} ?(.7)
\POS "domphib"; "codphib" **@{}; ?(.5) \ar@{=>}^{\phi_4} ?(.7)} \]
commutes also.
To be more precise, paste the naturality square for $\varepsilon$
onto the previous commuting diagram,
apply $S$ to the entire resulting diagram,
and then precompose this with $\eta_{e(f)}$; the result of all this,
because of the $2$-naturality of $\eta$ and adjunction equations,
is the commutative diagram just obtained.
That is, $G(f)\phi'=S(\tau)\phi_4$,
and so $\phi_4$ is the unique lifting of $\phi'$ through $S\tau$,
so that $G\phi'=\phi$.
To see that $\phi'$ is unique
suppose that $\psi:f{\rightarrow}g$ such that $G\psi=\phi$.
Then both $\varepsilon_{\Omega}E\phi'$ and $\varepsilon_{\Omega}E\psi$
classify $E\phi$, and so they are equal.
By adjointness one obtains $\phi'=\psi$ as required.
\end{proof}
\begin{example}\label{cdc-catpshf}
Let $\C$ be a small category
and consider the $2$-adjunction
\[ \xymatrix{{\CAT} \ar@<-1.5ex>[rr]_-{\Sp {\C}}
\save \POS?(.46)="bot" \restore
&& {\PSH {\C}} \ar@<-1.5ex>[ll]_-{E}
\save \POS?(.54)="top" \restore
\POS "top"; "bot" **@{}; ?(.5)*{\perp}} \]
given by
\[ \begin{array}{lcr}
E(X) = \op {\el(\op X)} && \Sp {\C}(Z) = [\op {(\C/C)},Z].
\end{array} \]
We shall now see that $\Sp {\C}$ preserves
classifying discrete opfibrations.
By theorem(\ref{cdc-from-2adj}) it suffices to show
that if $p:A{\rightarrow}B$ is a discrete opfibration,
then the naturality square
\[ \xymatrix{{E{\Sp {\C}}A} \ar[d]_{E{\Sp {\C}}p} \ar[r]^-{\varepsilon_A} & A \ar[d]^{p} \\
{E{\Sp {\C}}B} \ar[r]_-{\varepsilon_B} & B} \]
is a pullback.
To carry out this verification
we need an explicit description of $\varepsilon$.
For $A \in \CAT$ the category $E{\Sp {\C}}A$ has:
\begin{itemize}
\item  objects: pairs $(C,a)$ where $C \in \C$ and
$a:{\op {(\C/C)}}{\rightarrow}A$.
\item  arrows: from $(C_1,a_1){\rightarrow}(C_2,a_2)$ are pairs
$(f,\overline{f})$, where $f:C_2{\rightarrow}C_1$ and
\[ \TriTwoCell {\op {(\C/C_1)}} {\op {(\C/C_2)}} A {a_1}
{\op {f^{!}}} {a_2} {\overline{f}} \]
\end{itemize}
Then $\varepsilon_A(C,a)=a(1_C)$ and $\varepsilon_A(f,\overline{f})$
is $\overline{f}_{1_C}$.
To see that the above naturality square is a pullback on objects,
let $a \in A$ and $b:\op {(\C/C)}{\rightarrow}B$ such that $b(1_C)=pa$.
We must show that there is a unique
$\overline{a}:\op {(\C/C)}{\rightarrow}A$ such that
$\overline{a}(1_C)=a$ and $f\overline{a}=b$.
For $\gamma:C'{\rightarrow}C$ we have
$b(t_{\gamma}):pa{\rightarrow}b(\gamma)$ and so a unique
$\overline{a}(t_{\gamma}):a{\rightarrow}\overline{a}(\gamma)$
such that $p\overline{a}(t_{\gamma})=b(t_{\gamma})$.
Functoriality and uniqueness for $\overline{a}$ follows from
the uniqueness of the liftings involved in its description.
To see that the above naturality square is a pullback on arrows,
let $\alpha:a_1{\rightarrow}a_2$ in A
and $(\beta,\overline{\beta}):(C_1,b_1){\rightarrow}(C_2,b_2)$
in $E{\Sp {\C}}B$ be such that $\overline{\beta}_{1_C}=p\alpha$.
We must show that there is a unique
$\overline{\alpha}:\overline{a_1}{\rightarrow}
\overline{a_2}\beta^{!}$ such that $p\overline{\alpha}=\overline{\beta}$.
For $\gamma:C{\rightarrow}C_1$ the square on the left
\[ \TwoDiagRel
{\xymatrix{{pa_1} \ar[r]^-{b_1(t_{\gamma})} \ar[d]_{p\alpha}
& b_1{\gamma} \ar[d]^{\overline{\beta}_{\gamma}} \\
{pa_2} \ar[r]_-{b_2(t_{\gamma})} & {b_2\beta\gamma}}}
{}
{\xymatrix{{a_1} \ar[r]^-{\overline{a_1}(t_{\gamma})} \ar[d]_{\alpha}
& {\overline{a_1}(\gamma)} \ar[d]^{\overline{\alpha}_{\gamma}} \\
{a_2} \ar[r]_-{\overline{a_2}(t_{\gamma})}
& {\overline{a_2}(\beta\gamma)}}} \]
is commutative in $B$ since it is the naturality square
for $\overline{\beta}$ at $t_{\gamma}$,
and $\overline{\alpha}_{\gamma}$
is defined as the unique map in $A$ such that
the square on the right is commutative in $A$
and $p\overline{\alpha}_{\gamma}=\overline{\beta}_{\gamma}$.
Naturality and uniqueness for $\overline{\alpha}$ follows from
the uniqueness of the liftings involved in its description.
\end{example}
\begin{example}\label{cdc-globcatI}
Consider the case $\C=\G$, where $\G$ is the category
presented as follows:
the objects of $\G$ are natural numbers,
for each $n \in \N$ there are generating maps
\[ \xymatrix{n \ar@<1ex>[r]^-{\sigma} \ar@<-1ex>[r]_-{\tau} & {n+1}}, \]
and these maps satisfy $\sigma\tau=\tau\tau$ and $\tau\sigma=\sigma\sigma$.
The category $\PSh {\G}$ is the category of \emph{globular sets}
and $\PSH {\G}$ is the $2$-category of \emph{globular categories}
in the sense of \cite{Bat98}.
A functor $\op {(\G/n)}{\rightarrow}A$ is called an $n$-span
in $A$, because in the case $n=1$ such a functor is a span in $A$:
the category $\op {(\G/1)}$ is the poset which as a category looks like
\[ {\bullet} \leftarrow {\bullet} \rightarrow {\bullet} \]
In fact all the ${\G/n}$ are posets.
For example, the posets $\op {(\G/2)}$ and $\op {(\G/3)}$
look like
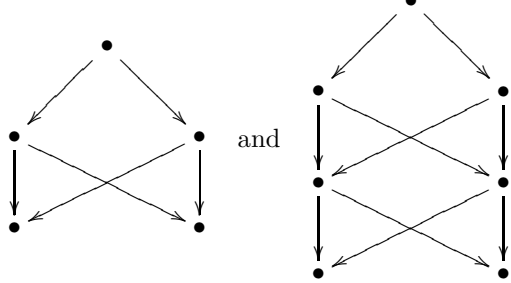
\[ \TwoDiagRel
{\xymatrix{& {\bullet} \ar[dr] \ar[dl] \\
{\bullet} \ar[d] \ar[drr] && {\bullet} \ar[d] \ar[dll] \\
{\bullet} && {\bullet}}}
{\textnormal{\,\,\, and \,\,\,}}
{\xymatrix{& {\bullet} \ar[dr] \ar[dl] \\
{\bullet} \ar[d] \ar[drr] && {\bullet} \ar[d] \ar[dll] \\
{\bullet} \ar[d] \ar[drr] && {\bullet} \ar[d] \ar[dll] \\
{\bullet} && {\bullet}}} \]
respectively.
The globular category ${\Sp {\C}}(\Set)$ plays a central role in \cite{Bat98},
it is the globular category of higher spans, and in \cite{Str00}
it is seen as the globular analogue of the category of sets.
In our language this last idea can be expressed as follows: by the previous
example ${\Sp {\C}}(U:\Set_{\bullet}{\rightarrow}\Set)$
is a classifying discrete opfibration.
Of course one may replace $U$ by any classifying discrete opfibration
in $\CAT$ (see example(\ref{cdc-cat}))
to get one in $\PSH {\G}$ in this way.
\end{example}
\begin{example}\label{cdc-globcatII}
The examples described here are likely to be important in the study
of the stabilisation hypothesis \cite{BD95} in the globular setting.
Consider the functor $(-+1):{\G}{\rightarrow}{\G}$
which acts as $n \mapsto n+1$ on objects, and whose
arrow map described by
\[ \TwoDiagRel
{\xymatrix{n \ar@<1ex>[r]^-{\sigma} \ar@<-1ex>[r]_-{\tau} & {n+1}}}
{\mapsto}
{\xymatrix{{n+1} \ar@<1ex>[r]^-{\sigma} \ar@<-1ex>[r]_-{\tau} & {n+2}}} \]
Restriction and right extension along $(-+1)$
provides an adjunction between endofunctors of $\PSh {\G}$,
and since these endofunctors are both right adjoints, they preserve
pullbacks and so may be regarded as acting on category objects
to provide a $2$-adjunction
\[ \xymatrix{{\PSH {\G}} \ar@<-1.5ex>[rr]_-{\Sigma}
\save \POS?="bot" \restore
&& {\PSH {\G}} \ar@<-1.5ex>[ll]_-{D}
\save \POS?="top" \restore
\POS "top"; "bot" **@{}; ?(.5)*{\perp}}. \]
Explicitly given a globular category $X$:
\[ \xymatrix{{X_0} & {X_1} \ar@<1ex>[l] \ar@<-1ex>[l]
& {X_2} \ar@<1ex>[l] \ar@<-1ex>[l]
& {X_3} \ar@<1ex>[l] \ar@<-1ex>[l]
& {...} \ar@<1ex>[l] \ar@<-1ex>[l]} \]
$DX$ is obtained by forgetting about $X_0$, and $\Sigma{X}$ is
\[ \xymatrix{{1} & {X_0} \ar@<1ex>[l] \ar@<-1ex>[l]
& {X_1} \ar@<1ex>[l] \ar@<-1ex>[l]
& {X_2} \ar@<1ex>[l] \ar@<-1ex>[l]
& {...} \ar@<1ex>[l] \ar@<-1ex>[l]} \]
The counit of this adjunction is an isomorphism,
and so $\Sigma$ preserves classifying discrete opfibrations
by theorem(\ref{cdc-from-2adj}).
Thus for any classifying discrete opfibration
$\tau:\Omega_{\bullet}{\rightarrow}{\Omega}$ in $\CAT$,
one obtains a sequence
\[ {\Sp {\G}}(\tau), {\Sigma}{\Sp {\G}}(\tau), {\Sigma}^2{\Sp {\G}}(\tau),
{\Sigma}^3{\Sp {\G}}(\tau), ... \]
of classifying discrete opfibrations in $\PSH {\G}$.
\end{example}
\begin{definition}
A \emph{2-topos} $(\ca K, \iop {(-)}, \tau)$
is a finitely complete cartesian closed
2-category $\ca K$ equipped with a duality involution $\iop {(-)}$
and a classifying discrete opfibration
$\tau:\Omega_{\bullet}{\rightarrow}{\Omega}$.
\end{definition}
We obtain our examples of 2-toposes from the above
examples of classifying discrete opfibrations,
since in each case the underlying $\ca K$ is of the form $\Cat(\ca E)$
for $\ca E$ a finitely complete category, and so comes
with a canonical duality involution by example(\ref{op-internal-CAT}).

\section{Yoneda structures from $2$-toposes}
\label{2-topos->yoneda}

Let $(\ca K, \iop {(-)}, \tau)$ be a 2-topos.
The purpose of this section is to
exhibit a good yoneda structure on $\ca K$.
This construction is due to Ross Street
\cite{Str74b} \cite{Str80}.

For $A$ in $\ca K$ we define
\[ \PSh A = [\iop A, \Omega] \]
and so we have a $2$-functor
$(\PSh {-}):\coop {\ca K}{\rightarrow}{\ca K}$.
For $A, B$ in $\ca K$ we denote by $G_{A,B}$ the composite
\[ \xymatrix{{\ca K(B,\PSh A){\iso}\ca K({\iop A} \times B,\Omega)}
\ar[rr]^-{G_{\tau,{\iop A \times B}}}
&& \DFib(1,\iop A \times B) \ar[rr]^{d_{1,A,B}} && \DFib(A,B)} \]
Clearly the $G_{A,B}$ are fully faithful
and pseudo-natural in $A$ and $B$.
A span
\[ S : A \leftarrow E \rightarrow B \]
is an \emph{attribute} when it is isomorphic to a span in the image
of $G_{A,B}$. For $A \in \ca K$ we denote by
\[ \epsilon_A : \xymatrix{A & {\PSh A_{\bullet}} \ar[l]_-{\pi_A}
\ar[r]^-{u_A} & {\PSh A}} \]
the attribute $G_{A,\PSh A}(1_{\PSh A})$.
\begin{lemma}\label{att-nat}
\begin{enumerate}
\item  For any attribute $S : A \leftarrow E \rightarrow B$
there is an arrow $f:B{\rightarrow}\PSh A$
unique up to isomorphism such that
$\rev f \comp \epsilon_A \iso S$.
\label{natB}
\item  Let $\alpha:A_1{\rightarrow}A_2$,
then $\rev {\PSh {\alpha}} \comp \epsilon_{A_1} \iso
\epsilon_{A_2} \comp \alpha$.
\label{natA}
\end{enumerate}
\end{lemma}
\begin{proof}
(\ref{natB}):
If $S$ is an attribute then by definition there is an
$f:B{\rightarrow}\PSh A$ such that $G_{A,B}(f){\iso}S$.
By pseudo-naturality of the $G_{A,B}$ there is an isomorphism
\[ \xymatrix{{\ca K(\PSh A,\PSh A)} \ar[d]_{-{\comp}f}
\save \POS?="domiso" \restore
\ar[r]^-{G_{A,\PSh A}} & {\DFib(A,\PSh A)} \ar[d]^{{\rev f}{\comp}-}
\save \POS?="codiso" \restore \\
{\ca K(B,\PSh A)} \ar[r]_-{G_{A,B}} & {\DFib(A,B)}
\POS "domiso"; "codiso" **@{}; ?*{\iso}} \]
whose component at $1_{\PSh A}$ is an isomorphism
$G_{A,B}(f){\iso}\rev f \comp \epsilon_A$.
If $g:B{\rightarrow}\PSh A$ also satisfies
$\rev g \comp \epsilon_A \iso S$,
then we have $G_{A,B}(f){\iso}G_{A,B}(g){\iso}S$
and so $f{\iso}g$ by the fully faithfulness of $G_{A,B}$.\\
(\ref{natA}):
The desired isomorphism arises by (\ref{natB})
and the component of $1_{\PSh A_2}$ of
\[ \xymatrix{{\ca K(\PSh A_2,\PSh A_2)} \ar[d]_{{\PSh {\alpha}}{\comp}-}
\save \POS?="domiso" \restore
\ar[r]^-{G_{A_2,\PSh A_2}} & {\DFib(A_2,\PSh A_2)} \ar[d]^{-{\comp}\alpha}
\save \POS?="codiso" \restore \\
{\ca K(\PSh A_2,\PSh A_1)} \ar[r]_-{G_{A_1,\PSh A_2}}
& {\DFib(A_1,\PSh A_2)}
\POS "domiso"; "codiso" **@{}; ?*{\iso}} \]
\end{proof}
Another way to express lemma(\ref{att-nat})(\ref{natB})
which is worth keeping in mind is that $G_{A,B}$ ``is''
span composition with $\epsilon_A$. More precisely we have
$G_{A,B} \iso \rev {(-)} \comp \epsilon_A$.
Define a map $f:A{\rightarrow}B$ to be admissible
when the span $f/B$ is an attribute.
Thus there is an arrow $B(f,1):B{\rightarrow}\PSh A$
unique up to isomorphism with the property that
$f/B \iso \rev {B(f,1)} \comp \epsilon_A$.
The admissible arrows form a right ideal: given $f:A{\rightarrow}B$
admissible and $g:C{\rightarrow}B$ we have the following
isomorphisms of spans
\begin{eqnarray*}
(fg)/B &\iso& f/B \comp g \iso \rev {B(f,1)} \comp \epsilon_A \comp g \\
&\iso& \rev {B(f,1)} \comp \rev {\PSh g} \comp \epsilon_C
\iso \rev {({\PSh g}B(f,1))} \comp \epsilon_C
\end{eqnarray*}
and so $fg$ is admissible and $B(fg,1) \iso {\PSh g}B(f,1)$.
In particular for an admissible object $A$ we denote $A(1,1)$
as $y_A : A{\rightarrow}{\PSh A}$.
\begin{proposition}\label{att-yoneda}
If $S \in \Span(A,B)$ is an attribute, $A$ is an admissible object
and $f:B{\rightarrow}\PSh A$ such that
$S \iso \rev f \comp \epsilon_A$;
then $S \iso y_A/f$.
\end{proposition}
\begin{proof}
By the definition of lax pullback
it suffices to prove that for any span
$S' : \xymatrix@1{A & E \ar[l]_-{d} \ar[r]^-{c} & B}$
there is a bijection between maps $S{\rightarrow}S'$ of spans
and 2-cells $y_Ad{\rightarrow}fc$, and that this bijection
is natural in $E$.
Noting that $S'{\iso}c{\comp}{\rev d}$ and $c \ladj \rev c$,
maps $S{\rightarrow}S'$ are in bijection with
maps $\rev d{\rightarrow}{\rev c}{\comp}S$.
By the coyoneda lemma (theorem(\ref{yon}))
such maps are in bijection with
$A/d{\rightarrow}{\rev c}{\comp}S$ in $\Span(A,B)$.
Now
\[ A/d \iso \rev d \comp A/A
\iso \rev d \comp \rev {y_A} \comp \epsilon_A
\iso \rev {(y_Ad)}  \comp \epsilon_A \]
and
\[ \rev c \comp S \iso \rev c \rev f \comp \epsilon_A
\iso \rev {(fc)} \comp \epsilon_A  \]
so we have a bijection between span maps $S'{\rightarrow}S$
and span maps
${\rev {(y_Ad)}}{\comp}{\epsilon_A}{\rightarrow}
{\rev {(fc)}}{\comp}{\epsilon_A}$.
Since $G_{A,B} \iso \rev {(-)} \comp \epsilon_A$
is fully faithful, the result follows.
\end{proof}
Thus when $A$ is admissible, the map $y_A$ is also
admissible: by proposition(\ref{att-yoneda}) we have
$\epsilon_A \iso y_A/{\PSh A}$ which is an attribute by definition,
and so we have $\PSh A(y_A,1) \iso 1_{\PSh A}$.

To complete the specification of the data for a yoneda structure
we must define for $f:A{\rightarrow}B$ with $A$ and $f$ admissible,
a $2$-cell $\chi^f:y_A{\rightarrow}B(f,1)f$.
For such an $f$ write $\varphi_f$ for the isomorphism
$\rev B(f,1) \comp \epsilon_A \iso f/B$.
Define $h_f:A/A{\rightarrow}f/f$ as the unique arrow such that
\[ \TwoDiagRel
{\xymatrix{{A/A} \ar@/_{1pc}/[ddr]
\save \POS?="domeq1" \restore
\ar[dr]|{h_f} \save \POS?="domeq2" \restore
\ar@/^{1pc}/[drr]
\save \POS?="codeq2" \restore
\\ & {f/f} \save \POS="codeq1" \restore \ar[d]
\save \POS?="domlam" \restore
\ar[r] & A \ar[d]^{f}
\save \POS?="codlam" \restore
\\ & A \ar[r]_{f} & B
\POS "domeq1"; "codeq1" **@{}; ?*{=}
\POS "domeq2"; "codeq2" **@{}; ?*{=}
\POS "domlam"; "codlam" **@{}; ?(.35) \ar@{=>} ?(.65)}}
{=}
{\xymatrix{{A/A} \ar[d]
\save \POS?="domlam" \restore
\ar[r] & A \ar[d]|{1_A}
\save \POS?="codlam" \restore
\ar@/^{1pc}/[ddr]^{f}
\save \POS?="codeq2" \restore
\\ A \ar@/_{1pc}/[drr]_{f}
\save \POS?="domeq1" \restore
\ar[r]|-{1_A} & A  \save \POS="domeq2" \restore \ar[dr]|{f}
\save \POS?="codeq1" \restore
\\ && B
\POS "domeq1"; "codeq1" **@{}; ?*{=}
\POS "domeq2"; "codeq2" **@{}; ?*{=}
\POS "domlam"; "codlam" **@{}; ?(.35) \ar@{=>} ?(.65)}} \]
where the 2-cells in the above diagram are lax pullback cells.
Then $\chi^f$
is defined as the unique $2$-cell such that
\begin{equation}\label{def-chi}
\xymatrix{{\rev {y_A} \comp \epsilon_A}
\ar[rr]^-{\rev {(\chi^f)} \comp \epsilon_A}
\ar[d]_{\varphi_{1_{\PSh A}}} &&
{\rev {(B(f,1)f)} \comp \epsilon_A
\iso \rev f \comp \rev {B(f,1)} \comp \epsilon_A}
\ar[d]^{\rev f \comp \varphi_{f}} \\ {A/A} \ar[rr]_{h_f} && {f/f}}
\end{equation}
commutes.
\begin{theorem}\label{yon-ax}
\cite{Str74b}
Let $A$ and $f:A{\rightarrow}B$ be admissible, and define $\chi^f$
as above.
\begin{enumerate}
\item  $\chi^f$ exhibits $f$ as an absolute left lifting of $y_A$
along $B(f,1)$.
\label{yon-lift}
\item  If
\[ \TriTwoCell A B {\PSh A} {y_A} f g {\phi} \]
which exhibits $f$ as an absolute left lifting of $y_A$ along $g$;
then $\phi$ exhibits $g$ as a pointwise left extension of $y_A$ along $f$.
\label{yon-ext}
\end{enumerate}
\end{theorem}
\begin{proof}
(\ref{yon-lift}):
For any span $\xymatrix@1{A & X \ar[l]_-{a} \ar[r]^-{b} & B}$
we must exhibit a bijection between $2$-cells
$y_Aa{\rightarrow}B(f,1)b$ and $2$-cells $fa{\rightarrow}b$
natural in $X$ which identifies $\chi^f$ and $1_f$.
Applying $G_{X,A}$ gives a bijection between $2$-cells
$y_Aa{\rightarrow}B(f,1)b$ and maps $A/a{\rightarrow}f/b$
of spans. By the coyoneda lemma such span maps are in bijection
with spans maps $\rev a{\rightarrow}{f/b}$ and thus with
$2$-cells $fa{\rightarrow}b$ by the definition of $f/b$.
This bijection is clearly natural in $X$, and tracing
$\chi^f$ through these bijections clearly produces $1_f$.
\\ (\ref{yon-ext}):
The proof proceeds in two stages: first we show that the statement is
true for $\phi=\chi^f$, and then in the general case we show that
$g \iso B(f,1)$.
For $b:X{\rightarrow}B$ we must show that pasting with
\begin{equation}\label{chilam}
\xymatrix{{f/b} \ar[rr]^-{q} \ar[d]_{p} \save \POS?="domlam" \restore
&& X \ar[d]^{b} \save \POS?="codlam" \restore \\
A \ar[rr]|{f} \ar[dr]_{y_A} \save \POS?="domphi" \restore
&& B \ar[dl]^{B(f,1)} \save \POS?="codphi" \restore \\ & {\PSh A}
\POS "domlam"; "codlam" **@{} ?(.42) \ar@{=>}^{\lambda} ?(.58)
\POS "domphi"; "codphi" **@{} ?(.35) \ar@{=>}^{\chi^f} ?(.65)}
\end{equation}
gives a bijection between 2-cells $B(f,1)b{\rightarrow}h$
and 2-cells $y_Ap{\rightarrow}hq$.
It suffices to exhibit a bijection between such $2$-cells
which is natural in $X$ and which sends $1_{B(f,1)b}$ to the above
composite.
Applying $G_{X,A}$ and noting that
$\rev {(B(f,1)b)} \comp \epsilon_A \iso q \comp \rev p$
exhibits a bijection between $2$-cells $B(f,1)b{\rightarrow}h$
and span maps $q{\comp}{\rev p}{\rightarrow}{\rev h}{\comp}{\epsilon_A}$,
and these are in bijection with span maps
${\rev p}{\rightarrow}{\rev q}{\comp}{\rev h}{\comp}{\epsilon_A}$
since $q \ladj \rev q$.
By the coyoneda lemma and since $\rev q \comp \rev h \iso \rev {(hq)}$,
these are in turn in bijection with span maps
$A/p{\rightarrow}{\rev {(hq)}}{\comp}{\epsilon_A}$.
Since $A/p \iso \rev {(y_Ap)}{\comp}{\epsilon_A}$,
$G_{A,B} \iso \rev {(-)} \comp \epsilon_A$
gives a bijection between
such maps
and with $2$-cells $y_Ap{\rightarrow}hq$ as required.
These correspondences are clearly natural
and map $1_{B(f,1)b}$ to the composite in (\ref{chilam}).
For the general case it now suffices to show
$g \comp \epsilon_A \iso f/B$.
It suffices to exhibit,
for each span
$S : \xymatrix@1{A & E \ar[l]_-{d} \ar[r]^-{c} & B}$,
a bijection between span maps $S{\rightarrow}g{\comp}{\epsilon_A}$
and $2$-cells $fd{\rightarrow}c$ natural in $E$.
Composition $\phi$ gives a bijection between $2$-cells $fd{\rightarrow}c$
and $2$-cells $y_Ad{\rightarrow}gc$ since $\phi$ is an absolute
left lift. Applying $\rev {(-)} \comp \epsilon_A$
and noting that $A/d \iso \rev {(y_Ad)} \comp \epsilon_A$
gives a bijection between $2$-cells $y_Ad{\rightarrow}gc$
and maps of spans
$A/d{\rightarrow}{\rev c}{\comp}{\rev g}{\comp}{\epsilon_A}$,
which by $c \ladj \rev c$ correspond to maps of spans
$S{\rightarrow}{\rev g}{\comp}{\epsilon_A}$
since $S \iso c \comp \rev d$.
\end{proof}
This last result ensures that the axioms for the intended yoneda structure
are indeed satisfied. We have proved:
\begin{corollary}\label{thm:2topos->yon}
Let $(\ca K, \iop {(-)}, \tau)$ be a 2-topos.
Then
\begin{itemize}
\item  $\PSh A = [\iop A,\Omega]$.
\item  $f:A{\rightarrow}B$ is admissible when $f/B$ is an attribute.
\item  When $A$ is admissible $y_A=A(1,1)$.
\item  For $f:A{\rightarrow}B$ such that $A$ and $f$ are admissible
define $\chi^f:y_A{\rightarrow}B(f,1)f$ as in (\ref{def-chi}) above.
\end{itemize}
specifies a good yoneda structure for $\ca K$.
\end{corollary}
%

\section{Characterising admissible maps}
\label{sec:char-adm}

When working with an example one would usually like
to characterise what the admissible maps are.
The following theorem is useful for this task.
Given a 2-topos with classifying discrete opfibration $\tau$,
we shall say that a map $f:A{\rightarrow}B$ is \emph{$\tau$-admissible}
when it is admissible in the sense of the yoneda structure on $\ca K$
induced by corollary(\ref{thm:2topos->yon}).
Similarly an object $A \in \ca K$ is \emph{$\tau$-small}
when it is small in the sense of this induced yoneda structure.
We shall also say that a map $f:A{\rightarrow}B$ in a 2-category
factors \emph{essentially through} a map $g:C{\rightarrow}B$
when there is a map $h:A{\rightarrow}C$ and an isomorphism
$f{\iso}gh$.
\begin{theorem}\label{prop:univ-enlarge}
Let $\ca K$ be a finitely complete cartesian closed 2-category
with duality involution $\iop {(-)}$.
Suppose that $\tau$ and $\tau'$ are classifying discrete opfibrations
fitting in a pullback square:
\[ \PbSq {\Omega_{\bullet}} {\Omega'_{\bullet}} {\Omega'}
{\Omega} {\tau} {} {\tau'} i \]
(see proposition(\ref{cdc-pb})).
For $f:A{\rightarrow}B$ we use the notation introduced above:
$\PSh A$, ${\PSh A}_{\bullet}$, $\epsilon_A$ and (when $f$
is $\tau$-admissible) $B(f,1)$; for these aspects of the 2-topos
$(\ca K, \iop {(-)}, \tau)$.
We use the notation
$\PSh A'$, ${\PSh A'}_{\bullet}$, $\epsilon'_A$ and (when $f$
is $\tau'$-admissible) $B'(f,1)$; for the corresponding
aspects of the 2-topos $(\ca K, \iop {(-)}, \tau')$.
If $f:A{\rightarrow}B$ is a $\tau'$-admissible map then the following
statements are equivalent.
\begin{enumerate}
\item  The adjoint transpose $\iop A \times B \rightarrow \Omega'$
of $B'(f,1)$ factors essentially through $i$.
\label{un-enl1}
\item  $B'(f,1)$ factors essentially through $[\iop A,i]$.
\label{un-enl2}
\item  $f$ is $\tau$-admissible.
\label{en-enl3}
\end{enumerate}
\end{theorem}
\begin{proof}
(\ref{un-enl1})$\iff$(\ref{un-enl2}) by the naturality of $\ca K$'s
cartesian closed structure.
In the diagram
\[ \xymatrix{
{{\ca K}([\iop A,\Omega'],[\iop A,\Omega'])}
\ar[r]^-{{\ca K}([\iop A,i],1)}
\ar[d] \save \POS?="domnat1" \restore &
{{\ca K}([\iop A,\Omega],[\iop A,\Omega'])}
\ar[d] \save \POS?="codnat1" \restore \save \POS?="codnat2" \restore &
{{\ca K}([\iop A,\Omega],[\iop A,\Omega])}
\ar[l]_-{{\ca K}(1,[\iop A,i])}
\ar[d] \save \POS?="domnat2" \restore \\
{{\ca K}(\iop A \times [\iop A,\Omega'],\Omega')} \ar[r]
\ar[d]_{G_{\tau'}} \save \POS?="dompsnat1" \restore &
{{\ca K}(\iop A \times [\iop A,\Omega],\Omega')}
\ar[dr]_{G_{\tau'}} \save \POS?="codpsnat1" \restore
\save \POS?="dombydef" \restore &
{{\ca K}(\iop A \times [\iop A,\Omega],\Omega)} \ar[l]
\ar[d]^{G_{\tau}} \save \POS?="codbydef" \restore \\
{\DFib(1,\iop A \times [\iop A,\Omega'])}
\ar[rr]^-{\rev {(\iop A \times [\iop A,i])} \comp (-)}
\ar[d] \save \POS?="dompsnat2" \restore &&
{\DFib(1,\iop A \times [\iop A,\Omega])}
\ar[d] \save \POS?="codpsnat2" \restore \\
{\DFib(A,[\iop A,\Omega'])}
\ar[rr]_-{\rev {[\iop A,i]} \comp (-)}
&& {\DFib(A,[\iop A,\Omega])}
\POS "domnat1"; "codnat1" **@{}; ?*{=}
\POS "domnat2"; "codnat2" **@{}; ?*{=}
\POS "dompsnat1"; "codpsnat1" **@{}; ?*{\iso}
\POS "dompsnat2"; "codpsnat2" **@{}; ?*{\iso}
\POS "dombydef"; "codbydef" **@{}; ?*{\iso}} \]
the top row of vertical arrows are the isomorphisms coming from
$\ca K$'s cartesian closed structure,
and the bottom row of vertical arrows are the equivalences
coming from the duality involution.
The top two squares commute by the naturality of the cartesian closed
structure. The triangle in the middle commutes up to isomorphism
by the definition of $G_{\tau'}$ and $G_{\tau}$
(we saw this also in proposition(\ref{cdc-pb})).
The other isomorphisms are pseudo naturality isomorphisms.
Tracing $1_{\PSh {A'}}$ and $1_{\PSh {A}}$ through this diagram,
and recalling the definitions of $\epsilon'_A$ and $\epsilon_A$,
we obtain an isomorphism of spans
\[ \epsilon_A \iso \rev {[\iop A,i]} \comp \epsilon'_A. \]
Suppose that $B'(f,1)$ factors essentially through $[\iop A,i]$,
so that there is a map $h:B{\rightarrow}{\PSh A}$ such that
$B'(f,1){\iso}[\iop A,i]h$. Then
\[ f/B \iso \rev h \comp \rev {[\iop A,i]} \comp \epsilon'_A
\iso f/B \iso \rev h \comp \epsilon_A \]
whence $f$ is $\tau$ admissible and $h \iso B(f,1)$.
Conversely if $f:A{\rightarrow}B$ is $\tau$-admissible
then by definition there is $B(f,1):B{\rightarrow}{\PSh A}$
such that $f/B \iso \rev {B(f,1)} \comp \epsilon_A$.
Thus $[\iop A,i]B(f,1) \iso B'(f,1)$ by lemma(\ref{att-nat}).
\end{proof}
\begin{example}\label{CAT-yonII}
Taking $\CAT$ with its usual duality involution,
and $\tau$ to be the forgetful functor
$U:\Set_{\bullet}{\rightarrow}{\Set}$,
corollary(\ref{thm:2topos->yon}) produces the yoneda structure
foreshadowed in example(\ref{CAT-yon}).
For $A \in \CAT$, ${\PSh A}_{\bullet}$ is the following
category of ``pointed presheaves'':
\begin{itemize}
\item  objects: are triples $(x,F,a)$ where $F \in {\PSh A}$,
$a \in A$, and $x \in F(a)$.
\item  arrows: an arrow $(x_1,F_1,a_1){\rightarrow}(x_2,F_2,a_2)$
is a pair $(\phi,\alpha)$
where $\phi:F_1{\rightarrow}F_2$ in $\PSh A$
and $\alpha:a_1{\rightarrow}a_2$ in $A$
such that $\phi_a(x_1)=F_2(\alpha)(x_2)$.
\end{itemize}
and
\[ \xymatrix{A & {{\PSh A}_{\bullet}} \ar[l]_-{\pi_A} \ar[r]^-{U_A}
& {\PSh A}} \]
are the obvious forgetful functors.
We stated in example(\ref{CAT-yon}) that $f:A{\rightarrow}B$ is admissible
iff $B(fa,b)$ is a small set for all $a \in A$ and $b \in B$.
To see this is the case let $\CAT'$ be a cartesian closed
2-category of categories which contains $\SET$ as an object
and $\CAT$ as a sub-2-category.
In the situation of theorem(\ref{prop:univ-enlarge})
with $\ca K=\CAT'$, $\Omega'=\SET$, $\Omega=\Set$
and $i:\Set{\rightarrow}{\SET}$ the inclusion,
notice that any $f:A{\rightarrow}B$ in $\CAT$
is $\tau'$-admissible because $B'(f,1)$ can be defined as the adjoint
transpose of
\[ \op A \times B \rightarrow \SET \]
given on objects by $(a,b){\mapsto}B(fa,b)$.
It is straight forward to observe directly that one has
an isomorphism of spans $f/B \iso \rev {B(f,1)} \comp \epsilon'_A$
and that there is a pullback square
\[ \xymatrix{{{\PSh A}_{\bullet}} \ar[r] \ar[d]_{U_A}
& {{\PSh A}'_{\bullet}} \ar[d]^{U'_A} \\
{\PSh A} \ar[r]_-{[\op A,i]} & {{\PSh A}'}} \]
in $\CAT'$.
Thus by theorem(\ref{prop:univ-enlarge})
$f$ is admissible iff $\forall a,b$, $B(fa,b)$
is in $\Set$ as claimed in example(\ref{CAT-yon}).
In particular $A \in \CAT$ is admissible iff its hom-sets are small.
Clearly a category $A$ equivalent to one with a small set of arrows
will be small in our sense: $A$ and $\PSh A$ are admissible.
The converse, that $A$ and $\PSh A$ admissible implies that $A$
is equivalent to a category which has a small set of arrows,
is a result of Freyd and Street \cite{FryStr}.
\end{example}
\begin{remark}\label{rem:size}
The results of \cite{FryStr}
apply also when ``small'' is taken to mean finite,
or more generally ``of cardinality less than $\lambda$''
where $\lambda$ is a regular cardinal.
Thus throughout this work, one could
replace $\Set$ by the category of sets of cardinality
less than $\lambda$.
\end{remark}
\begin{example}\label{psh-not-complete}
Here is an example of a 2-topos whose presheaves are \emph{not}
cocomplete.
Consider the classifying discrete opfibration obtained from
$U:\Set_{\bullet}{\rightarrow}\Set$ by pulling it back along
the inclusion of the full subcategory of $\Set$
consisting of the non-empty sets.
By theorem(\ref{prop:univ-enlarge})
a category is admissible for the resulting yoneda structure on $\CAT$
when it's hom-sets are non-empty and small.
The empty category $0$ is admissible (vacuously), and so
is $1 \iso \PSh 0$, so that $0$ is small.
Thus initial objects are small colimits for this yoneda structure.
However the category of non-empty sets lacks an initial object.
\end{example}
\begin{example}\label{CATPSH-yon}
In this example we characterise admissible maps and small objects
for example(\ref{cdc-catpshf}) where $\ca K=\PSH {\C}$
and $\Omega={\Sp {\C}}(\Set)$.
Unpacking the definitions and the duality involution
one obtains that for $A \in \PSH {\C}$ and $C \in \C$,
${\PSh A}_{\bullet}(C)$ is the category:
\begin{itemize}
\item  objects: are triples $(x,F,a)$ where $F \in {\PSh A}(C)$,
$a \in AC$, and $x \in F(a)(1_C)$.
\item  arrows: an arrow $(x_1,F_1,a_1){\rightarrow}(x_2,F_2,a_2)$
is a pair $(\phi,\alpha)$
where $\phi:F_1{\rightarrow}F_2$ in ${\PSh A}(C)$
and $\alpha:a_1{\rightarrow}a_2$ in $AC$
such that $\phi_{a,1_C}(x_1)=F_2(\alpha)_{1_C}(x_2)$.
\end{itemize}
To characterise the admissible maps we use
theorem(\ref{prop:univ-enlarge}) with $\ca K=[\op C,\CAT']$,
$\Omega'={\Sp {\C}}(\SET)$ via a version of example(\ref{cdc-catpshf})
with $\CAT'$ in place of $\CAT$,
and $\Omega={\Sp {\C}}(\Set)$.
As before notice that any $f:A{\rightarrow}B$ in $\PSH {\C}$
is $\tau'$-admissible,
and that the pullback square relating $\tau$ and $\tau'$ can be
obtained by applying $\Sp {\C}$ to the corresponding
pullback square found in the previous example.
Thus by theorem(\ref{prop:univ-enlarge})
it follows that $f$ is admissible iff
$\forall C \in {\C}, a \in AC, b \in BC$,
$B(C)(fa,b)$ is a small set.
Thus in particular $A \in \PSH {\C}$ is admissible
iff $\forall C$, $A(C)$ has small hom-sets.
If for all $C$, $A(C)$ is equivalent to a category with
a small set of arrows,
then since $\C$ is small, $A$ is small
in the sense that $A$ and $\PSh A$ are admissible.
For the converse note that
\begin{eqnarray*}
{\PSh A}(C) &{\iso}& {\PSH {\C}}(C,[\iop A,\Sp {\C}(\Set)])
\iso {\PSH {\C}}(C \times \iop A,\Sp {\C}(\Set)) \\
&{\iso}& {\CAT}(\op {\el(C \times A)},\Set)
\end{eqnarray*}
and so if ${\PSh A}(C)$ is locally small for all $C$,
then by \cite{FryStr}
$\el(C \times A)$ is small for all $C$,
and this implies that $A(C)$ is small for all $C$.
\end{example}
\begin{example}\label{susp-GLOB-yon}
In this example we characterise the admissible maps
and small objects for example(\ref{cdc-globcatII})
where $\ca K=\PSH {\G}$ and $\Omega=\Sigma^k\Sp {\G}(\Set)$.
We shall now write $\tau:\Omega_{\bullet}{\rightarrow}\Omega$
for the classifying discrete opfibration $\Sp {\G}(U)$
in $\PSH {\G}$, where $U:\Set_{\bullet}{\rightarrow}\Set$
is the forgetful functor.
We saw that $\Sigma$ preserves classifying discrete opfibrations.
By the characterisation of the $\tau$-admissible arrows and
$\tau$-small objects given in the previous example,
$\Sigma$ clearly preserves $\tau$-admissible arrows and
$\tau$-small objects.
Noting that $\Sigma{y_1}$ is $\tau$-admissible and
writing $i_1$ for the map $\Sigma{\Omega}(\Sigma{y_1},1)$, we have
\[ \TriTwoCell 1 {\Sigma{\Omega}} {\Omega} {y_1} {\Sigma{y_1}}
{i_1} {\chi^{\Sigma{y_1}}} \]
exhibiting $i_1$ as a pointwise
left extension of $y_1$ along $\Sigma{y_1}$,
and $\Sigma{y_1}$ as an absolute left lifting
of $y_1$ along $i_1${\footnotemark{\footnotetext{
Since $1$ is $\tau$-admissible we saw in
example(\ref{one-not-admissible})
that $t_{\Omega} \ladj y_1$ and so $y_1$ is a fully faithful
right adjoint. Fully faithful right adjoints are preserved by all
2-functors and so $\Sigma{y_1}$ is fully faithful,
whence $\chi^{\Sigma{y_1}}$ is in fact invertible
by proposition(\ref{ple-ff}).}}}.
By proposition(\ref{att-yoneda}) there is a 2-cell
\[ \LaxSq {\Omega_{\bullet}} {\Omega} 1 {\Omega} {} {\tau} {1} {y_1}
{\lambda} \]
which exhibits $\Omega_{\bullet}$ as $y_1/{\Omega}$ and since
$\Sigma$ is a right adjoint it preserves this lax pullback.
Define $\lambda'$ as follows:
\[ \TwoDiagRel
{\LaxSq {\Sigma\Omega_{\bullet}} {\Sigma\Omega} 1 {\Omega} {}
{\Sigma\tau} {i_1} {y_1} {\lambda'}}
{=}
{\xymatrix{{\Sigma\Omega_{\bullet}} \ar[rr]^-{\Sigma\tau}
\ar[d] \save \POS?="domslam" \restore && {\Sigma\Omega}
\ar[d]^{1} \save \POS?="codslam" \restore \\
1 \ar[rr]^-{\Sigma{y_1}} \ar[dr]_{y_1} \save \POS?="domchi" \restore
&& {\Sigma\Omega} \ar[dl]^{i_1} \save \POS?="codchi" \restore \\ & {\Omega}
\POS "domslam"; "codslam" **@{}; ?(.42) \ar@{=>}^{\Sigma\lambda} ?(.58)
\POS "domchi"; "codchi" **@{};
?(.35) \ar@{=>}^{\chi^{\Sigma{y_1}}} ?(.65)}} \]
Since $\Sigma\lambda$ is a lax pullback and $\chi^{\Sigma{y_1}}$
is an absolute left lifting,
it follows easily that $\lambda'$ exhibits $\Sigma\Omega_{\bullet}$
as $y_1/i_1$. By the universal property for $\lambda$,
there is a unique $i_{1,\bullet}$ fitting in a commutative square
\begin{equation}\label{sig-om-pb}
\xymatrix{{\Sigma\Omega_{\bullet}} \ar[r]^-{i_{1,\bullet}}
\ar[d]_{\Sigma\tau} & {\Omega_{\bullet}} \ar[d]^{\tau} \\
{\Sigma\Omega} \ar[r]_-{i_1} & {\Omega}}
\end{equation}
such that $\lambda{i_{1,\bullet}}=\lambda'$, and since $\lambda'$
is a lax pullback,
(\ref{sig-om-pb}) is a pullback.
Denoting $\ca K = [\op {\G},\CAT']$, $\Omega'=\Sp {\G}(\SET)$
and $\tau'$ for the corresponding classifying
discrete opfibration in $\ca K$, we have
\[ \xymatrix{{\Sigma\Omega_{\bullet}} \ar[d]_{\Sigma\tau}
\save \POS?(.3)="lpb1" \restore \ar[r] \save \POS?(.3)="tpb1" \restore
& {\Omega_{\bullet}} \ar[d]|{\tau} \save \POS?(.3)="rpb1" \restore
\save \POS?(.3)="lpb2" \restore \ar[r] \save \POS?(.3)="tpb2" \restore
& {\Omega'_{\bullet}} \ar[d]^{\tau'} \save \POS?(.3)="rpb2" \restore \\
{\Sigma\Omega} \ar[r]_-{i_1} \save \POS?(.3)="bpb1" \restore
& {\Omega} \ar[r]_-{i} \save \POS?(.3)="bpb2" \restore
& {\Omega'}
\POS "rpb1"; "lpb1" **@{}; ?!{"bpb1";"tpb1"}="cpb1" **@{}; ? **@{-};
"tpb1"; "cpb1" **@{}; ? **@{-}
\POS "rpb2"; "lpb2" **@{}; ?!{"bpb2";"tpb2"}="cpb2" **@{}; ? **@{-};
"tpb2"; "cpb2" **@{}; ? **@{-}} \]
in $\ca K$, and as before, every $f:A{\rightarrow}B$ in $\PSH {\G}$
is $\tau'$-admissible as a map in $\ca K$.
By theorem(\ref{prop:univ-enlarge}) $f$ is $\Sigma\tau$-admissible
iff it is $\tau$-admissible and the adjoint transpose of $B(f,1)$
factors through $i_1$. This last condition says that
$\forall a \in A_0, b \in B_0$, the set $B_0(fa,b)$ is a singleton.
In particular, $A \in \PSH {\G}$ is $\Sigma\tau$-admissible
when $A_0$ is chaotic{\footnotemark{\footnotetext{
Recall that a category $X$ is \emph{chaotic} when
$\forall x,y \in X$ there is a unique arrow $x{\rightarrow}y$ in $X$.
Thus a chaotic category is either empty or equivalent to $1$.}}}
and the $A_n$ are locally small for $n>0$.
Using the left 2-adjoint $D$ of $\Sigma$ described
in example(\ref{cdc-globcatII}),
we can say these things a little more efficiently.
A map $f:A{\rightarrow}B$ of $\PSH {\G}$ is $\Sigma\tau$-admissible
iff $Df$ is $\tau$-admissible and
$\forall a \in A_0, b \in B_0$, the set $B_0(fa,b)$ is a singleton.
An object $A$ of $\PSH {\G}$ is $\Sigma\tau$-admissible
iff $DA$ is $\tau$-admissible and $A_0$ is chaotic.
By the adjointness $D \ladj \Sigma$, cartesian closedness,
and since $D$ preserves products and commutes with the duality involution
we have isomorphisms
\[ \PSH {\G}(X,[\iop A,\Sigma\Omega])
\iso \PSH {\G}(DX \times D{\iop A},\Omega)
\iso \PSH {\G}(X,\Sigma\PSh {DA}) \]
natural in $X$, and so by the 2-categorical yoneda lemma, we get
an isomorphism $[\iop A,\Sigma\Omega] \iso \Sigma\PSh {DA}$.
Thus, $A \in \PSH {\G}$ is $\Sigma\tau$-small
iff $DA$ and $\PSh {DA}$ are $\tau$-admissible and $A_0$
is chaotic, that is, iff $DA$ is $\tau$-small and $A_0$ is
chaotic.

Let $k \in \N$.
Then replacing $\Sigma$ by $\Sigma^k$, 
and thus $D$ by the left adjoint $D^k$ of $\Sigma^k$
in the above analysis,
one can characterise the $\Sigma^k\tau$-admissible maps and objects,
and the $\Sigma^k\tau$-small objects in exactly the same way.
The results of this are:
\begin{enumerate}
\item  $f:A{\rightarrow}B$ is $\Sigma^k\tau$-admissible iff
for $n<k$: $\forall a \in A_n, b \in B_n$ there is a unique
arrow $fa{\rightarrow}b$ in $B_n$;
and for $n{\geq}k$: $\forall a \in A_n, b \in B_n$,
$B_n(fa,b)$ is in $\Set$.
\item  $A$ is $\Sigma^k\tau$-admissible iff
for $n<k$ the $A_n$ are chaotic
and for $n{\geq}k$ the $A_n$ are locally small.
\item  $A$ is $\Sigma^k\tau$-small iff
for $n<k$ the $A_n$ are chaotic
and for $n{\geq}k$ the $A_n$ are small.
\end{enumerate}
Analogously to the case $k=1$ we denote the map
$\Sigma^k\Omega(\Sigma^ky_1,1)$ by $i_k$.
The reader may verify, by unpacking the necessary definitions,
that $i_k$ has the following explicit description.
For $n{\geq}k$ an object of $\Sigma^k\Omega_n$ is an $(n-k)$-span
of small sets. The effect of $i_k$ on such an $(n-k)$-span $X$
is the $n$-span obtained by putting $1$ (a terminal object of $\Set$)
in the bottom $k$-levels
and $X$ in the top $(n-k)$-levels.
For example for $k=2$ and $n=4$:
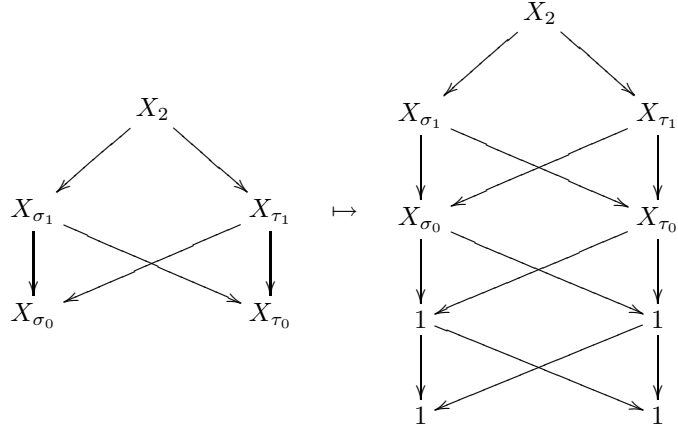
\[ \TwoDiagRel
{\xymatrix{& {X_2} \ar[dr] \ar[dl] \\
{X_{\sigma_1}} \ar[d] \ar[drr] && {X_{\tau_1}} \ar[d] \ar[dll] \\
{X_{\sigma_0}} && {X_{\tau_0}}}}
{\mapsto}
{\xymatrix{& {X_2} \ar[dr] \ar[dl] \\
{X_{\sigma_1}} \ar[d] \ar[drr] && {X_{\tau_1}} \ar[d] \ar[dll] \\
{X_{\sigma_0}} \ar[d] \ar[drr] && {X_{\tau_0}} \ar[d] \ar[dll] \\
{1} \ar[d] \ar[drr] && {1} \ar[d] \ar[dll] \\
{1} && {1}}} \]
and so $i_k$ is an ``internalisation'' of $\Sigma^k$.
In fact $\PSH {\G}(1,i_k)$ is exactly the endofunctor of $\PSh {\G}$
given by right kan extension along $(-+1)$
mentioned in example(\ref{cdc-globcatII}).
\end{example}

\section{Cocomplete presheaves for 2-toposes}
\label{sec:psh-coco}

This section is devoted to explaining why for our main
examples of 2-toposes, the resulting yoneda structures
do indeed have cocomplete presheaves.
The main result of this section,
theorem(\ref{pres-adm}),
relies on there being available
some dense sub-2-category of our 2-topos consisting of small objects.
For our examples this sub-2-category will be given by the representables.
\begin{lemma}\label{dense-all}
Let $\ca D \subseteq \ca K$ be a dense inclusion and
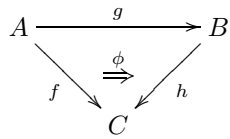
\[ \TriTwoCell A B C f g h {\phi} \]
be in $\ca K$.
If for all $D \in \ca D$ and $a:D{\rightarrow}A$,
${\phi}a$ exhibits $ga$ as a left lifting of $fa$ along $h$;
then ${\phi}$ exhibits $g$ as an absolute left lifting of $f$ along $h$.
\end{lemma}
\begin{proof}
Let $a:X{\rightarrow}A$, $b:X{\rightarrow}B$
and $\gamma:fa{\rightarrow}hb$.
We must exhibit unique ${\gamma}':ga{\rightarrow}b$
such that
\begin{equation}\label{wlem-eq1}
\TwoDiagRel
{\xymatrix{& X \ar[dl]_{a} \save \POS?="domgamp" \restore
\ar[dr]^{b} \save \POS?="codgamp" \restore \\
A \ar[rr]|-{g} \ar[dr]_{f} \save \POS?="domphi" \restore
&& B \ar[dl]^{h} \save \POS?="codphi" \restore \\ & C
\POS "domgamp"; "codgamp" **@{}; ?(.4) \ar@{=>}_{\gamma'} ?(.6)
\POS "domphi"; "codphi" **@{}; ?(.4) \ar@{=>}^{\phi} ?(.6)}}
{=} {*{\gamma}}
\end{equation}
For all $D \in \ca D$ and $x:D{\rightarrow}X$ there is by hypothesis
a unique $\gamma'_x$ such that
\begin{equation}\label{wlem-eq2}
\TwoDiagRel
{\xymatrix{& D \ar[dl]_{ax} \save \POS?="domgamp" \restore
\ar[dr]^{bx} \save \POS?="codgamp" \restore \\
A \ar[rr]|-{g} \ar[dr]_{f} \save \POS?="domphi" \restore
&& B \ar[dl]^{h} \save \POS?="codphi" \restore \\ & C
\POS "domgamp"; "codgamp" **@{}; ?(.4) \ar@{=>}_{\gamma'_x} ?(.6)
\POS "domphi"; "codphi" **@{}; ?(.4) \ar@{=>}^{\phi} ?(.6)}}
{=} {*{\gamma{x}}}
\end{equation}
and by uniqueness these $\gamma'_x$ are natural in $x$ and $D$:
in other words
\[ \TwoDiagRel
{\xymatrix{{gax_1} \ar[r]^-{\gamma'_{x_1}} \ar[d]_{ga\psi}
\save \POS?="domeq" \restore
& {bx_1} \ar[d]^{b\psi} \save \POS?="codeq" \restore
\\ {gax_1} \ar[r]_-{\gamma'_{x_2}} & {bx_1}
\POS "domeq"; "codeq" **@{}; ?*{=}}}
{}
{*{\gamma'_xd=\gamma'_{xd}}} \]
for all $\psi:x_1{\rightarrow}x_2$ and $d:D'{\rightarrow}D$.
Thus the $\gamma_x$ are in fact 2-natural in $D$.
By the density of $\ca D$ there is a unique ${\gamma}':ga{\rightarrow}b$
such that $\gamma'_x=\gamma'x$ and this $\gamma'$ clearly
satisfies (\ref{wlem-eq1}).
To see that $\gamma'$ is unique
let $\gamma'':ga{\rightarrow}b$
satisfies equation(\ref{wlem-eq1})
(with $\gamma''$ in place of $\gamma'$).
Then $\gamma''x$ satisfies the appropriate analogue of
equation(\ref{wlem-eq2})
for all $D \in \ca D$ and $x:D{\rightarrow}X$,
whence $\gamma''x=\gamma'x$, so that by the density of $\ca D$,
$\gamma''=\gamma'$.
\end{proof}
\begin{lemma}\label{dense-weights}
Let $\ca K$ be a finitely complete 2-category equipped with a
good yoneda structure and $\ca D \subseteq \ca K$ be a dense inclusion.
Suppose that
\[ \begin{array}{lccr}
{i:M{\rightarrow}\PSh C} &&& {f:C{\rightarrow}A}
\end{array} \]
in $\ca K$ where $C$ is small and $M$ and $f$ are admissible,
and that
\[ \TriTwoCell  M A {\PSh C} i k {A(f,1)} {\eta} \]
The following statements are equivalent:
\begin{enumerate}
\item  $\forall D \in \ca D$, $\forall m:D{\rightarrow}M$,
$\eta{m}$ exhibits $km$ as a left lifting of $im$ along $A(f,1)$.
\label{wts1}
\item  $\forall D \in \ca D$, $\forall m:D{\rightarrow}M$,
$\eta{m}$ exhibits $km$ as $\col(im,f)$.
\label{wts2}
\item  $\eta$ exhibits $k$ as $\col(i,f)$.
\label{wts3}
\end{enumerate}
\end{lemma}
\begin{proof}
Trivially (\ref{wts3})$\implies$(\ref{wts2})$\implies$(\ref{wts1})
and (\ref{wts1})$\implies$(\ref{wts3}) by lemma(\ref{dense-all}).
\end{proof}
\begin{theorem}\label{pres-adm}
Let
\[ \xymatrix{{\ca A} \ar@<-1.5ex>[rr]_-{S} \save \POS?="bot" \restore
&& {\ca B} \ar@<-1.5ex>[ll]_-{E} \save \POS?="top" \restore
\POS "top"; "bot" **@{}; ?(.5)*{\perp}} \]
be a $2$-adjunction, $\ca A$ and $\ca B$ be
finitely complete $2$-categories
each equipped with a good yoneda structure.
Suppose
\begin{enumerate}
\item  there is a dense inclusion $\ca D \subseteq \ca B$
such that all $D \in \ca D$ are admissible,
\label{hyp-1}
\item  for $\xymatrix@1{A \ar[r]^-{f} & B & C \ar[l]_-{g}}$ in $\ca A$,
$f/g$ is small when $A$ and $C$ are small and $B$ is admissible,
\label{hyp-5}
\item  $E$ preserves small and admissible objects,
pullbacks and opfibrations, and
\label{hyp-3}
\item  $S$ preserves admissible objects.
\label{hyp-4}
\end{enumerate}
If $A \in \ca A$ is admissible and cocomplete
then $SA$ is admissible and cocomplete.
\end{theorem}
\begin{proof}
Let $A \in \ca A$ be admissible and cocomplete.
Given
\[ \begin{array}{lccr}
i:M{\rightarrow}{\PSh C} &&& f:C{\rightarrow}SA
\end{array} \]
with $C$ small and $M$ admissible
we must produce $\eta:i{\rightarrow}SA(f,k)$
which exhibits $k$ as an absolute lifting of $i$ along $SA(f,1)$.
By corollary(\ref{weighted-as-conical}) such $\eta$ are in bijection with
\[ \LaxSq {y_C/i} M C {SA} p q k f {\phi} \]
exhibiting $k$ as a pointwise left extension of $fp$ along $q$
where $p$ and $q$ come from the defining lax pullback square for $y_C/i$.
By corollaries (\ref{weighted-as-conical}) and (\ref{dense-weights}),
and theorem(\ref{cofib-lex}), to say that $\phi$ is
such a pointwise left extension is to say that
$\forall D \in \ca D$ and $m:D{\rightarrow}M$, ${\phi}p_m$
exhibits $km$ as a left extension of $fpp_m$ along $q_m$
where
\begin{equation}\label{top-pb}
\PbSq {y_C/im} D M {y_C/i} {p_m} {q_m} m q
\end{equation}
A 2-cell $\phi$ corresponds by $E \ladj S$ to a 2-cell
\[ \LaxSq {E(y_C/i)} {EM} {EC} {A} {Ep} {Eq}
{\tilde{k}} {\tilde{f}} {\tilde{\phi}} \]
called the mate of $\phi$.
Notice that $EM$ is admissible since $E$ preserves admissible objects. 
By hypothesis (\ref{hyp-5}) and since $E$ preserves small objects,
$E(y_C/i)$ is small, and since $A$ is admissible cocomplete
there is a $\tilde{\phi}$ exhibiting $\tilde{k}$
as a pointwise left extension of $\tilde{f}Ep$ along $Eq$.
To see that the corresponding $\phi$ is a pointwise left extension
note that $\forall D \in \ca D$ and $m:D{\rightarrow}M$,
the mate of ${\phi}p_m$ is the composite obtained by applying
$E$ to (\ref{top-pb}) and mounting this on top of $\tilde{\phi}$.
Since $E$ preserves pullbacks and opfibrations this composite exhibits
$\tilde{\phi}Em$ as a left extension by theorem(\ref{cofib-lex}),
and so ${\phi}p_m$ exhibits $km$ as a left extension as required.
\end{proof}
\begin{example}\label{coco-psh-CATPSH}
Let $\C$ be a small category and consider the 2-topos
structure on $\PSH {\C}$ we began to analyze in example(\ref{CATPSH-yon}),
where $\Omega = \Sp {\C}(\Set)$.
We will now show that the resulting yoneda structure
has cocomplete presheaves.

For a small object $A \in \PSH {\C}$ the 2-adjunction
\[ \xymatrix{{\PSH {\C}} \ar@<-1.5ex>[rr]_-{[A,-]}
\save \POS?="bot" \restore
&& {\PSH {\C}} \ar@<-1.5ex>[ll]_-{(-) \times A}
\save \POS?="top" \restore
\POS "top"; "bot" **@{}; ?(.5)*{\perp}} \]
satisfies the hypotheses of theorem(\ref{pres-adm}) where
the objects of $\ca D$ are the representables.
All the hypotheses of theorem(\ref{pres-adm}) concerned with admissibility
or smallness follow easily from the 
characterisation of small and admissible objects
given in example(\ref{CATPSH-yon}).
In general, that is without any size hypothesis on $A$,
$(-){\times}A$ preserves pullbacks and opfibrations.
To see this first note that $(-){\times}A$ factors as
\[ \xymatrix{{\PSH {\C}} \ar[r]^-{t_A^*} & {\PSH {\C}/A}
\ar[r] & {\PSH {\C}}} \]
where $t_A^*$ is given by pulling back along $t_A$,
and the second leg of the factorisation takes the domain
of a map into $A$.
Thus it preserves pullbacks \cite{CJ}.
Given an opfibration $p:E{\rightarrow}B$, note that $p{\times}A$
is obtained by pulling back $p$ along
the projection $B{\times}A{\rightarrow}B$,
and so $p{\times}A$ is an opfibration.

Thus if $A$ is small and $X$ is admissible and cocomplete,
then $[A,X]$ is admissible and cocomplete.

In particular to see that $\PSH {\C}$ has cocomplete presheaves
it suffices to exhibit $\Omega$ as cocomplete.
We now consider the 2-adjunction
\[ \xymatrix{{\CAT} \ar@<-1.5ex>[rr]_-{\Sp {\C}}
\save \POS?(.46)="bot" \restore
&& {\PSH {\C}} \ar@<-1.5ex>[ll]_-{E}
\save \POS?(.54)="top" \restore
\POS "top"; "bot" **@{}; ?(.5)*{\perp}} \]
discussed in example(\ref{cdc-catpshf}),
and regard $\CAT$ as a 2-topos with classifying
discrete opfibration $U:\Set_{\bullet}{\rightarrow}\Set$.
By the characterisation of admissible and small objects
given in examples(\ref{CAT-yonII}) and (\ref{CATPSH-yon}),
this adjunction satisfies the conditions of theorem(\ref{pres-adm}).
Thus if $\ca E \in \CAT$ is locally small and cocomplete,
then $\Sp {\C}(\ca E)$ is admissible and cocomplete.

In the case $\C=\G$ a corollary of presheaf cocompleteness
and theorem(\ref{free-cocomp})
is that the globular category of higher spans of sets is the small
globular colimit completion of $1$.
\end{example}
\begin{example}\label{coco-psh-suspGLOB}
We now continue the analysis of example(\ref{susp-GLOB-yon})
and adopt the notation used there.
Moreover we adopt the following terminology:
an object $A \in \PSH {\G}$ is \emph{$\Sigma^k\tau$-cocomplete}
when it is cocomplete for the yoneda structure arising from
the 2-topos $(\PSH {\G},\iop {(-)},\Sigma^k\tau)$
by theorem(\ref{thm:2topos->yon}).
We shall now see that this yoneda structure also has cocomplete presheaves.

The hypotheses of theorem(\ref{pres-adm})
are trivially satisfied for the 2-adjunction $D \ladj \Sigma$,
where the yoneda structure for both $\ca A$ and $\ca B$
is that induced by $\tau$.
Thus if $A \in \PSH {\G}$ is $\tau$-admissible and $\tau$-cocomplete
then so is $\Sigma^kA$ for $k \in \N$.
In particular $\Sigma^k\Omega$ is $\tau$-cocomplete,
and so is $[\iop A,\Sigma^k\Omega]$ when $A$ is $\tau$-small
by the previous example.
From the characterisation of smalls
and admissibles given in example(\ref{susp-GLOB-yon}),
$\tau$-admissible implies $\Sigma^k\tau$-admissible
and $\tau$-cocomplete implies $\Sigma^k\tau$-cocomplete.
\end{example}

\section{The cartesian closedness of $\Omega$}\label{sec:omcc}

One of the facets of 2-topos theory is that $\Omega$,
the codomain of a classifying discrete opfibration
inherits quite a bit of natural structure.
A full discussion of this involves describing how the
theory of cartesian 2-monads interacts with 2-topos theory,
and is provided in \cite{OpII}. We shall now
describe how $\Omega$ is cartesian closed.

\begin{definition}\label{def:cpm-cc}
In a finitely complete 2-category $\ca K$,
an object $A$ is a \emph{cartesian pseudo monoid}
when the maps
\[ \xymatrix{1 & A \ar[l] \ar[r]^-{\Delta} & {A{\times}A}} \]
have right adjoints.
If for each $a:1{\rightarrow}A$ the map $(-{\times}a)$,
defined as the composite
\[ \xymatrix{{A{\iso}A{\times}1} \ar[r]^-{1_A{\times}a} &
{A{\times}A} \ar[r]^-{m} & A} \]
where $\Delta \ladj m$, has a right adjoint,
then $A$ is said to be \emph{cartesian closed}.
\end{definition}

In $\CAT$ a cartesian pseudo monoid is a category with
finite products, and an object of $\CAT$ is cartesian closed
in the sense of definition(\ref{def:cpm-cc})
when it is a cartesian closed category in the usual sense.
We now fix a 2-topos $(\ca K,\iop {(-)},\tau)$
and under some hypotheses clearly satisfied by
our main examples{\footnotemark{\footnotetext{
In particular examples(\ref{cdc-cat}),
(\ref{cdc-catpshf}), and (\ref{cdc-globcatII})
and the analogues of these obtained by replacing
$\Set$ everywhere by $\Set_{\lambda}$
for any regular cardinal $\lambda$.}}},
exhibit $\Omega$ as a cartesian closed object of $\ca K$.

\begin{proposition}
If $1$ is admissible then $y_1:1{\rightarrow}\Omega$
is right adjoint to $t_{\Omega}$.
\end{proposition}
\begin{proof}
Define $\eta$ by
\[ \TwoDiagRel
{\xymatrix{1 \ar[dr]_{y_1} \save \POS?="domid" \restore \ar[r]^-{y_1}
& {\Omega} \ar[d]|{1} \save \POS?="codid" \restore
\save \POS?="dometa" \restore \ar[r]^-{t_{\Omega}} & 1 \ar[dl]^{y_1}
\save \POS?="codeta" \restore \\ & {\Omega}
\POS "domid"; "codid" **@{} ?(.35) \ar@{=>}^{\id} ?(.65)
\POS "dometa"; "codeta" **@{} ?(.35) \ar@{=>}^{\eta} ?(.65)}}
{=} {*{\id}} \]
so that $\eta$ exhibits $y_1$ as a left extension along $t_{\Omega}$.
Note that this left extension is preserved by $t_{\Omega}$
and so $t_{\Omega} \ladj y_1$.
\end{proof}
\begin{example}\label{one-not-admissible}
We now exhibit an example of a 2-topos for which $1$ is not
admissible. 
Consider the classifying discrete opfibration obtained from
$U:\Set_{\bullet}{\rightarrow}\Set$ by pulling it back along
the inclusion of the full subcategory of $\Set$
consisting of those sets of cardinality at least $2$.
For the resulting yoneda structure on $\CAT$, $1$ is not
admissible since this full subcategory of $\Set$ lacks a terminal object.
\end{example}
The following result exhibits the multiplication of a cartesian
pseudo monoid structure on $\Omega$ in 2-topos theoretic terms.
A discrete fibration in $\ca K$ has \emph{small fibres}
when it is an attribute. Given a discrete fibration
$(d_1,E_1,c_1)$ from $A_1$ to $B_1$ and a discrete fibration
$(d_2,E_2,c_2)$ from $A_2$ to $B_2$, then
$(d_1{\times}d_2,E_1{\times}E_2,c_1{\times}c_2)$
is a discrete fibration from $A_1{\times}A_2$ to $B_1{\times}B_2$:
this is easily seen to be true in $\CAT$ by direct inspection,
and then in general by the representability of the notions involved.
Thus in particular $\tau \times \tau$ is a discrete opfibration.
\begin{proposition}\label{prop:om-cpm}
Suppose that $\tau \times \tau$ has small fibres.
Then a map
\[ m : \Omega \times \Omega \rightarrow \Omega \]
which classifies $\tau \times \tau$
is right adjoint to the diagonal map $\Delta_{\Omega}$.
\end{proposition}
\begin{proof}
Note that $\Omega_{\bullet}$
may be identified with the head of the span $y_1/{\Omega}$
by proposition(\ref{att-yoneda}).
For $x,z_1,z_2:X{\rightarrow}{\Omega}$
we must exhibit a bijection between 2-cells
$x{\rightarrow}m(z_1,z_2)$ and $x\Delta_{\Omega}{\rightarrow}(z_1,z_2)$
which is natural in $X$ and in $(z_1,z_2)$: then the 2-cell
$\eta:1{\rightarrow}m\Delta_{\Omega}$ corresponding
to the identity exhibits $\Delta_{\Omega}$ as an absolute
left lifting of $1_{\Omega}$ along $m$, and so is the unit
of an adjunction.
Maps $x{\rightarrow}m(z_1,z_2)$
are in bijection with maps of spans $y_1/x{\rightarrow}y_1/m(z_1,z_2)$
by the universal property of $\tau$.
Now $m(z_1,z_2)=m(z_1{\times}z_2)\Delta_X$ so that
$y_1/m(z_1,z_2){\iso}{\Delta_X^{\textnormal{rev}}}
{\comp}y_1/m(z_1{\times}z_2)$,
and so by the adjointness $\Delta_X{\ladj}{\Delta_X^{\textnormal{rev}}}$,
span maps $y_1/x{\rightarrow}y_1/m(z_1,z_2)$ are in bijection
with span maps $\Delta_X{\comp}y_1/x{\rightarrow}y_1/m(z_1{\times}z_2)$.
From the pullbacks
\[ \xymatrix{{{y_1/z_1}{\times}y_1/z_2}
\ar[d] \save \POS?(.3)="lpb2" \restore
\ar[r] \save \POS?(.3)="tpb2" \restore
& {{y_1/{\Omega}}{\times}{y_1/{\Omega}}}
\ar[d]^{\tau{\times}\tau} \save \POS?(.3)="rpb2" \restore
\save \POS?(.3)="lpb3" \restore \ar[r] \save \POS?(.3)="tpb3" \restore
& {{y_1/{\Omega}}} \ar[d]^{\tau} \save \POS?(.3)="rpb3" \restore \\
{X{\times}X} \ar[r]_-{z_1{\times}z_2} \save \POS?(.3)="bpb2" \restore
& {\Omega{\times}\Omega} \ar[r]_-{m}
\save \POS?(.3)="bpb3" \restore & {\Omega}
\POS "rpb2"; "lpb2" **@{}; ?!{"bpb2";"tpb2"}="cpb2" **@{}; ? **@{-};
"tpb2"; "cpb2" **@{}; ? **@{-}
\POS "rpb3"; "lpb3" **@{}; ?!{"bpb3";"tpb3"}="cpb3" **@{}; ? **@{-};
"tpb3"; "cpb3" **@{}; ? **@{-}} \]
one has $y_1/m(z_1{\times}z_2){\iso}y_1/z_1{\times}y_1/z_2$,
and so span maps $\Delta_X{\comp}y_1/x{\rightarrow}y_1/m(z_1{\times}z_2)$
are in bijection with
span maps $\Delta_X{\comp}y_1/x{\rightarrow}y_1/z_1{\times}y_1/z_2$,
which are in bijection with pairs of span maps
$(y_1/x{\rightarrow}y_1/z_1,y_1/x{\rightarrow}y_1/z_2)$,
and these are in bijection with pairs of 2-cells
$(x{\rightarrow}z_1,x{\rightarrow}z_2)$
by the universal property of $\tau$.
Finally these pairs of 2-cells are in bijection with 2-cells
$x\Delta_{\Omega}{\rightarrow}(z_1,z_2)$.
It is straightforward to verify that
these bijections are easily are respected by composition
with maps into $X$, and by composition with 2-cells out of $(z_1,z_2)$.
\end{proof}
In the main examples the hypotheses of these results
are easy to verify, but on the other hand, the fact
that $\Omega$ is a cartesian pseudo monoid
is easily verifiable without them: cartesian pseudo monoids
are preserved by finite product preserving 2-functors,
and our main examples are obtained from well-known examples in $\CAT$
by applying certain right 2-adjoints.

However cartesian closed objects are not
preserved by right 2-adjoints{\footnotemark{\footnotetext{
For an example note that the functor category
$[\Sigma{\N},\Set_f]$ is not cartesian closed
where $\Sigma\N$ is the monoid of natural numbers under addition
regarded as a one object category, and $\Set_f$ is the
category of finite sets. Thus the representable
2-functor $\CAT(\Sigma{\N},-)$ doesn't preserve cartesian
closed objects.}}}
and so proposition(\ref{prop:om-cpm})
is useful to us because it specifies \emph{how} the
cartesian pseudo monoid structure of $\Omega$ is specified
in 2-topos theoretic terms. We exploit this in the proof
of theorem(\ref{thm:om-cc}). First we require a lemma,
which is an analogue of the yoneda lemma for 2-sided discrete
fibrations (theorem(\ref{yon})).
\begin{lemma}\label{an-of-yon}
Let $\ca K$ be a finitely complete 2-category and $f:A{\rightarrow}B$
be in $\ca K$.
Recall the definition of the map $i_f:A{\rightarrow}f/B$
prior to theorem(\ref{yon}).
For any span $(d_1,E_1,c_1)$ from $X$ to $Z$
and discrete fibration $(d_2,E_2,c_2)$ from $A{\times}X$ to $B{\times}Z$,
composition with
\[ i_f \times \id : f \times E_1 \rightarrow f/B \times E_1 \]
determines a bijection between maps of spans
$f/B{\times}E_1{\rightarrow}E_2$ and maps of spans
$f{\times}E_1{\rightarrow}E_2$.
\end{lemma}
\begin{proof}
By the representability of the notions involved
it suffices to prove this result for the case $\ca K=\CAT$.
In this case
the head of the span $f/B \times E_1$ can be described as follows:
\begin{itemize}
\item  objects are 4-tuples
$(e,a,\beta:fa{\rightarrow}b,b)$
where $e \in E_1$.
\item  an arrow
\[ (\varepsilon,\alpha,\beta) :
(e_1,a_1,\beta_1,b_1) \rightarrow (e_2,a_2,\beta_2,b_2) \]
consists of maps $\varepsilon:e_1{\rightarrow}e_2$,
$\alpha:a_1{\rightarrow}a_2$
and $\beta:b_1{\rightarrow}b_2$,
such that $\beta\beta_1=\beta_2f(\alpha)$.
\end{itemize}
and the left and right legs of the span send $(\varepsilon,\alpha,\beta)$
described above to $(\alpha,d_1\varepsilon)$ and $(\beta,c_1\varepsilon)$
respectively. The image of $i_f \times \id$ is the full subcategory
given by the $(e,a,\beta,b)$
such that $\beta=\id$.
Let $\phi:f/B{\times}E_1{\rightarrow}E_2$.
For any object $(e,a,\beta,b)$ we have a map
\[ (1,1,\beta) : (e,a,\id,fa) \rightarrow (e,a,\beta,b) \]
and since $d_2\phi(1,1,\beta)=\id$ and $c_2\phi(1,1,\beta)=(\beta,\id)$,
$\phi(e,a,\beta,b)$ and $\phi(1,1,\beta)$
are defined uniquely by $(\phi(e,a,1_{fa},fa),\beta)$
since $E_2$ is a discrete fibration.
For any arrow $(\varepsilon,\alpha,\beta)$ as above, we have
a commutative square
\[ \xymatrix{{\phi(e_1,a_1,\id,fa_1)} \ar[r]^-{\phi(1,1,\beta_1)}
\ar[d]_{\phi(\varepsilon,\alpha,f\alpha)}
& {\phi(e_1,a_1,\beta_1,b_1)} \ar[d]^{\phi(\varepsilon,\alpha,\beta)} \\
{\phi(e_2,a_2,\id,fa_2)} \ar[r]_-{\phi(1,1,\beta_2)}
& {\phi(e_2,a_2,\beta_2,b_2)}} \]
in $E_2$, but from the proof of theorem(\ref{thm-2-sided-fib})
we know that $\phi(1,1,\beta_1)$ is $d_2$-opcartesian
and so the rest of the above square is determined uniquely by
$\phi(\varepsilon,\alpha,f\alpha)$.
\end{proof}
\begin{theorem}\label{thm:om-cc}
If $1$ is small, $\Omega$ is cocomplete
and $\tau \times \tau$ has small fibres,
then $\Omega$ is cartesian closed.
\end{theorem}
\begin{proof}
For $x:1{\rightarrow}\Omega$ it suffices to provide
\[ \TriTwoCell 1 {\Omega} {\Omega} x {y_1} {-{\times}x} \psi \]
which exhibits $-{\times}x$ as a left extension of $x$ along $y_1$:
since $\Omega$ is cocomplete, the left extension of $x$ along $y_1$
exists and is left adjoint to $\Omega(x,1)$
by theorem(\ref{free-cocomp}).
We do this by exhibiting a bijection between 2-cells
$x{\rightarrow}py_1$ and 2-cells $(-{\times}x){\rightarrow}p$
which is natural in $p$.
By the universal property of $\tau$ 2-cells
$x{\rightarrow}py_1$ are in bijection with span maps
$y_1/x{\rightarrow}y_1/py_1$.
Since $y_1/py_1{\iso}y_1^{\textnormal{rev}}{\comp}y_1/p$ and
$y_1{\ladj}y_1^{\textnormal{rev}}$, these span maps are in bijection
with span maps $y_1{\comp}y_1/x{\rightarrow}y_1/p$.
We have span isomorphisms
\[ y_1 \comp y_1/x \iso (y_1{\times}\id) \comp (\id{\times}y_1/x)
\iso y_1 \times y_1/x \]
and so span maps $y_1{\comp}y_1/x{\rightarrow}y_1/p$
are in bijection with span maps $y_1{\times}y_1/x{\rightarrow}y_1/p$
which by lemma(\ref{an-of-yon}) are in bijection with span maps
$y_1/{\Omega}{\times}y_1/x{\rightarrow}y_1/p$.
From the pullbacks
\[ \xymatrix{{{y_1/{\Omega}}{\times}y_1/x}
\ar[d] \save \POS?(.3)="lpb2" \restore
\ar[r] \save \POS?(.3)="tpb2" \restore
& {{y_1/{\Omega}}{\times}{y_1/{\Omega}}}
\ar[d]^{\tau{\times}\tau} \save \POS?(.3)="rpb2" \restore
\save \POS?(.3)="lpb3" \restore \ar[r] \save \POS?(.3)="tpb3" \restore
& {{y_1/{\Omega}}} \ar[d]^{\tau} \save \POS?(.3)="rpb3" \restore \\
{\Omega{\times}1} \ar[r]_-{\id{\times}x} \save \POS?(.3)="bpb2" \restore
& {\Omega{\times}\Omega} \ar[r]_-{m}
\save \POS?(.3)="bpb3" \restore & {\Omega}
\POS "rpb2"; "lpb2" **@{}; ?!{"bpb2";"tpb2"}="cpb2" **@{}; ? **@{-};
"tpb2"; "cpb2" **@{}; ? **@{-}
\POS "rpb3"; "lpb3" **@{}; ?!{"bpb3";"tpb3"}="cpb3" **@{}; ? **@{-};
"tpb3"; "cpb3" **@{}; ? **@{-}} \]
one has $y_1/{\Omega}{\times}y_1/x{\iso}y_1/m(\id{\times}x)$,
and so span maps
$y_1/{\Omega}{\times}y_1/x{\rightarrow}y_1/p$ are in bijection with
span maps $y_1/m(\id{\times}x){\rightarrow}y_1/p$,
and these are in bijection with 2-cells $(-{\times}x){\rightarrow}p$
by the universal property of $\tau$.
It is straight forward to verify that the bijections just described
are natural in $p$.
\end{proof}

\section{Acknowledgements}

This work was completed while the author was a postdoctoral
fellow at the University of Ottawa and at UQAM,
and I am indebted to these mathematics departments,
Phil Scott and Andr{\'e} Joyal for their support.
I would also like to thank Ross Street and Susan Niefield
for comments and corrections to an earlier version of this paper.

\bibliography{highercats}

\providecommand{\bysame}{\leavevmode\hbox to3em{\hrulefill}\thinspace}
\providecommand{\MR}{\relax\ifhmode\unskip\space\fi MR }
\providecommand{\MRhref}[2]{%
  \href{http://www.ams.org/mathscinet-getitem?mr=#1}{#2}
}
\providecommand{\href}[2]{#2}
\begin{thebibliography}{Bat98b}

\bibitem[B\'85]{Ben85}
J.~B\'{e}nabou, \emph{Fibred categories and the foundation of naive category
  theory}, Journal of Symbolic Logic \textbf{50} (1985), 10--37.

\bibitem[Bat98a]{Bat98b}
M.~Batanin, \emph{Computads for finitary monads on globular sets}, Contemporary
  Mathematics \textbf{230} (1998), 37--57.

\bibitem[Bat98b]{Bat98}
\bysame, \emph{Monoidal globular categories as a natural environment for the
  theory of weak $n$-categories}, Advances in Mathematics \textbf{136} (1998),
  39--103.

\bibitem[Bat02]{Bat02}
\bysame, \emph{The {Eckmann-Hilton} argument, higher operads and
  {$E_{n}$-spaces}}, arXiv:math.CT/0207281, 2002.

\bibitem[Bat03]{Bat03}
\bysame, \emph{The combinatorics of iterated loop spaces},
  arXiv:math.CT/0301221, 2003.

\bibitem[BC04]{BC04}
J.~Baez and A.~Crans, \emph{Higher-dimensional algebra vi: Lie 2-algebras},
  Theory and applications of categories \textbf{12} (2004), 492--528.

\bibitem[BD95]{BD95}
J.~Baez and J.~Dolan, \emph{Higher-dimensional algebra and topological quantum
  field theory}, Journal Math.Phys \textbf{36} (1995), 6073--6105.

\bibitem[BD98]{BD98a}
\bysame, \emph{Higher-dimensional algebra iii: $n$-categories and the algebra
  of opetopes}, Advances in Mathematics \textbf{135} (1998), 145--206.

\bibitem[BL04]{BL04}
J.~Baez and A.~Lauda, \emph{Higher-dimensional algebra v: 2-groups}, Theory and
  applications of categories \textbf{12} (2004), 423--491.

\bibitem[Bou74]{Bou74}
D.~Bourn, \emph{Sur les ditopos}, C. R. Acad. Sci. Paris \textbf{279} (1974),
  911--913.

\bibitem[BS05]{BS05}
J.~Baez and U.~Schreiber, \emph{Higher gauge theory}, arXiv:math.DG/0511710,
  2005.

\bibitem[CJ95]{CJ}
A.~Carboni and P.T. Johnstone, \emph{Connected limits, familial
  representability and {Artin} glueing}, Mathematical Structures in Computer
  Science \textbf{5} (1995), 441--459.

\bibitem[Ehr58]{Ehr58}
C.~Ehresmann, \emph{Gattungen von lokalen strukturen}, Jber. Deutsch. Math.
  Verein \textbf{60} (1958), 49--77.

\bibitem[FS95]{FryStr}
P.~Freyd and R.~Street, \emph{On the size of categories}, Theory and
  applications of categories \textbf{1} (1995), 174--181.

\bibitem[Gra66]{Gray66}
J.W. Gray, \emph{Fibred and cofibred categories}, Proc. Conference on
  Categorical Algebra at La Jolla, Springer, 1966, pp.~21--83.

\bibitem[Gro70]{Gro70}
A.~Grothendieck, \emph{Cat\'{e}gories fibr\'{e}es et descente}, Lecture Notes
  in Math. \textbf{224} (1970), 145--194.

\bibitem[Her99]{Her99}
C.~Hermida, \emph{Some properties of fib as a fibred $2$-category}, JPAA
  \textbf{134(1)} (1999), 83--109.

\bibitem[Joh02]{JohElpI}
P.T. Johnstone, \emph{Sketches of an elephant: A topos theory compendium},
  vol.~1, Oxford logic guides, no.~43, Oxford science publications, 2002.

\bibitem[Kel82]{Kel82}
G.M. Kelly, \emph{Basic concepts of enriched category theory, lms lecture note
  series}, vol.~64, Cambridge University Press, 1982.

\bibitem[KS74]{KS74}
G.M. Kelly and R.~Street, \emph{Review of the elements of 2-categories},
  Lecture Notes in Math. \textbf{420} (1974), 75--103.

\bibitem[Law70]{LawDct}
F.W. Lawvere, \emph{Equality in hyperdoctrines and comprehension schema as an
  adjoint functor}, Proceedings of the American Mathematical Society Symposium
  on Pure Mathematics XVII, 1970, pp.~1--14.

\bibitem[MM91]{MM}
S.~{Mac Lane} and I.~Moerdijk, \emph{Sheaves in geometry and logic},
  Springer-Verlag, 1991.

\bibitem[Pen74]{Pen74}
J.~Penon, \emph{Sous-cat\'{e}gories classifi\'{e}e}, C. R. Acad. Sci. Paris
  \textbf{278} (1974), 475--477.

\bibitem[Str74a]{Str74b}
R.~Street, \emph{Elementary cosmoi}, Lecture Notes in Math. \textbf{420}
  (1974), 134--180.

\bibitem[Str74b]{Str74}
\bysame, \emph{Fibrations and yoneda's lemma in a $2$-category}, Lecture Notes
  in Math. \textbf{420} (1974), 104--133.

\bibitem[Str76]{Str76}
\bysame, \emph{Limits indexed by category-valued $2$-functors}, J. Pure Appl.
  Algebra \textbf{8} (1976), 149--181.

\bibitem[Str80a]{Str80}
\bysame, \emph{Cosmoi of internal categories}, Trans. Amer. Math. Soc.
  \textbf{258} (1980), 271--318.

\bibitem[Str80b]{Str80b}
\bysame, \emph{Fibrations in bicategories}, Cahiers Topologie G\'{e}om.
  Differentielle \textbf{21} (1980), 111--160.

\bibitem[Str00]{Str00}
\bysame, \emph{The petit topos of globular sets}, J. Pure Appl. Algebra
  \textbf{154} (2000), 299--315.

\bibitem[SW78]{SW78}
R.~Street and R.F.C. Walters, \emph{Yoneda structures on $2$-categories},
  J.Algebra \textbf{50} (1978), 350--379.

\bibitem[TV05]{TV05}
B.~T\"{o}en and G.~Vezzosi, \emph{Homotopical algebraic geometry {I}: Topos
  theory}, Advances in Mathematics \textbf{193(2)} (2005), 257--372.

\bibitem[Web]{OpII}
M.~Weber, \emph{Operads within a monoidal pseudo algebra {II}}, in preparation.

\bibitem[Web05]{OpI}
\bysame, \emph{Operads within monoidal pseudo algebras}, Applied Categorical
  Structures \textbf{13} (2005).

\end{thebibliography}
\end{document}